\def\Resetstrings{%Clears all strings before processing reference listing.
    \def\present{ }\let\bgroup={\let\egroup=}%primitive TeX
    \def\Astr{}\def\astr{}\def\Atest{}\def\atest{}%
    \def\Bstr{}\def\bstr{}\def\Btest{}\def\btest{}%
    \def\Cstr{}\def\cstr{}\def\Ctest{}\def\ctest{}%
    \def\Dstr{}\def\dstr{}\def\Dtest{}\def\dtest{}%
    \def\Estr{}\def\estr{}\def\Etest{}\def\etest{}%
    \def\Fstr{}\def\fstr{}\def\Ftest{}\def\ftest{}%
    \def\Gstr{}\def\gstr{}\def\Gtest{}\def\gtest{}%
    \def\Hstr{}\def\hstr{}\def\Htest{}\def\htest{}%
    \def\Istr{}\def\istr{}\def\Itest{}\def\itest{}%
    \def\Jstr{}\def\jstr{}\def\Jtest{}\def\jtest{}%
    \def\Kstr{}\def\kstr{}\def\Ktest{}\def\ktest{}%
    \def\Lstr{}\def\lstr{}\def\Ltest{}\def\ltest{}%
    \def\Mstr{}\def\mstr{}\def\Mtest{}\def\mtest{}%
    \def\Nstr{}\def\nstr{}\def\Ntest{}\def\ntest{}%
    \def\Ostr{}\def\ostr{}\def\Otest{}\def\otest{}%
    \def\Pstr{}\def\pstr{}\def\Ptest{}\def\ptest{}%
    \def\Qstr{}\def\qstr{}\def\Qtest{}\def\qtest{}%
    \def\Rstr{}\def\rstr{}\def\Rtest{}\def\rtest{}%
    \def\Sstr{}\def\sstr{}\def\Stest{}\def\stest{}%
    \def\Tstr{}\def\tstr{}\def\Ttest{}\def\ttest{}%
    \def\Ustr{}\def\ustr{}\def\Utest{}\def\utest{}%
    \def\Vstr{}\def\vstr{}\def\Vtest{}\def\vtest{}%
    \def\Wstr{}\def\wstr{}\def\Wtest{}\def\wtest{}%
    \def\Xstr{}\def\xstr{}\def\Xtest{}\def\xtest{}%
    \def\Ystr{}\def\ystr{}\def\Ytest{}\def\ytest{}%
}
\Resetstrings

\def\Refformat{%Determines the kind of reference by the presence or
         \if\Jtest\present
             {\if\Vtest\present\journalarticleformat
                  \else\conferencereportformat\fi}
            \else\if\Btest\present\bookarticleformat
               \else\if\Rtest\present\technicalreportformat
                  \else\if\Itest\present\bookformat
                     \else\otherformat\fi\fi\fi\fi}

\def\Rpunct{%Default punctuation control strings if the punctuation
   \def\Lspace{ }%
   \def\Lperiod{ }%  .
   \def\Lcomma{ }%    ,
   \def\Lquest{ }%     ?
   \def\Lcolon{ }%   :
   \def\Lscolon{ }%   ;
   \def\Lbang{ }%      !
   \def\Lquote{ }%   '
   \def\Lqquote{ }%   "
   \def\Lrquote{ }%    `
   \def\Rspace{}%
   \def\Rperiod{.}%  .
   \def\Rcomma{,}%    ,
   \def\Rquest{?}%     ?
   \def\Rcolon{:}%   :
   \def\Rscolon{;}%   ;
   \def\Rbang{!}%      !
   \def\Rquote{'}%   '
   \def\Rqquote{"}%   "
   \def\Rrquote{`}%    `
   }

\def\Lpunct{%Default punctuation control strings if the punctuation
   \def\Lspace{}%
   \def\Lperiod{\unskip.}%  .
   \def\Lcomma{\unskip,}%    ,
   \def\Lquest{\unskip?}%     ?
   \def\Lcolon{\unskip:}%   :
   \def\Lscolon{\unskip;}%   ;
   \def\Lbang{\unskip!}%      !
   \def\Lquote{\unskip'}%   '
   \def\Lqquote{\unskip"}%   "
   \def\Lrquote{\unskip`}%    `
   \def\Rspace{\spacefactor=1000}%
   \def\Rperiod{\spacefactor=3000}%  .
   \def\Rcomma{\spacefactor=1250}%    ,
   \def\Rquest{\spacefactor=3000}%     ?
   \def\Rcolon{\spacefactor=2000}%   :
   \def\Rscolon{\spacefactor=1250}%   ;
   \def\Rbang{\spacefactor=3000}%      !
   \def\Rquote{\spacefactor=1000}%   '
   \def\Rqquote{\spacefactor=1000}%   "
   \def\Rrquote{\spacefactor=1000}%    `
   }

\def\Refstd{%Standard control strings for formatting bibliography listings.
     \def\Acomma{\unskip, }%between multiple author names
     \def\Aand{\unskip\ and }%between two author names
     \def\Aandd{\unskip\ and }%between last two of multiple author names
     \def\Ecomma{\unskip, }%between multiple editor names
     \def\Eand{\unskip\ and }%between two editor names
     \def\Eandd{\unskip\ and }%between last two of multiple author names
     \def\acomma{\unskip, }%same for authors of reviewed material
     \def\aand{\unskip\ and }%same for authors of reviewed material
     \def\aandd{\unskip\ and }%same for authors of reviewed material
     \def\ecomma{\unskip, }%same for translators
     \def\eand{\unskip\ and }%same for translators
     \def\eandd{\unskip\ and }%same for translators
     \def\Namecomma{\unskip, }%same for citations using authors' names
     \def\Nameand{\unskip\ and }%same for citations using authors' names
     \def\Nameandd{\unskip\ and }%same for citations using authors' names
     \def\Revcomma{\unskip, }%between last and first name of reversed name
     \def\Initper{.\ }%punctuation after initial
     \def\Initgap{\dimen0=\spaceskip\divide\dimen0 by 2\hskip-\dimen0}%
   }

\def\Smallcapsaand{%Smallcaps redefinition of \Aand and \Aandd for \Refstd
     \def\Aand{\unskip\bgroup{\Smallcapsfont\ AND }\egroup}%
     \def\Aandd{\unskip\bgroup{\Smallcapsfont\ AND }\egroup}%
     \def\eand{\unskip\bgroup\Smallcapsfont\ AND \egroup}%
     \def\eandd{\unskip\bgroup\Smallcapsfont\ AND \egroup}%
   }

\def\Smallcapseand{%Smallcaps redefinition of \Eand, \Eeand, etc for Refstd
     \def\Eand{\unskip\bgroup\Smallcapsfont\ AND \egroup}%
     \def\Eandd{\unskip\bgroup\Smallcapsfont\ AND \egroup}%
     \def\aand{\unskip\bgroup\Smallcapsfont\ AND \egroup}%
     \def\aandd{\unskip\bgroup\Smallcapsfont\ AND \egroup}%
   }

   \def\Citefont{}%citations
   \def\ACitefont{}%alternate citations
   \def\Authfont{}%authors
   \def\Titlefont{}%titles
   \def\Tomefont{\sl}%journals or books
   \def\Volfont{}%volume number of journal
   \def\Flagfont{}%citation flag
   \def\Reffont{\rm}%set at beginning of reference listing
   \def\Smallcapsfont{\sevenrm}%small caps
   \def\Flagstyle#1{\hangindent\parindent\indent\hbox to0pt%flag style
       {\hss[{\Flagfont#1}]\kern.5em}\ignorespaces}%        for references

\def\Underlinemark{\vrule height .7pt depth 0pt width 3pc}%for replacing
\def\Citebrackets{\Rpunct%defaults for putting citations in brackets [].
   \def\Lcitemark{\def\Cfont{\Citefont}[\bgroup\Cfont}%mark at left of citation
   \def\Rcitemark{\egroup]}%mark at right of citation
   \def\LAcitemark{\def\Cfont{\ACitefont}\bgroup\ACitefont}%
   \def\RAcitemark{\egroup}%mark at right of alternate citation
   \def\LIcitemark{\egroup}%mark at left of insertion in citation
   \def\RIcitemark{\bgroup\Cfont}%mark at right of insertion in citation
   \def\Citehyphen{\egroup--\bgroup\Cfont}%separater for string of citations
   \def\Citecomma{\egroup,\hskip0pt\bgroup\Cfont}%
   \def\Citebreak{}%mark between parts of citation (e.g. author\Citebreak date)
   }

\def\Citeparen{\Rpunct%defaults for putting citations in parenthesis ().
   \def\Lcitemark{\def\Cfont{\Citefont}(\bgroup\Cfont}%mark at left of citation
   \def\Rcitemark{\egroup)}%mark at right of citation
   \def\LAcitemark{\def\Cfont{\ACitefont}\bgroup\ACitefont}%
   \def\RAcitemark{\egroup}%mark at right of alternate citation
   \def\LIcitemark{\egroup}%mark at left of insertion in citation
   \def\RIcitemark{\bgroup\Cfont}%mark at right of insertion in citation
   \def\Citehyphen{\egroup--\bgroup\Cfont}%separater for string of citations
   \def\Citecomma{\egroup,\hskip0pt\bgroup\Cfont}%
   \def\Citebreak{}%mark between parts of citation (e.g. author\Citebreak date)
   }

\def\Citesuper{\Lpunct%defaults for making superscript citations
   \def\Lcitemark{\def\Cfont{\Citefont}\raise1ex\hbox\bgroup\bgroup\Cfont}%
   \def\Rcitemark{\egroup\egroup}%mark at right of citation
   \def\LAcitemark{\def\Cfont{\ACitefont}\bgroup\ACitefont}%
   \def\RAcitemark{\egroup}%mark at right of alternate citation
   \def\LIcitemark{\egroup\egroup}%mark at left of insertion in citation
   \def\RIcitemark{\raise1ex\hbox\bgroup\bgroup\Cfont}%
   \def\Citehyphen{\egroup--\bgroup\Cfont}%separater for string of citations
   \def\Citecomma{\egroup,\hskip0pt\bgroup%
      \Cfont}%separater for multiple citations
   \def\Citebreak{}%mark between parts of citation (e.g. author\Citebreak date)
   }

\def\Citenamedate{\Rpunct%defaults for making name-date citations
   \def\Lcitemark{%mark at left of citation--also sets internal punctuation
      \def\Citebreak{\egroup\ [\bgroup\Citefont}%separater in citation
      \def\Citecomma{\egroup]; %between multiple citations
         \bgroup\let\uchyph=1\Citefont}(\bgroup\let\uchyph=1\Citefont}%
   \def\Rcitemark{\egroup])}%mark at right of citation
   \def\LAcitemark{%mark at left of alternate citation
      \def\Citebreak{\egroup\ [\bgroup\Citefont}\def\Citecomma{\egroup], %
         \bgroup\ACitefont }\bgroup\let\uchyph=1\ACitefont}%
   \def\RAcitemark{\egroup]}%mark at right of alternate citation
  \def\Citehyphen{\egroup--\bgroup\Citefont}%separater for string of citations
   \def\LIcitemark{\egroup}%mark at left of insertion in citation
   \def\RIcitemark{\bgroup\Citefont}%mark at right of insertion in citation
   }
\def\annotations{}
\def\Flagstyle#1{\par\noindent\hbox to 30pt{[#1]\hfil}%
\hangindent=30pt\hangafter=1}

\def\ifundefined#1{\expandafter\ifx\csname#1\endcsname\relax}
\ifundefined{thesis}\def\bye{\end{document}}\fi

\def\endauthor{:\kern6pt}

\Refstd\Citebrackets%set general formats for reference list and citations
\def\Citefont{\rm}\def\Titlefont{\Reffont}\def\Volfont{\bf}\def\Tomefont{\it}%redefine some fonts

\def\journalarticleformat{\Reffont\let\uchyph=1\parindent=1.25pc\def\Comma{}%

\sfcode`\.=1000\sfcode`\?=1000\sfcode`\!=1000\sfcode`\:=1000\sfcode`\;=1000\sfcode`\,=1000%\frenchspacing
                \par\vfil\penalty-200\vfilneg%\filbreak
      \if\Ftest\present\Flagstyle\Fstr\fi%

\if\Atest\present\bgroup\Authfont\Astr\egroup\def\Comma{\unskip\endauthor}\fi%
        \if\Ttest\present\Comma\bgroup\Titlefont\Tstr\egroup\def\Comma{, }\fi%
         \if\etest\present\hskip.2em(\bgroup\estr\egroup)\def\Comma{\unskip,
}\fi%
          \if\Jtest\present\Comma\bgroup\Tomefont\Jstr\/\egroup\def\Comma{,
}\fi%

\if\Vtest\present\if\Jtest\present\hskip.2em\else\Comma\fi\bgroup\Volfont\Vstr\egroup\def\Comma{, }\fi%
            \if\Dtest\present\hskip.2em(\bgroup\Dstr\egroup)\def\Comma{, }\fi%
             \if\Ptest\present\bgroup, \Pstr\egroup\def\Comma{, }\fi%
              \if\ttest\present\Comma\bgroup\Titlefont\tstr\egroup\def\Comma{,
}\fi%

\if\jtest\present\Comma\bgroup\Tomefont\jstr\/\egroup\def\Comma{, }\fi%

\if\vtest,\present\if\jtest\present\hskip.2em\else\Comma\fi\bgroup\Volfont\vstr\egroup\def\Comma{, }\fi%
                 \if\dtest\present\hskip.2em(\bgroup\dstr\egroup)\def\Comma{,
}\fi%
                  \if\ptest\present\bgroup, \pstr\egroup\def\Comma{, }\fi%
                   \if\Gtest\present{\Comma Gov't ordering no.
}\bgroup\Gstr\egroup\def\Comma{, }\fi%
                    \if\Mtest\present\Comma MR
\#\bgroup\Mstr\egroup\def\Comma{, }\fi%
                     \if\Otest\present{\Comma\bgroup\Ostr\egroup.}\else{.}\fi%
\annotations  %  MODIFICATION
\vskip5ptplus1ptminus1pt}%\smallskip (but 3 is now 5)  MODIFICATION

\def\conferencereportformat{\Reffont\let\uchyph=1\parindent=1.25pc\def\Comma{}%

\sfcode`\.=1000\sfcode`\?=1000\sfcode`\!=1000\sfcode`\:=1000\sfcode`\;=1000\sfcode`\,=1000%\frenchspacing
                \par\vfil\penalty-200\vfilneg%\filbreak
      \if\Ftest\present\Flagstyle\Fstr\fi%

\if\Atest\present\bgroup\Authfont\Astr\egroup\def\Comma{\unskip\endauthor}\fi%
        \if\Ttest\present\Comma\bgroup\Titlefont\Tstr\egroup\def\Comma{, }\fi%
         \if\Jtest\present\Comma\bgroup\Tomefont\Jstr\/\egroup\def\Comma{,
}\fi%
          \if\Ctest\present\Comma\bgroup\Cstr\egroup\def\Comma{, }\fi%
           \if\Dtest\present\hskip.2em(\bgroup\Dstr\egroup)\def\Comma{, }\fi%
            \if\Mtest\present\Comma MR \#\bgroup\Mstr\egroup\def\Comma{, }\fi%
             \if\Otest\present{\Comma\bgroup\Ostr\egroup.}\else{.}\fi%
\annotations  %  MODIFICATION
\vskip5ptplus1ptminus1pt}%\smallskip (but 3 is now 5)  MODIFICATION

\def\bookarticleformat{\Reffont\let\uchyph=1\parindent=1.25pc\def\Comma{}%

\sfcode`\.=1000\sfcode`\?=1000\sfcode`\!=1000\sfcode`\:=1000\sfcode`\;=1000\sfcode`\,=1000%\frenchspacing
                \par\vfil\penalty-200\vfilneg%\filbreak
      \if\Ftest\present\Flagstyle\Fstr\fi%

\if\Atest\present\bgroup\Authfont\Astr\egroup\def\Comma{\unskip\endauthor}\fi%
        \if\Ttest\present\Comma\bgroup\Titlefont\Tstr\egroup\def\Comma{, }\fi%
         \if\etest\present\hskip.2em(\bgroup\estr\egroup)\def\Comma{\unskip,
}\fi%
          \if\Btest\present\Comma in
\bgroup\Tomefont\Bstr\/\egroup\def\Comma{\unskip, }\fi%
           \if\otest\present\Comma\ \bgroup\ostr\egroup\def\Comma{, }\fi%
            \if\Etest\present\Comma\bgroup\Estr\egroup\unskip,
\ifnum\Ecnt>1eds.\else ed.\fi\def\Comma{, }\fi%
             \if\Stest\present\Comma\bgroup\Sstr\egroup\def\Comma{, }\fi%
              \if\Vtest\present\ \bgroup\bf\Vstr\egroup\def\Comma{, }\fi%
               \if\Ntest\present\Comma no. \bgroup\Nstr\egroup\def\Comma{,
}\fi%
                \if\Itest\present\Comma\bgroup\Istr\egroup\def\Comma{, }\fi%
                 \if\Ctest\present\ (\bgroup\Cstr\egroup)\def\Comma{, }\fi%
                  \if\Dtest\present\Comma\bgroup\Dstr\egroup\def\Comma{, }\fi%
                   \if\Ptest\present\Comma pp.\ \Pstr\def\Comma{, }\fi%
\if\ttest\present\Comma\bgroup\Titlefont\Tstr\egroup\def\Comma{, }\fi%
                     \if\btest\present\Comma in
\bgroup\Tomefont\bstr\egroup\def\Comma{, }\fi%
                       \if\atest\present\Comma\bgroup\astr\egroup\unskip,
\if\acnt\present eds.\else ed.\fi\def\Comma{, }\fi%
                        \if\stest\present\Comma\bgroup\sstr\egroup\def\Comma{,
}\fi%
                         \if\vtest\present\Comma vol.
\bgroup\vstr\egroup\def\Comma{, }\fi%
                          \if\ntest\present\Comma no.
\bgroup\nstr\egroup\def\Comma{, }\fi%

\if\itest\present\Comma\bgroup\istr\egroup\def\Comma{, }\fi%

\if\ctest\present\Comma\bgroup\cstr\egroup\def\Comma{, }\fi%

\if\dtest\present\Comma\bgroup\dstr\egroup\def\Comma{, }\fi%
                              \if\ptest\present\Comma\pstr\def\Comma{, }\fi%
                               \if\Gtest\present{\Comma Gov't ordering no.
}\bgroup\Gstr\egroup\def\Comma{, }\fi%
                                \if\Mtest\present\Comma MR
\#\bgroup\Mstr\egroup\def\Comma{, }\fi%

\if\Otest\present{\Comma\bgroup\Ostr\egroup.}\else{.}\fi%
\annotations  %  MODIFICATION
\vskip5ptplus1ptminus1pt}%\smallskip (but 3 is now 5)  MODIFICATION

\def\bookformat{\Reffont\let\uchyph=1\parindent=1.25pc\def\Comma{}%

\sfcode`\.=1000\sfcode`\?=1000\sfcode`\!=1000\sfcode`\:=1000\sfcode`\;=1000\sfcode`\,=1000%\frenchspacing
                \par\vfil\penalty-200\vfilneg%\filbreak
      \if\Ftest\present\Flagstyle\Fstr\fi%

\if\Atest\present\bgroup\Authfont\Astr\egroup\def\Comma{\unskip\endauthor}%

\else\if\Etest\present\bgroup\def\Eand{\Aand}\def\Eandd{\Aandd}\Authfont\Estr\egroup\unskip, \ifnum\Ecnt>1eds.\else ed.\fi\def\Comma{, }%

\else\if\Itest\present\bgroup\Authfont\Istr\egroup\def\Comma{, }\fi\fi\fi%

\if\Ttest\present\Comma\bgroup\Tomefont\Tstr\/\egroup\def\Comma{\unskip, }%

\else\if\Btest\present\Comma\bgroup\Titlefont\Bstr\/\egroup\def\Comma{\unskip,
}\fi\fi%
            \if\otest\present\Comma\ \bgroup\ostr\egroup\def\Comma{, }\fi%

\if\etest\present\hskip.2em(\bgroup\estr\egroup)\def\Comma{\unskip, }\fi%
              \if\Stest\present\Comma\bgroup\Sstr\egroup\def\Comma{, }\fi%
               \if\Vtest\present\ \bgroup\bf\Vstr\egroup\def\Comma{, }\fi%
                \if\Ntest\present\Comma no. \bgroup\Nstr\egroup\def\Comma{,
}\fi%
                 \if\Atest\present\if\Itest\present
                         \Comma\bgroup\Istr\egroup\def\Comma{\unskip, }\fi%
                      \else\if\Etest\present\if\Itest\present
                              \Comma\bgroup\Istr\egroup\def\Comma{\unskip,
}\fi\fi\fi%
                     \if\Ctest\present\ (\bgroup\Cstr\egroup)\def\Comma{, }\fi%
                      \if\Dtest\present\Comma\bgroup\Dstr\egroup\def\Comma{,
}\fi%
                      \if\Ptest\present\bgroup, pp.\ \Pstr\egroup\def\Comma{,
}\fi% MODIFICATION

\if\ttest\present\Comma\bgroup\Titlefont\tstr\egroup\def\Comma{, }%

\else\if\btest\present\Comma\bgroup\Titlefont\bstr\egroup\def\Comma{, }\fi\fi%

\if\stest\present\Comma\bgroup\sstr\egroup\def\Comma{, }\fi%
                           \if\vtest\present\Comma vol.
\bgroup\vstr\egroup\def\Comma{, }\fi%
                            \if\ntest\present\Comma no.
\bgroup\nstr\egroup\def\Comma{, }\fi%

\if\itest\present\Comma\bgroup\istr\egroup\def\Comma{, }\fi%

\if\ctest\present\Comma\bgroup\cstr\egroup\def\Comma{, }\fi%

\if\dtest\present\Comma\bgroup\dstr\egroup\def\Comma{, }\fi%
                                \if\Gtest\present{\Comma Gov't ordering no.
}\bgroup\Gstr\egroup\def\Comma{, }\fi%
                                 \if\Mtest\present\Comma MR
\#\bgroup\Mstr\egroup\def\Comma{, }\fi%

\if\Otest\present{\Comma\bgroup\Ostr\egroup.}\else{.}\fi%
\annotations  %  MODIFICATION
\vskip5ptplus1ptminus1pt}%\smallskip (but 3 is now 5)  MODIFICATION

\def\technicalreportformat{\Reffont\let\uchyph=1\parindent=1.25pc\def\Comma{}%

\sfcode`\.=1000\sfcode`\?=1000\sfcode`\!=1000\sfcode`\:=1000\sfcode`\;=1000\sfcode`\,=1000%\frenchspacing
                \par\vfil\penalty-200\vfilneg%\filbreak
      \if\Ftest\present\Flagstyle\Fstr\fi%

\if\Atest\present\bgroup\Authfont\Astr\egroup\def\Comma{\unskip\endauthor}%

\else\if\Etest\present\bgroup\def\Eand{\Aand}\def\Eandd{\Aandd}\Authfont\Estr\egroup\unskip, \ifnum\Ecnt>1eds.\else ed.\fi\def\Comma{, }%

\else\if\Itest\present\bgroup\Authfont\Istr\egroup\def\Comma{, }\fi\fi\fi%
          \if\Ttest\present\Comma\bgroup\Titlefont\Tstr\egroup\def\Comma{,
}\fi%
           \if\Atest\present\if\Itest\present
                   \Comma\bgroup\Istr\egroup\def\Comma{, }\fi%
                \else\if\Etest\present\if\Itest\present
                        \Comma\bgroup\Istr\egroup\def\Comma{, }\fi\fi\fi%
            \if\Rtest\present\Comma\bgroup\Rstr\egroup\def\Comma{, }\fi%
             \if\Ctest\present\ (\bgroup\Cstr\egroup)\def\Comma{, }\fi%
              \if\Dtest\present\Comma\bgroup\Dstr\egroup\def\Comma{, }\fi%
              \if\Ptest\present\bgroup, \Pstr\egroup\def\Comma{, }\fi%
               \if\ttest\present\Comma\bgroup\Titlefont\tstr\egroup\def\Comma{,
}\fi%
                \if\itest\present\Comma\bgroup\istr\egroup\def\Comma{, }\fi%
                 \if\rtest\present\Comma\bgroup\rstr\egroup\def\Comma{, }\fi%
                  \if\ctest\present\Comma\bgroup\cstr\egroup\def\Comma{, }\fi%
                   \if\dtest\present\Comma\bgroup\dstr\egroup\def\Comma{, }\fi%
                    \if\Gtest\present{\Comma Gov't ordering no.
}\bgroup\Gstr\egroup\def\Comma{, }\fi%
                     \if\Mtest\present\Comma MR
\#\bgroup\Mstr\egroup\def\Comma{, }\fi%
                      \if\Otest\present{\Comma\bgroup\Ostr\egroup.}\else{.}\fi%
\annotations  %  MODIFICATION
\vskip5ptplus1ptminus1pt}%\smallskip (but 3 is now 5)  MODIFICATION

\def\otherformat{\Reffont\let\uchyph=1\parindent=1.25pc\def\Comma{}%

\sfcode`\.=1000\sfcode`\?=1000\sfcode`\!=1000\sfcode`\:=1000\sfcode`\;=1000\sfcode`\,=1000%\frenchspacing
                \par\vfil\penalty-200\vfilneg%\filbreak
      \if\Ftest\present\Flagstyle\Fstr\fi%

\if\Atest\present\bgroup\Authfont\Astr\egroup\def\Comma{\unskip\endauthor}%

\else\if\Etest\present\bgroup\def\Eand{\Aand}\def\Eandd{\Aandd}\Authfont\Estr\egroup\unskip, \ifnum\Ecnt>1eds.\else ed.\fi\def\Comma{, }%

\else\if\Itest\present\bgroup\Authfont\Istr\egroup\def\Comma{, }\fi\fi\fi%
          \if\Ttest\present\Comma\bgroup\Titlefont\Tstr\egroup\def\Comma{,
}\fi%
            \if\Atest\present\if\Itest\present
                    \Comma\bgroup\Istr\egroup\def\Comma{, }\fi%
                 \else\if\Etest\present\if\Itest\present
                         \Comma\bgroup\Istr\egroup\def\Comma{, }\fi\fi\fi%
                 \if\Ctest\present\ (\bgroup\Cstr\egroup)\def\Comma{, }\fi%
                  \if\Dtest\present\Comma\bgroup\Dstr\egroup\def\Comma{, }\fi%
                  \if\Ptest\present\bgroup, \Pstr\egroup\def\Comma{, }\fi%
                   \if\Gtest\present{\Comma Gov't ordering no.
}\bgroup\Gstr\egroup\def\Comma{, }\fi%
                    \if\Mtest\present\Comma MR
\#\bgroup\Mstr\egroup\def\Comma{, }\fi%
                     \if\Otest\present{\Comma\bgroup\Ostr\egroup.}\else{.}\fi%
\annotations  %  MODIFICATION
\vskip5ptplus1ptminus1pt}%\smallskip (but 3 is now 5)  MODIFICATION
\catcode`\@=11
\def\matrixx#1{\null\,\vcenter{\normalbaselines\m@th
    \ialign{\hfil$\displaystyle ##$\hfil&&\quad\hfil$\displaystyle
##$\hfil\crcr
    \mathstrut\crcr\noalign{\kern-\baselineskip}
    #1\crcr\mathstrut\crcr\noalign{\kern-\baselineskip}}}\,}
\catcode`\@=12

\catcode`\@=11
\def\ATstartsection{\@startsection}
\def\resetATfont{\reset@font}
\def\zAT{\z@}
\new@mathgroup\msa@group
\new@mathgroup\msb@group
\def\RIfM@{\relax\protect\ifmmode}
\def\fraktur{\mathfrak}   

\catcode`\@=12

\def\hexnumberZ#1{\ifnum#1<10 \number#1\else
 \ifnum#1=10 A\else\ifnum#1=11 B\else\ifnum#1=12 C\else
 \ifnum#1=13 D\else\ifnum#1=14 E\else\ifnum#1=15 F\fi\fi\fi\fi\fi\fi\fi}

\newcounter{alphactr}

\newcounter{romanctr}

\newcounter{arabicctr}

\def\geom{g\'eom\'etrie}
\def\EXPOSE{expos\'e}
\def\SGA{S\'em\-in\-aire de G\'eom\'e\-trie Al\-g\'e\-bri\-que}
\def\PLUS{+}
\def\TIMES{*}

\def\begindiagram{\def\normalbaselines{\baselineskip20pt
      \lineskip3pt \lineskiplimit3pt}}

\catcode`\@=11
\def\operatoratfont{\operator@font}
\let\@afterindentfalse\@afterindenttrue
\@afterindenttrue
\catcode`\@=12
\def\op{\operatoratfont}

\catcode`\@=11
\def\newtheoremrm#1{\@ifnextchar[{\@othmrm{#1}}{\@nthmrm{#1}}}
\def\@nthmrm#1#2{\@ifnextchar[{\@xnthmrm{#1}{#2}}{\@ynthmrm{#1}{#2}}}
\def\@xnthmrm#1#2[#3]{\expandafter\@ifdefinable\csname #1\endcsname%
{\@definecounter{#1}\@addtoreset{#1}{#3%
}\expandafter\xdef\csname the#1\endcsname{\expandafter\noexpand%
\csname the#3\endcsname . \noexpand\arabic{#1}}\global\@namedef{#1%
}{\@thmrm{#1}{#2}}\global\@namedef{end#1}{\endtrivlist}}}
\def\@ynthmrm#1#2{\expandafter\@ifdefinable\csname #1\endcsname%
{\@definecounter {#1}\expandafter\xdef\csname the#1\endcsname%
{\noexpand\arabic{#1}}\global\@namedef{#1}{\@thmrm{#1}{#2}%
}\global\@namedef{end#1}{\endtrivlist}}}
\def\@othmrm#1[#2]#3{\expandafter\@ifdefinable\csname #1\endcsname
{\global\@namedef{the#1}{\@nameuse{the#2}}\global\@namedef{#1%
}{\@thmrm{#2}{#3}}\global\@namedef{end#1}{\endtrivlist}}}
\def\@thmrm#1#2{\refstepcounter{#1}\@ifnextchar[{\@ythmrm{#1}{#2}%
}{\@xthmrm{#1}{#2}}}
\def\@xthmrm#1#2{\trivlist \item[\hskip \labelsep{\bf #2\ \csname
the#1\endcsname}]\ignorespaces}
\def\@ythmrm#1#2[#3]{\trivlist \item[\hskip \labelsep{\bf #2\ \csname
the#1\endcsname\ (#3)}]\ignorespaces}
\catcode`\@=12

\newtheoremrm{remark}[theorem]{Remark}

\def\many#1{\manyxx#1 \cdots \manyxx#1}
\def\manyxx#1{\if#1\PLUS{+}\else\if#1\TIMES{*}\fi\fi}

\def\xhmode#1{\relax\ifmmode{#1}\else{\snap\hbox{$#1$}}\fi}
\def\snap{\hskip 0pt plus 4cm\penalty1000\hskip 0pt plus -4cm}
\def\PP{\xmode{\Bbb P\kern1pt}}
\def\email{jaffe{\kern0.5pt}@{\kern0.5pt}cpthree.unl.edu}
\font\bbtwo=msbm10 scaled\magstep2
\def\CP{{\Bbb C}\kern1pt{\Bbb P}}
\def\C{{\Bbb C}\kern1pt}
\def\Q{{\Bbb Q}\kern1pt}
\def\P#1{{\PP^{#1}}}
\def\N{\xmode{\Bbb N}}
\def\Z{\xmode{\Bbb Z}}
\def\A{\xmode{\Bbb A\kern1pt}}
\def\O{{\cal O}}
\def\T{{\mathbf{T}}}                 %  tangent space
 \def\Cl{\mathop{\operatoratfont Cl}\nolimits}
\def\Ker{\mathop{\operatoratfont Ker}\nolimits}
\def\NS{\mathop{\operatoratfont NS}\nolimits}
\def\Pic{\mathop{\operatoratfont Pic}\nolimits}
\def\Proj{\mathop{\operatoratfont Proj}\nolimits}
\def\Sing{\mathop{\operatoratfont Sing}\nolimits}
\def\Spec{\mathop{\operatoratfont Spec}\nolimits}
\def\blowup{\mathop{\operatoratfont blowup}\nolimits}
\def\codim{\mathop{\operatoratfont codim}\nolimits}
\def\div{\mathop{\operatoratfont div}\nolimits}
\def\dlog{\mathop{\operatoratfont dlog}\nolimits}
\def\edge{\mathop{\operatoratfont edge}\nolimits}
\def\edges{\mathop{\operatoratfont edges}\nolimits}
\def\genus{\mathop{\operatoratfont genus}\nolimits}
\def\lcm{\mathop{\operatoratfont lcm}\nolimits}
\def\length{\mathop{\operatoratfont length}\nolimits}
\def\ord{\mathop{\operatoratfont ord}\nolimits}
\def\order{\mathop{\operatoratfont order}\nolimits}
\def\rank{\mathop{\operatoratfont rank}\nolimits}
\def\rem{\mathop{\operatoratfont rem}\nolimits}
\def\type{\mathop{\operatoratfont type}\nolimits}
\def\vertices{\mathop{\operatoratfont vertices}\nolimits}
\def\CHAR{\mathop{\operatoratfont char \kern1pt}\nolimits}
\def\DEF{\mathop{\operatoratfont def \kern1pt}\nolimits}
\def\mp[[#1||#2||#3]]{\snap\hbox{$#1:\kern3pt #2\ \to\ #3$}}
\def\mapx[[#1||#2]]{\snap\hbox{$#1\ \to\ #2$}}
\def\dmap[[#1||#2||#3]]%
{     \setbox0=\hbox{$\displaystyle #1:\kern3pt #2\ \ \mapE\ \ \ #3$%
}     \ifdim\wd0<\hsize$$#1:\kern3pt #2\ \ \mapE\ \ \ #3$$\else%
     \begin{flushleft} $\displaystyle #1:\kern3pt #2$
     \end{flushleft} \begin{flushright} $\displaystyle\mapE{}\ \ \ #3$
     \end{flushright}\fi}
\def\dmapx[[#1||#2]]%
{     \setbox0=\hbox{$\displaystyle #1\ \ \mapE\ \ \ #2$%
}     \ifdim\wd0<\hsize$$#1\ \ \mapE\ \ \ #2$$\else%
     \begin{flushleft} $\displaystyle #1$ \end{flushleft}
     \begin{flushright} $\displaystyle\mapE{}\ \ \ #2$ \end{flushright}\fi}
\def\disomorx[[#1||#2]]{$${#1}\ \ \iso\ \ {#2}$$}
\def\inverselim{\displaystyle{\lim_{\overleftarrow{\hphantom{\lim}}}}}
\def\mapE#1{{\smash{\mathop{\longrightarrow}\limits^{#1}}}}
\def\smapE#1{{\smash{\mathop{\kern-10pt\rightarrow\kern-10pt}\limits^{#1}}}}
\def\mapS#1{\Big\downarrow\rlap{$\vcenter{\hbox{$\scriptstyle{#1}$}}$}}
\def\diagramx#1{{$$\begindiagram\matrixx{#1}$$}}
\def\splitdiagram#1#2{
    \begin{flushleft}$\displaystyle\begindiagram\matrixx{#1}$\end{flushleft}
    \begin{flushright}$\displaystyle\begindiagram\matrixx{#2}$\end{flushright}}
\def\splitdisplay#1#2{
     \begin{flushleft}$\displaystyle{#1}$\end{flushleft}
     \begin{flushright}$\displaystyle{#2}$\end{flushright}}
\def\sescomma#1#2#3%
{  \setbox0=\hbox{$\matrixx{0&\mapE{}&#1&\mapE{}&#2&\mapE{}&#3&\mapE{}&0,}$%
}  \ifdim\wd0<\hsize\diagramx{0&\mapE{}&#1&\mapE{}&#2&\mapE{}&#3&\mapE{}&
  0,}\else\setbox0=\hbox{$\matrixx{0&\smapE{}&#1&\smapE{}&#2&\smapE{}&#3&
  \smapE{}&0,}$}\ifdim\wd0<\hsize\diagramx{0&\smapE{}&#1&\smapE{}&#2&\smapE{}&
  #3&\smapE{}&0,}\else\splitdiagram{0&\mapE{}&#1&\mapE{}&#2}{\mapE{}&#3
  &\mapE{}&0,}\fi\fi}
\def\les#1#2#3%
{  \setbox0=\hbox{$\matrixx{0&\mapE{}&#1&\mapE{}&#2&\mapE{}&#3}$%
}  \ifdim\wd0<\hsize\diagramx{0&\mapE{}&#1&\mapE{}&#2&\mapE{}&
  #3}\else\setbox0=\hbox{$\matrixx{0&\smapE{}&#1&\smapE{}&#2&\smapE{}&
  #3}$}\ifdim\wd0<\hsize\diagramx{0&\smapE{}&#1&\smapE{}&#2&\smapE{}&
  #3}\else\splitdiagram{0&\mapE{}&#1&\mapE{}&#2}{\mapE{}&#3}\fi\fi}
\def\sesmaps#1#2#3#4#5{\setbox0=\hbox{$\matrixx{0&\mapE{}&#1&\mapE{#2}&#3&
     \mapE{#4}&#5&\mapE{}&0}$}\ifdim\wd0<\hsize\diagramx{0&\mapE{}&#1&
     \mapE{#2}&#3&\mapE{#4}&#5&\mapE{}&0}\else\setbox0=\hbox{$\matrixx{0&
     \smapE{}&#1&\smapE{#2}&#3&\smapE{#4}&#5&\smapE{}&
     0}$}\ifdim\wd0<\hsize\diagramx{0&\smapE{}&#1&\smapE{#2}&#3&\smapE{#4}&#5
     &\smapE{}&0}\else\splitdiagram{0&\mapE{}&#1&\mapE{#2}&#3}{\mapE{#4}&#5&
     \mapE{}&0}\fi\fi}
\def\squareSE#1#2#3#4{\diagramx{#1&\mapE{}&#2\cr
		      \mapS{}&&\mapS{}\cr #3&\mapE{}&#4\cr}}
\def\makenull#1{{\setbox0=\hbox{#1}\wd0=0pt\box0}}
\def\nulldot{\makenull{.}}
\def\o#1{\if#1\PLUS\oplus\else\if#1\TIMES\otimes\fi\fi}
\def\manyo#1{\o#1 \cdots \o#1}
\def\tC{{\tilde{C}}}
\def\tS{{\tilde{S}}}
\def\tT{{\tilde{T}}}
\def\tY{{\tilde{Y}}}
\def\ltF{{\tilde{\lowercase{F}}}}
\def\ltX{{\tilde{\lowercase{X}}}}
\def\oA{{\overline{A}}}
\def\oC{{\overline{C}}}
\def\oS{{\overline{S}}}
\def\loA{{\overline{\lowercase{A}}}}
\def\loK{{\overline{\lowercase{K}}}}
\def\shA{{\cal{A}}}
\def\shD{{\cal{D}}}
\def\shJ{{\cal{J}}}
\def\shL{{\cal{L}}}
\def\hA{{\hat{A}}}
\def\hC{{\hat{C}}}
\def\hS{{\hat{S}}}
\def\lbC{{\mathbf{\lowercase{C}}}}
\def\lbD{{\mathbf{\lowercase{D}}}}
\def\lbE{{\mathbf{\lowercase{E}}}}
\def\lbH{{\mathbf{\lowercase{H}}}}
\def\lbR{{\mathbf{\lowercase{R}}}}
\def\lfM{{\xmode{{\fraktur{\lowercase{M}}}}}}
\def\lfP{{\xmode{{\fraktur{\lowercase{P}}}}}}
\def\hlfP{{\xmode{{\hat{\fraktur{\lowercase{P}}}}}}}
\newenvironment{warning}{\trivlist \item[\hskip \labelsep{\bf
Warning!\kern1pt}]}{\endtrivlist}
\def\STCI{set-theoretic complete intersection}
\def\STCIs{set-theoretic complete intersections}
\def\RDP{rational double point}
\def\IFF{if and only if}
\def\CM{Cohen-Macaulay}
\def\LHS{left hand side}
\def\RHS{right hand side}
\def\WMA{we may assume}
\def\WMAT{we may assume that}
\def\th#1{\xmode{{#1}^{\operatoratfont th}}}
\def\IN{\subset}
\def\notIN{\not\IN}
\def\setof#1{\xmode{\{#1\}}}
\def\iso{\cong}
\def\qed{{\hfill$\square$}}
\def\formulaqed#1{$$#1\eqno\square$$}
\def\nor#1{#1_{\operatoratfont nor}}
\def\RED#1{#1_{\operatoratfont red}}
\def\half{{1\over2}}
\def\vec#1#2#3{#1_{#2}, \ldots, #1_{#3}}
\def\svec#1#2#3{#1_{#2} \many+ #1_{#3}}
\def\xmode#1{\relax\ifmmode{#1}\else{$#1$}\fi}
\def\cat#1{\xhmode{\ll\negthinspace\hbox{#1}\negthinspace\gg}}
\def\dfunx[[ #1 || #2 ]]{\dmapx[[ \cat{#1} || \cat{#2} ]]}
\def\pref#1{(\ref{#1})}
\def\block#1{\section{#1}}
\def\abs#1{{\vert{#1}\vert}}
\def\floor#1{\lfloor#1\rfloor}
\def\makeaddress{
     \vskip 0.15in
     \par\noindent {\footnotesize Department of Mathematics and Statistics,
                                  University of Nebraska}
     \par\noindent {\footnotesize Lincoln, NE 68588-0323, USA\ \ (\email)}}
\def\paidforby{ \def\thefootnote{\fnsymbol{footnote}}
     \par\noindent David B. Jaffe\protect\footnote{Partially supported by
                   the National Science Foundation.}
     \makeaddress\def\thefootnote{\arabic{footnote}}\setcounter{footnote}{0}}
\def\br#1{{[#1]}}
\def\inn#1{\langle #1 \rangle}
\newenvironment{proof}{\trivlist \item[\hskip \labelsep{\sc
Proof.\kern1pt}]}{\endtrivlist}%\proof
 \newenvironment{proofnodot}{\trivlist \item[\hskip \labelsep{\sc
Proof}]}{\endtrivlist}%\proofnodot
 \newenvironment{sketch}{\trivlist \item[\hskip \labelsep{\sc
Sketch.\kern1pt}]}{\endtrivlist}%\sketch

\newenvironment{romanlist}{\begin{list}{(\roman{romanctr})}{\usecounter{romanctr}}}{\end{list}}%\romanlist

\newenvironment{arabiclist}{\begin{list}{(\arabic{arabicctr})}{\usecounter{arabicctr}}}{\end{list}}%\arabiclist
 \newenvironment{definition}{\trivlist \item[\hskip
       \labelsep{\bf Definition.\kern1pt}]}{\endtrivlist}
\newenvironment{examples}{\trivlist \item[\hskip
       \labelsep{\bf Examples.\kern1pt}]}{\endtrivlist}
\documentclass[12pt]{article}\usepackage{amssymb}
\newtheorem{theorem}{Theorem}[section]
 \newtheorem{lemma}[theorem]{Lemma}%\lemma
 \newtheorem{corollary}[theorem]{Corollary}%\corollary
 \newtheorem{prop}[theorem]{Proposition}%\prop
 \newtheorem{exampleth}[theorem]{Example}
\newenvironment{example}{\begin{exampleth}\fontshape{n}\selectfont}{\end{exampleth}}

\catcode`\@=11
\def\tableofcontents{\section*{Contents\markboth{CONTENTS}{CONTENTS}}
 \@starttoc{toc}}
\def\l@part#1#2{\addpenalty{\@secpenalty}
 \addvspace{2.25em plus 1pt} \begingroup
 \@tempdima 3em \parindent \z@ \rightskip \@pnumwidth \parfillskip
-\@pnumwidth
 {\footnotesize \rm \leavevmode #1\hfil \hbox to\@pnumwidth{\hss \ }}\par
 \nobreak \endgroup}
\def\l@special{\@dottedtocline{1}{0.0em}{2.3em}}
\def\l@section{\@dottedtocline{2}{1.5em}{2.3em}}  % 0 was 2, 0.0 was 2.3
\def\l@subsubsection{\@dottedtocline{3}{3.8em}{3.2em}}
\def\l@paragraph{\@dottedtocline{4}{7.0em}{4.1em}}
\def\l@subparagraph{\@dottedtocline{5}{10em}{5em}}
\def\listoffigures{\section*{List of Figures\markboth
 {LIST OF FIGURES}{LIST OF FIGURES}}\@starttoc{lof}}
\def\l@figure{\@dottedtocline{1}{1.5em}{2.3em}}
\def\listoftables{\section*{List of Tables\markboth
 {LIST OF TABLES}{LIST OF TABLES}}\@starttoc{lot}}
\let\l@table\l@figure
\catcode`\@=12

\begin{document}
{\par\noindent\Large\bf Applications of iterated curve blowup to}
\vskip 0.05in
{\par\noindent\Large\bf set theoretic complete intersections in
                        $\hbox{{\bbtwo P}}^3$}
\vskip 0.15in
\paidforby
{\footnotesize\tableofcontents}

\newpage

\section*{Introduction}
\addcontentsline{toc}{special}{Introduction}
\def\thefootnote{\arabic{footnote}}\setcounter{footnote}{0}

\indent
We describe some new results on the set-theoretic complete intersection
problem for projective space curves.  Fix an algebraically closed ground
field $k$.  Let $S, T \IN \P3$ be surfaces.
Suppose that $S \cap T$ is set-theoretically a smooth curve $C$ of degree
$d$ and genus $g$.  For purposes of the introduction, we label the main
results as A, B, Q, X, I, II, and III.\footnote{The actual numbering in the
text is A = \ref{thmA}, B = \ref{thmB}, Q = \ref{thmQ}, X = \ref{thmX},
I = \ref{thmI}, II = \ref{thmII}, III = \ref{thmIII}.}
The results I, II, and III are more technical than A, B, Q, and X.

Suppose that $S$ and $T$ have no common singular points.  We discover that
this requirement imposes severe limitations.  Indeed, theorem (X) asserts that
if $C$ is not a complete intersection, then
$\deg(S), \deg(T) < 2d^4$.  Fixing $(d,g)$, one can in fact form a finite
list of all possible pairs $(\deg(S),\deg(T))$, which is much shorter than
the list implied by theorem (X).  For instance, when $(d,g) = (4,0)$, and
assuming for simplicity that $\deg(S) \leq \deg(T)$, we find that
$$(\deg(S), \deg(T)) \in \{ (3,4), (3,8), (4,4), (4,7), (6,26), (9,48),
     (10,28),$$
$$(12,18), (13,16), (17,220), (18,118), (19,84), (20,67), (22,50), (28,33)
\}.$$
Very little is known about which of these degree pairs actually correspond
to surface pairs $(S,T)$.

Suppose that $S$ and $T$ have only rational singularities, and that the
ground field $k$ has characteristic zero.  We continue to assume that
$S$ and $T$ have no common singular points.  Under these conditions, we prove
(A) that $d \leq g + 3$.  (The actual statement is somewhat stronger.)

Suppose that $S$ is normal, and that $d > \deg(S)$.  Make no assumptions about
how the singularities of $S$ and $T$ meet.  Assume that $\CHAR(k) = 0$.
We show (Q) that $C$ is linearly normal.  In particular, it follows by
Riemann-Roch that $d \leq g+3$.

Suppose that $S$ is a quartic surface having only rational singularities.
Allow $T$ to be an arbitrary surface, and make no assumptions about how the
singularities of $S$ and $T$ meet.  Assume that $\CHAR(k) = 0$.  Under these
conditions, we prove (B) that $C$ is linearly normal.

In other papers\Lspace \Lcitemark 17\Rcitemark \Rspace{},\Lspace \Lcitemark
18\Rcitemark \Rspace{}, we
have proved the following complementary results (in characteristic zero):
if $S$ is has only ordinary nodes as singularities, or is a cone, or
has degree $\leq 3$, then $d \leq g + 3$.  It is conceivable (in characteristic
zero) that this
inequality is valid without any restrictions whatsoever on $S$ and $T$, or
even that $C$ is always linearly normal.
Examples of smooth \STCI\ curves in $\C\P3$ have been constructed by
Gallarati\Lspace \Lcitemark 8\Rcitemark \Rspace{}, Catanese\Lspace \Lcitemark
3\Rcitemark \Rspace{},
Rao (\Lcitemark 27\Rcitemark \ prop.\ 14), and the author\Lspace \Lcitemark
19\Rcitemark \Rspace{}.

To explain the results (I), (II), and (III), and to describe the methods by
which we prove (A), (B), and (X), there are two key ideas which must be
discussed%
.\footnote{We also give an alternate proof of the key ingredient of (X),
which is independent of the main machine of this paper.}
Both of these ideas have to do with the iterated blowing up of curves.

The first idea has to do with certain invariants $p_i = p_i(S,C)$
$(i \in \N)$ which we associate to a pair $(S,C)$ consisting of an abstract
surface $S$ and a smooth curve $C$ on $S$ such that $C \notIN \Sing(S)$.  Let
\mp[[ \pi || \tS || S ]] be the blowup along $C$.  Then $p_1(S,C)$ is the
sum of the multiplicities of the exceptional curves.  (See \S\ref{measure} for
details.)  Moreover, $\pi$ admits a unique section $\tC$ over $C$, so we
can define $p_2(S,C) = p_1(\tS,\tC)$, $p_3(S,C) = p_2(\tS,\tC)$, and so forth.
We refer to the sequence $(p_1, p_2, \ldots)$ as the {\it type\/} of $(S,C)$.
It is a sum of local contributions, one for each singular point of $S$ along
$C$, and it is a rather mysterious measure of how singular $S$ is along $C$.
The type depends not only on the particular species of singular points of
$S$ which lie on $C$, but also on the way in which $C$ passes through those
points.  For example, the local contribution to the type coming from an
$A_3$ singularity is either $(1,1,1,0,\ldots)$ or $(2,0,\ldots)$, depending
on how $C$ passes through the singular point.

The second idea is the following construction.  For this we assume
(as in the first paragraph) that $C = S \cap T$ (in $\P3$) and that
$C \notIN \Sing(S)$, $C \notIN \Sing(T)$.  Other than this, no restrictions
are necessary on the singularities of $S$ and $T$.
Let $Y_1$ denote the blowup of $\P3$ along $C$.  Let $S_1, T_1 \IN Y_1$
denote the strict transforms of $S$ and $T$ respectively.  Let $E_1 \IN Y_1$
be the exceptional divisor, which is a ruled surface over $\C$.  Then
$S_1 \cap E_1$ is a curve $C_1$ (mapping isomorphically onto $C$), together
with some rulings.  The total number of rulings, counted with multiplicities,
is $p_1(S,C)$.  Now let $Y_2$ be the blowup of $Y_1$ along $C_1$.  Let
$S_2, T_2, E_2 \IN Y_2$ be as above.  Then $S_2 \cap E_2$ is a curve $C_2$
plus $p_2(S,C)$ rulings.  Iterate this construction $n$ times, where $n$
is the multiplicity of intersection of $S$ and $T$ along $C$.  Then
$S_n \cap T_n$ is a union of strict transforms of rulings.  This fact leads
us to theorems (I) and (II), which are statements about the numbers $p_i$.
Theorem (III) is also such a statement, but it does not depend on the
construction we have just described.

We describe theorems (I), (II), and (III).  These depend on the data
$(s,t,d,g)$, where $s = \deg(S)$, $t = \deg(T)$.  To make this description
as simple as possible, we restrict our attention here to the special case
where $(s,t,d,g) = (4,4,4,0)$.

Theorem (I) has the hypothesis that $\Sing(S) \cap \Sing(T) = \emptyset$.
Its conclusion (applied to our special case) is that:
$$p_1 = p_2 = p_3 = 8.$%
$Theorem (I) is used in the proofs of (A) and (X).

Theorem (II) has no additional hypotheses.
Its conclusion (applied to our special case) is that:
\begin{eqnarray*}
p_1 & \geq & 8; \\
2p_1 + p_2 & \geq & 24; \\
8p_1 + 3p_2 + p_3 & \geq & 96.
\end{eqnarray*}
Theorem (II) is not used in the proofs of (A) or (B).

Theorem (III) has the hypotheses that $S$ has only rational singularities,
and that $\CHAR(k) = 0$.  Its conclusion (applied to our special case) is that:
$$\half p_1 + {1\over6} p_2 + {1\over12}p_3 \many+
              {1 \over k(k+1)}p_k + \cdots \geq 6.$%
$Theorem (III), or actually a minor variant of it, is used in the proof of (B).

We now mention some open problems and possible ways to improve upon the
results in this paper.

\par\noindent{\bf 1.}  Let $(S,C)$ be the local scheme $S$ of a normal surface
singularity, together with a smooth curve $C$ on $S$.  There are three
fundamental invariants of $(S,C)$ which are utilized in this paper.  Firstly,
there is the type of $(S,C)$.  Secondly,
there is the order of $(S,C)$, i.e.\ the smallest positive integer $n$ such
that $\O_S(nC)$ is Cartier.  Thirdly, there is $\Delta(S,C)$, which we
describe in \S\ref{def-section}.  What relationships exist between these three
invariants?  What is their relationship to the Milnor fiber?

\par\noindent{\bf 2.}
We suspect that (A), (B), and (III) are valid over an arbitrary
algebraically closed field.  There are significant difficulties in proving
this which we have not explored fully.  The proofs of (I) and (II) do not
depend on the characteristic.

\par\noindent{\bf 3.}  The proofs of (A) and (B) use a bound
\pref{bound-formula} on the number of exceptional curves in a minimal
resolution
for a surface $S \IN \C\P3$ having only rational singularities.  Formulate
and prove a suitable generalization for arbitrary normal surfaces.

\par\noindent{\bf 4.}  Construct examples of surfaces $S, T \IN \C\P3$,
having only rational singularities, meeting set-theoretically along a smooth
curve $C$, such that $\Sing(S) \cap \Sing(T) \not= \emptyset$.  The
only example we know of is where $\deg(S) = \deg(T) = 2$, and $C$ is a line.

\par\noindent{\bf 5.}  The generic hypothesis that $C \notIN \Sing(T)$ can
probably be eliminated.

\vspace*{0.1in}
\par\noindent{\footnotesize{\it Acknowledgements.}  I thank Dave Morrison
for helpful comments, and Juan Migliore for raising the issue of linear
normality of \STCIs.}
\section*{Conventions}
\addcontentsline{toc}{special}{Conventions}

\begin{arabiclist}

\item We fix an algebraically closed field $k$.

\item A {\it curve\/} [resp.\ {\it surface}] [resp.\ {\it three-fold}] is an
excellent $k$-scheme
such that every maximal chain of irreducible proper closed subsets has length
one [resp.\ two] [resp.\ three].  We make the following additional assumptions:
\begin{itemize}
\item all curves are reduced and irreducible;
\item in part III and the introduction, all surfaces are reduced and
      irreducible.
\end{itemize}

\item A surface {\it embeds in codimension one\/} if it can be exhibited
as an effective Cartier divisor on a regular three-fold.

\item A {\it variety} is an integral separated scheme of finite-type over $k$.

\item If $X$ and $Y$ are schemes, then the notation $X \IN Y$ carries the
implicit assumption that $X$ is a {\it closed subscheme\/} of $Y$.

\item In several situations, we use {\it bracketed exponents\/} to denote
{\it repetition\/} in sequences, and we drop trailing zeros, where appropriate.
For example,
$$(2,1^\br{4}) = (2,1,1,1,1,0, \ldots)$%
$and
$$(3^\br{\infty}) = (3,3,\ldots).$$

\item We use the Grothendieck convention regarding projective space
bundles.

\item For any variety $V$, we let $A^k(V)$ denote the group of
codimension $k$ cycles on $V$, modulo {\it algebraic\/} equivalence.  When
$d = \dim(V)$ and $V$ is complete, we identify $A^d(V)$ with $\Z$.

\end{arabiclist}

\vspace{0.25in}

\part{Local geometry of smooth curves on singular surfaces}

\block{Definitions}\label{def-section}

\par\indent\indent We define the category of {\it surface-curve pairs}.
(Sometimes, we use the shorthand term {\it pair\/} for a surface-curve pair.)
An {\it object\/} $(S,C)$
in this category consists of a surface $S$, together with a curve $C \IN S$,
such that $C$ is a regular scheme and $C \notIN \Sing(S)$.  A {\it morphism\/}
\mp[[ f || (S',C') || (S,C) ]] is a pair $(S' \mapE{} S, C' \mapE{} C)$ of
morphisms of $k$-schemes, such that the diagram:
\squareSE{C'}{C}{S'}{S%
}commutes, and such that if \mp[[ \phi || C' || C \times_S S' ]] is the
induced map, then $\phi \times_S \Spec \O_{S,C}$ is an isomorphism.  Such
a morphism $f$ is {\it cartesian\/} if $\phi$ is an isomorphism.  Most
properties
of morphisms of schemes also make sense as properties of morphisms in this
category: the properties are to be interpreted as properties of the morphism
\mapx[[ S' || S ]].

Let $(S,C)$ be a surface-curve pair.  We say that:
\begin{itemize}
\item $(S,C)$ is {\it geometric\/} if $S$ is a variety;
\item $(S,C)$ is {\it local\/} if $S$ is a local scheme;
\item $(S,C)$ is {\it local-geometric\/} if $S$ is a local scheme,
      essentially of finite type over $k$.
\end{itemize}

To give a local surface-curve pair $(S,C)$ is equivalent to giving the
data $(A,\lfP)$, consisting of an excellent local $k$-algebra $A$, of
pure dimension two,
together with a height one prime $\lfP \IN A$ such that $A_\lfP$ and
$A/\lfP$ are regular.  We write $(S,C) = \Spec(A,\lfP)$ to denote this
correspondence.

There are two operations on surface-curve pairs which we will be using.
Firstly, if $(S,C)$ is a local surface-curve pair, then the {\it completion\/}
$(\hS,\hC)$ makes sense and is also a local surface-curve pair.  Indeed,
if $(S,C) = \Spec(A,\lfP)$, then $A/\lfP$ is regular, and so $\hA/\hlfP$
is regular, since it equals $\widehat{A/\lfP}$, and the completion of a
regular local ring is regular.  The reader may also check easily that
$\hA_\hlfP$ is regular.  Moreover, $\hA$ is excellent, since any noetherian
complete local ring is excellent.  Note also: there is a canonical morphism
\mapx[[ (S,C) || (\hS,\hC) ]].

Secondly, for any surface-curve pair $(S,C)$, one can define the {\it blowup\/}
$(\tS,\tC)$ of $(S,C)$.  This is done by letting \mp[[ \pi || \tS || S ]] be
the blowup of $S$ along $C$, and by letting $\tC$ be the unique section
of $\pi$ over $C$, which exists e.g.\ by (\Lcitemark 11\Rcitemark \ 7.3.5).
There is a
canonical morphism \mapx[[ (\tS, \tC) || (S,C) ]].

Two local surface-curve pairs are {\it analytically isomorphic\/} if their
completions are isomorphic.

If $(S,C)$ is a surface-curve pair, and $p \in C$, we let $(S,C)_p$ denote
the corresponding local surface-curve pair.  A {\it configuration\/} is
an element of the free abelian monoid on the set of analytic isomorphism
classes of local-geometric pairs.

Let $(S,C)$ be a geometric surface-curve pair.  We may associate the
configuration:
$$\sum_{p \in \Sing(S) \cap C} [(S,C)_p]$%
$to $(S,C)$.  On occasion, we shall identify $(S,C)$ with the associated
configuration.

We are interested in invariants of a geometric surface-curve pair $(S,C)$
which depend only on the associated configuration.  There are four
such invariants which we shall consider:

\begin{arabiclist}

\item The {\it order\/} of $(S,C)$ is the smallest $n \in \N$ such that
$\O_S(nC)$ is Cartier, or else $\infty$ if $\O_S(nC)$ is not Cartier for
all $n \in \N$.  If $\Sing(S) \cap C = \setof{\vec p1k}$, then
$$\order(S,C) = \lcm\setof{\order(S,C)_{p_1}, \ldots, \order(S,C)_{p_k}},$%
$so the computation of the order is a purely local problem.  Moreover, at least
if $S$ is normal, the order depends only on the associated configuration.
Indeed, in that case, if $(S,C)$ is a local-geometric pair, then $S$ is
excellent, so $\hS$ is normal, and so by (\Lcitemark 6\Rcitemark \ 6.12) one
knows that the
canonical map \mapx[[ \Cl(S) || \Cl(\hS) ]] is injective.  Hence
$\order(S,C) = \order(\hS,\hC)$.

\item The {\it type\/} of $(S,C)$, which is the sequence $(p_i)_{i \in \N}$
discussed in the introduction, and studied in \S\ref{measure}.

\item Assume that $S$ is normal.  We define an invariant $\Delta(S,C) \in \Q$.
Let \mp[[ \pi ||  \tS || S ]] be a minimal resolution, and let
$\vec E1n \IN \tS$ be the exceptional curves.  Let $\tC \IN \tS$ be the strict
transform of $C$.  According to (\Lcitemark 26\Rcitemark \ p.\ 241),
there is a unique $\Q$-divisor $E = \sum a_i E_i$ such that
$(\tC + E) \cdot E_i = 0$ for all $i$.  We define
$\Delta(S,C) = -E^2$.  Then $\Delta(S,C)$ is independent of $\pi$.  If $S$
is projective, then $\Delta(S,C) = C^2 - \tC^2$, where $C^2$ is
defined in (\Lcitemark 26\Rcitemark \ p.\ 241).

\item Assume that $S$ has only rational double points along $C$.  Let
$\Sigma(S,C)$ equal the number of exceptional curves in the minimal
resolution of those singularities of $S$ which lie on $C$.

\end{arabiclist}
\block{The type of a surface-curve pair}\label{measure}

\par\indent\indent We define the {\it type\/} of a surface-curve pair, and
show that it is an analytic invariant, at least when the surface embeds in
codimension one.

\begin{definition}
Let $(S,C)$ be a surface-curve pair.  Let $(\tS,\tC)$ be the blowup of
$(S,C)$.  Let $\vec E1n \IN \tS$ be the (reduced)
exceptional curves.  We define numbers $p_i(S,C)$, for each $i \in \N$.
Define:
$$p_1(S,C) = \sum_{i=1}^n \length \O_{\pi^{-1}(C), E_i},$%
$where \mp[[ \pi || \tS || S ]] is the blowup map.
For $i \geq 2$, recursively define $p_i(S,C)$ by:
$$p_{i+1}(S,C) = p_i(\tS, \tC).$%
$The {\it type\/} of $(S,C)$ is the sequence $(p_1, p_2, \ldots)$.
\end{definition}

It is clear that the computation of the $p_i$ may be reduced to the computation
of the $p_i$ when $(S,C)$ is a local pair.

\begin{remark}\label{goober-peas}
Let $(S,C)$ be a surface-curve pair, and assume that $S$ embeds in codimension
one.  Then we have $S \IN T$ for some smooth three-fold $T$.  Let
\mp[[ \pi_S || \tS || S ]] and \mp[[ \pi_T || \tT || T ]] be the blowups
of $S$ and $T$ along $C$.  Then $\pi_S^{-1}(C) \iso \tS \cap E$, where
$E \IN \tT$ is the exceptional divisor.  This fact plays an absolutely
central role in our type computations.
\end{remark}

\begin{remark}
We do not know for which $(S,C)$ we have $p_1(S,C) \geq p_2(S,C)$, and hence
that $p_k(S,C) \geq p_{k+1}(S,C)$ for all $k \geq 1$.  Conceivably, these
inequalities may hold whenever $S$ embeds in codimension one, or even
whenever $S$ is \CM.  By explicit calculation,
we shall find in \pref{fantastico} that the inequalities hold if $S$ has only
rational double points along $C$.  However, as \pref{type-ex-2} shows,
for some $(S,C)$ one has $p_1(S,C) < p_2(S,C)$.
\end{remark}

\begin{remark}
We consider the following general question.  Let $(S,C)$ be a local-geometric
pair.  Assume that $S$ is not smooth.
Let $p \in S$ be the closed point.  Let $(\tS,\tC)$ be the blowup of
$S$ along $C$.  Let \mp[[ \pi || \tS || S ]] be the blowup map.
What is the structure of $X = \RED{\pi^{-1}(p)}$?  If
$S$ embeds in codimension one, then $X$ will be a $\P1$.  Weird things can
happen if $S$ is not \CM.  For example, in \pref{type-ex-2}, $X$ is isomorphic
to $\Proj \C[s,t,u] / (s^3 - t^2u)$, which is a rational curve with a cusp.
In \pref{type-ex-3}, $X$ is the disjoint union of a point and several copies
of $\P1$, which do not meet $\tC$.  The isolated point of $X$ is the unique
point of $\tC$ lying over $p$.  Assuming only that $S$ is \CM, we do not know
if $X$ is always isomorphic to $\P1$, or even if it is always connected.
However: if $S$ embeds in codimension two, then $X$ embeds in $\P2$.
\end{remark}

We will prove \pref{analytic-invariant} that the type of a local-geometric
pair $(S,C)$ is an analytic invariant, provided that $S$ embeds in
codimension one.  There are some preliminaries.

\begin{lemma}\label{formal-woof-1}
Let \mp[[ f || A || B ]] be a flat, formally smooth homomorphism of Artin
local rings.  Assume that $A$ contains a field.  Then $\length(A) =
\length(B)$.
\end{lemma}

\begin{proof}
Let $K$ and $L$ be the residue fields of $A$ and $B$.  Let \mp[[ i || K || L ]]
be the induced map.  Let $\lfM$ be the maximal ideal of $A$.  Then the map
\mapx[[ K = A/\lfM || B/\lfM B ]] is formally smooth, so $B/\lfM B = L$
and so $i$ is formally smooth.  A theorem of Cohen
(\Lcitemark 23\Rcitemark \ 28.J) implies that $A$ contains a
coefficient field, which we also denote by $K$.  Since $i$ is formally smooth,
$L/K$ is a separable field extension.  It follows by the cited theorem that
we may find a coefficient field $L$ for $B$ which contains $f(K)$.

Let $\oA = A \o*_K L$.  Then $f$ factors as:
\diagramx{A&\mapE{h}&\oA&\mapE{g}&B.%
}Since $i$ is formally smooth, so is $h$.
Since both $h$ and $g \circ h$ are formally smooth it follows by
(\Lcitemark 12\Rcitemark \ 17.1.4) that $g$ is formally smooth.

Clearly $B$ is a finite $\oA$-module.  In particular, $g$ is of finite-type,
so $g$ is smooth.  Since $g$ is smooth of relative dimension zero, $g$ is
\'etale.  By (\Lcitemark 12\Rcitemark \ 18.1.2), the obvious functor:
\dfunx[[ \'etale $\oA$-schemes || \'etale $L$-schemes ]]%
is an equivalence of categories, so $g$ is an isomorphism.  Hence
$B \iso A \o*_K L$.  Hence $\length(A) = \length(B)$.  \qed
\end{proof}

\begin{corollary}\label{formal-woof-2}
Let \mp[[ f || X' || X ]] be a flat, formally smooth morphism of irreducible
noetherian schemes.  Assume that $X$ is defined over a field.  Let $\eta$
and $\eta'$ be the generic points of $X$ and $X'$.  Then:
$$\length \O_{X',\eta'} = \length \O_{X,\eta}.$$
\end{corollary}

\begin{remark}
We do not know if \pref{formal-woof-1} and \pref{formal-woof-2} are true
without the hypothesis of being ``defined over a field''.
\end{remark}

\begin{lemma}\label{flat-is-cartesian}
Let \mp[[ f || (S_1,C_1) || (S_2,C_2) ]] be a flat morphism
of surface-curve pairs.  Then $f$ is cartesian.
\end{lemma}

\begin{proof}
We must show that the induced map \mp[[ \phi || C_1 || C_2 \times_{S_1} S_2 ]]
is an isomorphism.  It suffices to show that $C_1 = C_2 \times_{S_1} S_2$
as closed subschemes of $S_2$.  Let \mp[[ \pi || C_2 \times_{S_1} S_2 || C_2 ]]
be the projection map.  Because $\pi$ is flat, any irreducible component of
$C_2 \times_{S_1} S_2$ must dominate $C_2$.
(See e.g.{\ }\Lcitemark 14\Rcitemark \ III 9.7.)  But
$\phi \times_{S_2} \Spec \O_{S_2,C_2}$ is an isomorphism, so it follows that
$C_2 \times_{S_1} S_2$ is irreducible.  Since $C_1 = C_2 \times_{S_1} S_2$
at their generic points, they are equal as closed subschemes of $S_2$.  \qed
\end{proof}

\begin{prop}\label{formal-woof-3}
Let \mp[[ f || (S_1,C_1) || (S_2,C_2) ]] be a formally smooth, flat morphism
of surface-curve pairs.  Assume that the induced map
\mapx[[ C_1 || C_2 ]] is bijective.  Assume that $S_1$ and $S_2$ embed in
codimension one.  Then $(S_1,C_1)$ and $(S_2,C_2)$ have the same type.
\end{prop}

\begin{proof}
The subscript $i$ will always vary through the set $\setof{1,2}$.  Because
of our hypothesis on the map \mapx[[ C_1 || C_2 ]],
we may assume that $S_1, S_2$ are local schemes and that
$f$ is a local morphism.  Let \mp[[ \pi_i || (\tS_i, \tC_i) || (S_i, C_i) ]]
be the blowup maps.  By the universal property of blowing up, and because
$f$ is cartesian by \pref{flat-is-cartesian}, we obtain
a map \mp[[ \ltF || \tS_1 || \tS_2 ]] which makes the diagram:
\diagramx{\tS_1&\mapE{\ltF}&\tS_2\cr
         \mapS{\pi_1}&&\mapS{\pi_2}\cr
         S_1&\mapE{f}&S_2\cr%
}commute.  Furthermore, using the flatness of $f$, we see that this diagram
is cartesian and as a consequence that $\ltF$ is flat and formally smooth.
Since $S_1$ and $S_2$ embed in codimension one, so do $\tS_1$ and $\tS_2$.
Since $p_{k+1}(S_i,C_i) = p_k(\tS_i,\tC_i)$ for all $k \geq 1$,
the proof of the proposition will follow
if we can show that $p_1(S_1,C_1) = p_1(S_2,C_2)$.

Let $x_i \in S_i$ be the unique closed points.  It is clear that
$\pi_1^{-1}(x_1)$ maps onto $\pi_2^{-1}(x_2)$.  Moreover, $\ltF(\tC_1) =
\tC_2$.
Let $\ltX_i = \pi_i^{-1}(x_i) \cap \tC_i$.  Then $\ltF(\ltX_1) = \ltX_2$.
Let $P_i = C_i \times_{S_i} \tS_i$.  A little thought shows that there is
a cartesian diagram:
\squareSE{P_1}{P_2}{\tS_1}{\tS_2\nulldot%
}Since $S_1$ and $S_2$ embed in codimension one, $P_i = \tC_i \cup E_i$,
where $E_i \iso \P1$ and $\tC_i \cap E_i = \ltX_i$.  The equality
$p_1(S_1,C_1) = p_1(S_2,C_2)$ can then be deduced from \pref{formal-woof-2}.
\qed
\end{proof}

If $(S,C)$ is a local-geometric pair, then $S$ is excellent, so the completion
map \mapx[[ \hS || S ]] is formally smooth.  Hence we have:

\begin{corollary}\label{analytic-invariant}
If two local-geometric surface-curve pairs embed in codimension one and are
analytically isomorphic, then they have the same type.
\end{corollary}

\begin{remark}
We do not know if \pref{formal-woof-3} and \pref{analytic-invariant} are
true without the hypothesis of ``embedding in codimension one''.
\end{remark}
\block{Examples}

\par\indent\indent
We give three examples which illustrate type computations and pathological
aspects of blowing up.
Cf.\ \pref{fantastico}, where rational double points are dealt with.

The first example illustrates a general conjecture which we cannot yet make
precise: amongst surfaces of given degree in $\P3$, those which occur in
positive characteristic can have ``larger'' type than those which occur in
characteristic zero.  Of course, the type also depends on the choice of a curve
on the surface.

More specifically, the example shows that in characteristic two, a quartic
surface (together with a suitably chosen curve) can have $p_1 = p_2 = p_3 = 8$.
We expect that this cannot happen in characteristic zero.  If so, it would
follow from (\ref{thmI} = ``I'') that a smooth quartic rational curve
$C \IN \C\PP^3$ cannot be the \STCI\ of two quartic surfaces, unless $C$ is
contained in the singular locus of one of the surfaces.

On the other hand, Hartshorne
\Lcitemark 15\Rcitemark \Rspace{} and
Samuel\Lspace \Lcitemark 30\Rcitemark \Rspace{}  have shown that in positive
characteristic,
the monomial rational quartic curve $C \IN \P3$ is a \STCI.
See \pref{examplex-char-two} for additional comments.

The second two examples have to do with pairs $(S,C)$ in which $S$ is not
Cohen-Macaulay.  These seem to be of some intrinsic interest, but
have no direct relevance to the problem of \STCIs\ in $\P3$.  Example two
might be viewed as a statement about
the properties of the singularity at the vertex of the cone over a
space curve.  It would be very nice to understand better the connection
between this singularity and the properties of the space curve.

\begin{prop}\label{example-char-two}
Let $k$ be an algebraically closed field of characteristic two.
Let $S \IN \P3$ be the cuspidal cone given by $y^4 - x^3w = 0$.  Let
$C \IN S$ be the smooth rational quartic curve given by
$$(s,t) \mapsto (x,y,z,w) = (s^4, s^3t, st^3, t^4).$%
$Then $C$ meets $\Sing(S)$ at the unique point $(0,0,0,1)$, and the
type of $(S,C)$ is $(8,8,8)$.
\end{prop}

\begin{proof}
We will calculate
in the category of affine varieties, so we will replace $S$ by an affine
variety, and when we refer to a {\it blowup}, we will actually mean a correctly
chosen affine piece of the blowup.  We let $(S_n, C_n)$ denote the \th{n}
iterated blowup of $(S,C)$.

The assertion that $C \cap \Sing(S) = \setof{(0,0,0,1)}$ is easily checked.
Taking the affine piece at $w = 1$, we find that $S$ is given by
$y^4 = x^3$ and that $C$ is given by $x = yz$ and $y = z^3$.
Making the change of variable $x \mapsto x + yz$, followed by
$y \mapsto y+z^3$, we obtain the new equation:
$$y^4 + z^{12} = (x+yz+z^4)^3$%
$for $S$ and the equation $x = y = 0$ for $C$.

Blow up $S$ along $C$,
formally substituting $xy$ for $x$.  Then $S_1$ is given by:
$$y^3 = x^3y^2 + x^2y^2z + x^2yz^4 + xy^2z^2 + y^2z^3 + yz^6 + xz^8 + z^9.$%
$Intersecting with the exceptional divisor, as in \pref{goober-peas},
corresponds to setting $y = 0$.  We obtain $z^8(z + x) = 0$, which tells us
that $p_1(S,C) = 8$, and that $C_1$ is given by $y = 0$, $z+x=0$.  Making the
change of variable $z \mapsto z - x$, we obtain the new equation:
$$y^3 = y^2z^3 + yz^6 + x^4yz^2 + x^8z + z^9$%
$for $S_1$, and the equation $y = z = 0$ for $C_1$.

Blow up $S_1$ along
$C_1$, formally substituting $zy$ for $z$.  Then $S_2$ is given by:
$$y^2 = y^4z^3 + y^6z^6 + x^4y^2z^2 + y^8z^9 + x^8z.$%
$Setting $y = 0$, we obtain $x^8z = 0$, which tells us that $p_2(S,C) = 8$,
and that $C_2$ is given by $y = z = 0$.

Blow up $S_2$ along $C_2$, formally substituting $zy$ for $z$.
Then the blown up surface $S_3$ is given by:
$$y = y^6z^3 + y^{11}z^6 + x^4y^3z^2 + y^{16}z^9 + x^8z.$%
$Setting $y = 0$, we obtain $x^8z = 0$, which tells us that $p_3(S,C) = 8$,
and that the new curve $C_3$ is given by $y = z = 0$.  One can check
that $S_3$ is smooth along $C_3$, so $p_k(S,C) = 0$ for
all $k > 3$.  \qed
\end{proof}

\begin{prop}\label{type-ex-2}
Let $S \IN \A^4 = \Spec \C[x,y,z,w]$ be the cone over the monomial
quartic curve:
$$(s,t) \mapsto (x,y,z,w) = (s^4, s^3t, st^3, t^4)$%
$in $\P3$.  Let $C \IN S$ be the ruling given by $y = z = w = 0$.
Then $\type(S,C) = (1,2)$.  Moreover, if \mp[[ \pi || S_1 || S ]]
denotes the blowup along $C$, and $p \in S$ denotes the unique singular point,
then $\RED{\pi^{-1}(p)} \iso \Proj \C[s,t,u](s^3-t^2u)$.
\end{prop}

\begin{proof}
One sees that $S$ is given by the equations $yz = xw$, $x^2z = y^3$,
$z^3 = yw^2$, and $y^2w = xz^2$.  The blowup $S_1$ of $S$ along $C$ is
obtained\footnote{In this situation, where $S$ does not embed in codimension
one, it is apparently necessary to look at all of the affine pieces of the
blowup.  The details of this are left to the reader.  These calculations
are greatly facilitated by the use of a computer program such as Macaulay.}
by formally substituting $z = sy$, $w = ty$.  Then $S_1$ is
given by $sy = tx$, $sx^2 = y^2$, $s^3 = t^2$, and $ty = s^2x$.  Then:
\disomorx[[ \pi_1^{-1}(C) || \Spec \C[x,y,s,t]/(y,tx,sx^2,s^3-t^2,s^2x). ]]%
Set-theoretically,
$$\pi_1^{-1}(C) = V(s,t,y) \cup V(x,y, s^3-t^2).$%
$We have $C_1 = V(s,t,y)$.  Thus:
$$p_1(S,C) = \length \C[x,s,t]/(tx,sx^2,s^3-t^2,s^2x)_{(x,s^3-t^2)},$%
$which equals one.

The blowup $S_2$ of $S_1$ along $C_1$ is obtained by formally substituting
$s = ay$, $t = by$.  Then $S_2$ is given by $ax^2 = y$
and $b = a^2x$.  Let \mp[[ \pi_2 || S_2 || S_1 ]] be the
blowup map.  Then:
\begin{eqnarray*}
\pi_2^{-1}(C_1) & \iso & \Spec \C[a,b,x] / (ax^2,b-a^2x)\\
                & \iso & \Spec \C[a,x] / (ax^2).
\end{eqnarray*}
This implies that $p_2(S,C) = 2$.  Since $S_2$ is smooth, we see that
$\type(S,C)$ is as claimed.
\end{proof}

\begin{example}\label{type-ex-3}
Let $S$ be a smooth surface, which is a variety.  Fix $n \geq 2$, and let
$\vec p1n \in S$
be distinct (closed) points.  Let \mp[[ \pi || S || \oS ]] be the morphism
which pinches $\vec p1n$ together, yielding $p \in \oS$.  Let $C \IN S$
be a smooth curve passing through $p_1$ but not through $\vec p2n$.  Then
$\oC = \pi(C)$ is smooth.  Let \mp[[ f || X || \oS ]] be the blowup along
$\oC$.  Then $X$ is obtained from $S$ by blowing up $\vec p2n$.  Hence
$f^{-1}(\oC)$ is isomorphic to the disjoint union of $\oC$ with $n-1$
copies of $\P1$.  The type of $(S,C)$ is $(n-1)$.
\end{example}
\block{Classification of rational double point pairs}\label{class}

\par\indent\indent In this section, we assume that $k$ has characteristic
zero.  We describe a classification of local-geometric pairs $(S,C)$, up to
analytic isomorphism, where $S$ is the local scheme of a rational double
point singularity.

Let $(S,C)$ be a local-geometric pair corresponding to a \RDP.  As is
well-known, such objects $S$ are classified (up to analytic isomorphism)
by A-D-E Dynkin diagrams.  Let $\tS$ be the minimal resolution of $S$, and
let $\tC \IN \tS$ be the strict transform of $C$.  Let $\vec E1n \IN \tS$
be the exceptional curves, numbered as in
(\Lcitemark 18\Rcitemark \ p.\ 167).
Then $\tC$ meets a unique exceptional curve $E_k$, and we have
$\tC \cdot E_k = 1$.  Moreover, there are some restrictions on $k$, depending
on $S$.  (See\Lspace \Lcitemark 18\Rcitemark \Rspace{}\ 2.2.)

In this way, we are able to define certain local-geometric pairs $A_{n,k}$,
$D_{n,k}$, and $E_{n,k}$.  In fact, one can show\Lspace \Lcitemark 20\Rcitemark
\Rspace{}
that these pairs are well-defined, up to analytic isomorphism.  We have:

\begin{theorem}\label{the-conjecture}
Let $(S,C)$ be a local-geometric
pair, where $S$ is the local scheme of a rational double point
singularity.  Then $(S,C)$ is analytically isomorphic to a unique member of
the following list of local pairs:
\begin{itemize}
\item $A_{n,k}$ (for some positive integers $n,k$ with $k \leq (n+1)/2$);
\item $D_{n,1}$ (for some integer $n \geq 4$);
\item $D_{n,n}$ (for some integer $n \geq 5$);
\item $E_{6,1}$;
\item $E_{7,1}$.
\end{itemize}
\end{theorem}

\begin{remark}
Equations for these pairs may be found in the proof of \pref{fantastico}.
\end{remark}
\block{Invariants of rational double point configurations}\label{invrdp}

\par\indent\indent In this section, we assume that $k$ has characteristic zero.
We will calculate the type of $(S,C)$ in the case where $S$ is the local
scheme of a rational double point singularity.
This depends not only on $S$, but also on $C$.  Note that if $S$ is the
local scheme of any rational singularity, and $S$ embeds in a nonsingular
three-fold, then $S$ ``is'' a rational double point.

For each pair of positive integers $(n,k)$ with $k \leq n$, we define a
sequence $\phi(n,k)$ of integers, via the following recursive definition:
$$\phi(n,k) = \cases{ \phi(n,n-k+1),&if $k > {n+1\over2}$;\cr
                      (k),&if $k = {n+1\over2}$;\cr
                      (k,\phi(n-k,k)),&if $k < {n+1\over2}$.}$%
$
\begin{examples}
\
\begin{arabiclist}
\item $\phi(n,1) = (1^\br{n})$ for all $n \geq 1$;
\item $\phi(rk,k) = (k^\br{r-1},1^\br{k})$ for all $k \geq 1$, $r \geq 1$
      (generalizing 1);
\item $\phi(rk-1,k) = (k^\br{r-1})$ for all $r \geq 2$, $k \geq 1$
      (also generalizing 1);
\item $\phi(10,4) = (4,3,1^\br{3}).$
\end{arabiclist}
\end{examples}

Let $a, b \in \N$.  For each integer $n \geq 0$, we define the \th{n}
{\it iterated remainder\/} on division of $a$ by $b$, denoted $\rem_n(a,b)$.
Let $\rem_0(a,b) = b$, and let $\rem_1(a,b)$ be the usual remainder.  For
$n \geq 2$, define:
$$\rem_n(a,b) = \cases{ \rem_1(\rem_{n-2}(a,b), \rem_{n-1}(a,b)),
                        &if $\rem_{n-1}(a,b) \not= 0$;\cr
                        0,&if $\rem_{n-1}(a,b) = 0$.}$%
$
\par Let $a,b \in \N$.  For each integer $n \geq 1$, we define the \th{n}
{\it iterated quotient\/} of $a$ by $b$, denoted $\div_n(a,b)$.  Let
$\div_1(a,b) = \floor{a/b}$.  For $n \geq 2$, define:
$$\div_n(a,b) = \cases{ \div_1(\rem_{n-2}(a,b), \rem_{n-1}(a,b)),
                        &if $\rem_{n-1}(a,b) \not= 0$;\cr
                        0,&if $\rem_{n-1}(a,b) = 0$.}$%
$
\begin{prop}\label{key-rem}
Fix $k, n \in \N$ with $k \leq (n+1)/2$.  Let $t$ be the largest integer such
that $\rem_t(n-k+1,k) \not= 0$.  Let $r_i = \rem_i(n-k+1,k)$,
$d_i = \div_i(n-k+1,k)$, for various $i$.  Then:
$$\phi(n,k) = (r_0^\br{d_1}, r_1^\br{d_2}, \ldots, r_t^\br{d_{t+1}}).$%
$\end{prop}

\begin{sketch}
Define $r_{-1} = n-k+1$.  One shows that for all $p \geq 0$,
$$\phi(r_{p-1}+r_p-1, r_p) =
     \cases{(r_p^\br{d_{p+1}}, \phi(r_p + r_{p+1}-1, r_{p+1})),
            &if $r_{p+1} \not= 0$;\cr
            (r_p^\br{d_{p+1}}),&if $r_{p+1} = 0$.\cr}$%
$The result then follows by induction.  \qed
\end{sketch}

\begin{prop}\label{fantastico}
The type of $A_{n,k}$ is $\phi(n,k)$.  The type of $D_{n,1}$ is $(2)$.  We
have:
$${\begindiagram
\type(D_{n,n}) = \cases{({n\over2}),&if $n$ is even\kern1.5pt$;$\cr
                          ({n-1\over2}, 1^\br{n-1}),&
                               if $n$ is odd.\cr}}$%
$The type of $E_{6,1}$ is $(2,2)$.  The type of $E_{7,1}$ is $(3)$.
\end{prop}

\begin{proof}
We let $(S,C)$ correspond to the given pair.
The comments in the first paragraph of the proof of \pref{example-char-two}
apply equally well here.
We make use of the explicit resolutions of rational double points given
in the appendix to\Lspace \Lcitemark 25\Rcitemark \Rspace{}.

First we consider the $A_{n,k}$ case.  (We allow $1 \leq k \leq n$.)  Then
$S$ is given by $xy - z^{n+1} = 0$, and $C$ is given parametrically by
$x = u^k$, $y = u^{n-k+1}$, $z = u$.  [In terms of the notation used in
\Lcitemark 25\Rcitemark \Rspace{}, this may be seen as the image of
$V(u_k=1) \IN W_k$.]  After making the change of variable
$x \mapsto x + z^k$ and $y \mapsto y + z^{n-k+1}$, we find that $S$ is
given by:
$$xy + yz^k + xz^{n-k+1} = 0,\eqno(*)$%
$and that $C$ is given by $x = y = 0$.  From now on, we assume that
$k \leq {n+1\over2}$.  Blow-up along $C$,
formally substituting $yx$ for $y$.  Then $S_1$ is given by:
$$xy + yz^k + z^{n-k+1} = 0.$%
$Intersecting with the exceptional divisor, as in \pref{goober-peas},
corresponds to setting $x = 0$.  We obtain $z^k(y + z^{n-2k+1}) = 0$.  This
tells us that $p_1(S,C) = k$ and that $C_1$ is given by $x = 0$ and
$y + z^{n-2k+1} = 0$.  After making the change of variable
$y \mapsto y - z^{n-2k+1}$, and thence $y \mapsto -y$, we obtain the
equation:
$$xy + yz^k + xz^{n-2k+1} = 0\eqno(**)$%
$for $S_1$, and the equation $x = y = 0$ for $C_1$.  If $n-2k+1=0$, then
$S_1$ is smooth along $C_1$, and we are done.  Otherwise, compare
$(*)$ with $(**)$, to complete the $A_{n,k}$ case.

Now we deal with the case $D_{n,1}$.  In terms of the notation used
in\Lspace \Lcitemark 25\Rcitemark \Rspace{}, $C$ is the image of
$V(v_0 = 0) \IN W_0$.
Following\Lspace \Lcitemark 25\Rcitemark \Rspace{}, we would have two cases
($n$ even, $n$ odd), but in fact these two cases are identical in this
situation, after interchanging variables $(x \leftrightarrow y)$.
We find that $S$ is given by:
$$x^2z + y^2 - z^{n-1} = 0,$%
$and that $C$ is given by $y = z = 0$.  Blow up along $C$,
formally substituting $zy$ for $z$.  Then $S_1$ is given by:
$$x^2z + y - y^{n-2}z^{n-1} = 0.$%
$Setting $y = 0$, we obtain $x^2z = 0$.  This tells us that $p_1(S,C) = 2$
and that $C_1$ is given by $y = z = 0$.  An easy calculation shows that
$S_1$ is smooth.  The result for $D_{n,1}$ follows.

Now we deal with the case $D_{n,n}$.  In terms of the notation used in
\Lcitemark 25\Rcitemark \Rspace{}, $C$ is the image of
$V(u_n = 0) \IN W_n$.  We may take the same equation for $S$ as we did
in the case $D_{n,1}$.  There are two cases:

Case I: $n$ is even.  Then $C$ is given parametrically by $x = u^{(n-2)/2}$,
$y = 0$, $z = u$.  After making the change of variable
$x \mapsto x + z^{(n-2)/2}$, we find that $S$ is given by:
$$x^2z + y^2 + 2xz^{n/2} = 0,$%
$and that $C$ is given by $x = y = 0$.  Blow up along $C$,
formally substituting $xy$ for $x$.  Then $S_1$ is given by:
$$x^2yz + y + 2xz^{n/2} = 0.$%
$Setting $y = 0$, we obtain $xz^{n/2} = 0$.  This tells us that
$p_1(S,C) = n/2$ and that $C_1$ is given by $x = y = 0$.  On checks that
$S_1$ is smooth.

Case II: $n$ is odd.  Then $C$ is given parametrically by $x = 0$,
$y = u^{(n-1)/2}$, $z = u$.  (In this case, $x$ and $y$ are interchanged from
the notation in\Lspace \Lcitemark 25\Rcitemark \Rspace{}.)  After making the
change of variable $y \mapsto y + z^{(n-1)/2}$, we find that $S$ is given by:
$$x^2z + y^2 + 2yz^{(n-1)/2} = 0,$%
$and that $C$ is given by $x = y = 0$.  Blow up along $C$,
formally substituting $yx$ for $y$.  Then $S_1$ is given by:
$$xz + xy^2 + 2yz^{(n-1)/2} = 0.\eqno(*{*}*)$%
$Setting $x = 0$, we obtain $yz^{(n-1)/2} = 0$.  This tells us that:
$$p_1(S,C) = (n-1)/2$%
$and that $C_1$ is given by $x = y = 0$.  Now blow-up
along $C_1$, formally substituting $xy$ for $x$.  Then the blow-up
$S_2$ is given by:
$$xz + xy^2 + 2z^{(n-1)/2} = 0.$%
$Setting $y = 0$, we obtain $z(x + 2z^{(n-1)/2 - 1}) = 0$.  This tells us
that $p_2(S,C) = 1$ and that $C_2$ is given by
($y = 0$ and $x + 2z^{(n-1)/2 - 1} = 0$).  After making the change of
variable $x \mapsto x - 2z^{(n-1)/2 - 1}$, we find that $S_2$ is given by:
$$xz + xy^2 - 2y^2z^{(n-1)/2 - 1} = 0,$%
$and that $C_2$ is given by $x = y = 0$.

Now blow up along $C_2$, formally substituting $xy$ for $x$.  Then
the blown up surface $S_3$ is given by:
$$xz + xy^2 - 2yz^{(n-1)/2-1} = 0,$%
$$C_3$ is given by $x = y = 0$, and $p_3(S,C) = 1$.  Replacing $y$ by
$-y$, we may assume that $S_3$ is given by:
$$xz + xy^2 + 2yz^{(n-1)/2-1} = 0.$%
$This looks like $(*{*}*)$, except that $n$ is now replaced by $n-2$.
Note that if $n = 1$, then $(*{*}*)$ is smooth.  A little
thought shows that the asserted type of $D_{n,n}$ is correct.  A posteriori,
we see that $(S_1,C_1) = A_{n-1,1}$.  A direct proof of this assertion would
of course simplify the proof.

For both $E_{6,1}$ and $E_{7,1}$, we may choose any smooth curve for $C$.

For $E_{6,1}$, $S$ is given by $x^2 - y^3 - z^4 = 0$, and $C$ is given by
$y = 0$, $x + z^2 = 0$.  After making the change of variable
$x \mapsto x - z^2$, we obtain the new equation $x^2 - 2xz^2 - y^3 = 0$
for $S$.  Then $C$ is given by $x = y = 0$.  Blow up along $C$, substituting
$xy$ for $x$.  The equation for $S_1$ is:
$$x^2y - 2xz^2 - y^2 = 0.$%
$Setting $y = 0$, we obtain $xz^2 = 0$.  Hence $p_1(S,C) = 2$, and $C_1$ is
given by $x = y = 0$.  Blow up $S_1$ along $C_1$, substituting $xy$ for $x$.
The equation for $S_2$ is $x^2y^2 - 2xz^2 - y = 0$.  Substituting
$y = 0$, we obtain $xz^2 = 0$.  Hence $p_2(S,C) = 2$.  One checks that
$S_2$ is smooth, so $p_k(S,C) = 0$ for all $k > 2$.

For $E_{7,1}$, $S$ is given by $x^2 + y^3 - yz^3 = 0$, and $C$ is given by
$x = y = 0$.  Blow up along $C$, substituting $yx$ for $y$.  The equation
for $S_1$ is $x + x^2y^3 - yz^3 = 0$.  Setting $x = 0$, we obtain $yz^3 = 0$.
Hence $p_1(S,C) = 3$.  As $S_1$ is smooth, we see that the type of $E_{7,1}$
is as claimed.  \qed
\end{proof}

\begin{warning}
Amongst the \RDP\ local-geometric pairs, those of the kind $(S,C) = D_{n,n}$
with $n$ odd ($n \geq 5$) are highly atypical.  The following phenomena happen
only for these special pairs:
\begin{romanlist}
\item $\Sigma(S,C) < \sum_{i=1}^\infty p_i(S,C)$;
\item $p_r(S,C) \not= 0$ for some $r > \order(S,C)$: see \pref{interesting}.
\end{romanlist}
\end{warning}

The calculation in the proposition allows one to compute not just the
type of a rational double point, but also the precise sequence of
(analytic equivalence classes of) local surface-curve pairs which
arise under successive blowups:

\begin{itemize}\label{blowupAnk}
\item $\blowup(A_{n,k})
     \cases{\hbox{is smooth},&if $k = {n+1\over2}$;\cr
            = A_{n-k,n-2k+1},&if ${n-k+1\over2} < k < {n+1\over2}$;\cr
            = A_{n-k,k},&if $k \leq {n-k+1\over2}$;}$
\item $\blowup(D_{n,1})$ is smooth;
\item $\blowup(D_{n,n})$ is smooth, if $n$ is even;
\item $\blowup(D_{n,n}) = A_{n-1,1},$ if $n$ is odd;
\item $\blowup(E_{6,1}) = A_{3,2}$;
\item $\blowup(E_{7,1})$ is smooth.
\end{itemize}

We now calculate $\order(S,C)$, where $S$ is the local scheme of a
rational double point.

\begin{prop}\label{order-calc}
The order of $A_{n,k}$ is the order of $\loK$ in $\Z/(n+1)\Z$.
The order of $D_{n,1}$ is $2$.  The order of $D_{n,n}$ is $2$ if $n$
is even, and it is $4$ if $n$ is odd.  The order of $E_{6,1}$ is $3$.
The order of $E_{7,1}$ is $2$.
\end{prop}

\begin{proof}
Let $(S,C)$ correspond to the given pair.
Some of the orders ($E_{6,1}$, $E_{7,1}$, $D_{n,1}$ ($n$ even) and
$D_{n,n}$ ($n$ even)) can be computed immediately if one knows the
abstract group $\Cl(S)$.  A list of these groups may be found in
(\Lcitemark 22\Rcitemark \ p.\ 258).  We do not use this approach.

Let $\tS$ be the minimal
resolution of $S$.  Let $\vec E1n \IN \tS$ be the exceptional curves.
Let $\tC \IN \tS$ denote the strict transform of $C$.  There is a unique
$\Q$-divisor $E = a_1E_1 \many+ a_nE_n$ such that
$\tC \cdot E_i = -E \cdot E_i$ for all $i$.
(See\Lspace \Lcitemark 26\Rcitemark \Rspace{}.)  The total transform $\oC \IN
\tS$
of $C$ is $\tC + E$; it is a $\Q$-divisor.  According to
(\Lcitemark 26\Rcitemark \ p.\ 242), $\oC$ is integral
(i.e.\ $E$ is integral, i.e.\ $\vec a1n \in \Z$) \IFF\ ``$C$ is locally
analytically equivalent to zero''.  Since $S$ is a rational double point,
this is equivalent to $[C] = 0$ in $\Cl(S)$.  As this discussion applies
not just to $C$, but also to positive integer multiples of $C$, we see
that the order of $(S,C)$
is the least positive integer $N$ such that $N(\vec a1n) \in \Z^n$.
Let $M$ be the inverse of the self-intersection matrix of the $E_i$.
Then $(\vec a1n)$ is the \th{k} column of $M$.  This may be computed from
an explicit formula for $M$, which one may find in
(\Lcitemark 18\Rcitemark \ p.\ 169).

In case $(S,C) = A_{n,k}$, one finds that:
$$a_i = \cases{-k(n-i+1)/(n+1),&if $i \geq k$;\cr
               ki/(n+1) - i,&if $i \leq k$.\cr}$%
$From this we calculate that $a_1 = k/(n+1) - 1$.  The proof for $A_{n,k}$
follows.

In case $(S,C) = D_{n,1}$, one finds that:
$$a_i = \cases{-1,&if $i \leq n-2$;\cr
               -1/2,&if $n-1 \leq i \leq n$.\cr}$%
$Hence $\order(D_{n,1}) = 2$.

In case $(S,C) = D_{n,n}$, one finds that:
$$a_i = \cases{-i/2,&if $i \leq n-2$;\cr
               -(n-2)/4,&if $i = n-1$;\cr
               -n/4,&if $i = n$.\cr}$%
$Hence $\order(D_{n,n})$ is as claimed.

We now deal with the two exceptional cases.  The inverses of the
self-intersection matrices do not appear in
\Lcitemark 18\Rcitemark \Rspace{}, and we omit them here for lack of space.
In case $(S,C) = E_{6,1}$, one finds that:
$$(\vec a1n) = (-{4\over3}, -{5\over3}, -2, -1, -{4\over3}, -{2\over3}),$%
$and in case $(S,C) = E_{7,1}$, one finds that:
\formulaqed{(\vec a1n) = (-{3\over2}, -2, -{5\over2}, -3, -{3\over2}, -2, -1).}
\end{proof}

It is interesting to note that for $D_{n,n}$ ($n$ odd), one has
$p_r(S,C) \not= 0$ for some $r > \order(S,C)$.  This does not occur for
the other \RDP\ pairs, as we shall see in \pref{interesting}.

\begin{lemma}\label{gcdgcd}
Let $k$ and $N$ be positive integers, with $k < N$.  Assume that
$k \nmid N$.  Then:
$$\floor{N/k} \leq (N-k) / \gcd(k,N).$$
\end{lemma}

\begin{proof}
First suppose that $k > N/2$.  Then $\floor{N/k} = 1$, so we must show that
$\gcd(k,N) \leq N-k$.  Indeed, if $x | k$ and $x | N$, then $x | (N-k)$,
so this is clear.

Hence \WMA\ that $k \leq N/2$.  Since $k \nmid N$, $\gcd(k,N) \leq k/2$.
Therefore it suffices to show that $(k/2)(N/k) \leq N-k$.  This follows
from $k \leq N/2$.  \qed
\end{proof}

\begin{corollary}\label{orderbound}
Let $k$ and $N$ be positive integers, with $k \leq N$.  Let $t$ be the
smallest positive integer such that $\rem_t(N,k) = 0$.  Let
$d_i = \div_i(N,k)$, for $i = 1, \ldots, t$.  Then:
$$\svec d1t \leq N/\gcd(k,N).$$
\end{corollary}

\begin{proof}
The case $k = N$ is clear, so \WMA\ that $k < N$.  If $t = 1$, the result
is clear.  Let $r_1 = \rem_1(N,k)$.  We may assume that $r_1 \not= 0$.
By induction on $t$, \WMA\ that:
$$\svec d2t \leq k/\gcd(r_1,k).$%
$Therefore it suffices to show that:
$$d_1 + {k \over \gcd(r_1,k)} \leq {N \over \gcd(k,N)}.$%
$One sees that $\gcd(r_1,k) = \gcd(k,N)$.  Therefore it suffices to show
that $d_1 \leq (N-k) / \gcd(k,N)$.  This follows from \pref{gcdgcd}.  \qed
\end{proof}

\begin{prop}\label{interesting}
Let $(S,C)$ be a local-geometric
pair corresponding to a \RDP.  Assume that $(S,C) \not= D_{n,n}$ for any
odd integer $n \geq 5$.  Then $p_r(S,C) = 0$ for all $r \geq \order(S,C)$.
\end{prop}

\begin{proof}
We utilize \pref{fantastico} and \pref{order-calc}.  The only nontrivial
case is $(S,C) = A_{n,k}$.  We may assume that $k \leq (n+1)/2$.
For any $a,b \in \N$, let $o(a,b)$ denote the
order of $\loA$ in $\Z/b\Z$.  In the notation of \pref{key-rem}, we must show
that:
$$\svec d1{t+1} < o(k,n+1).$%
$Translating to the notation of \pref{orderbound}
($N = n+1$), both $d_1$ and $t$ change by $1$.  The statement we need is:
$$(d_1 - 1) + \svec d2t < o(k,N).$%
$Since $o(k,N) = N/\gcd(k,N)$, this does follow from \pref{orderbound}.  \qed
\end{proof}

The content of the following proposition may be found in
(\Lcitemark 18\Rcitemark \ proof of 2.3, pp.\ 169-170).

\begin{prop}\label{delta-formulas}
We have:
\begin{eqnarray*}
\Delta(A_{n,k}) & = & k(n-k+1)/(n+1)\\
\Delta(D_{n,1}) & = & 1\\
\Delta(D_{n,n}) & = & n/4\\
\Delta(E_{6,1}) & = & 4/3\\
\Delta(E_{7,1}) & = & 3/2.
\end{eqnarray*}
\end{prop}
\block{Technical lemmas on rational double points}

\par\indent\indent In this section we assume that $k$ has characteristic
zero.  We prove various technical relationships between the
invariants of \RDP\ local-geometric pairs.  We use these results in part III.
The result \pref{rdp-one} appears to be of intrinsic interest.

\begin{prop}\label{potato-1}
Let $(S,C)$ be a local-geometric pair corresponding to a rational
double point.  Write:
$$\type(S,C) = (n_1^\br{k_1}, \ldots, n_r^\br{k_r})$%
$with $n_1 > \cdots > n_r \geq 1$ and $k_i \geq 1$ for each $i$.  Assume that
$r > 1$.  Then $k_r > 1$ and $n_r | n_{r-1}$.
\end{prop}

\begin{proof}
We use \pref{fantastico}.  The proposition is clear if $(S,C)$ is of
species $D$ or $E$.  Therefore we may assume that $(S,C)$ is of species
$A$.  In the notation of \pref{key-rem}, we may write:
$$\type(S,C) = (r_0^\br{d_1}, \ldots, r_t^\br{d_{t+1}}).$%
$Since $r_{t+1} = 0$, we have $r_{t-1} = r_t d_{t+1}$.  Hence $r_t | r_{t-1}$.
Hence $n_r | n_{r-1}$.  Since $r_{t-1} > r_t$, $d_{t+1} > 1$.  Hence $k_r > 1$.
\qed
\end{proof}

\begin{prop}\label{potato-2}
If $\type(A_{n,k}) = \type(A_{n',k'})$, where $k \leq (n+1)/2$ and
$k' \leq (n'+1)/2$, then $n = n'$ and $k = k'$.
\end{prop}

\begin{proof}
We use \pref{fantastico}.
Write $\type(A_{n,k}) = (r_0^\br{d_1}, \ldots, r_t^\br{d_{t+1}})$, as in
\pref{key-rem}.  Then $n = n' = r_0(d_1+1) + r_1 - 1$, and $k = k' = r_0$.
\qed
\end{proof}

\begin{lemma}\label{yechh}
Fix positive integers $k$ and $N$ with $k \leq N/2$.  Let $t$ be the smallest
positive integer such that $\rem_t(N,k) = 0$.  Let $r_i = \rem_i(N,k)$,
$d_i = \div_i(N,k)$, for various $i$.  Then:
\begin{romanlist}
\item If $t = 1$, then $r_0d_1^{-1} = k^2/N$.
\item If $t = 2$, then $(r_0-r_1)d_1^{-1} + r_1(d_1+d_2)^{-1} \leq k^2/N$.
\item If $t \geq 3$, then:
$$(r_0-r_1)d_1^{-1} + (r_1-r_2)(d_1+d_2)^{-1} +
                      r_2(d_1+d_2+1)^{-1} \leq k^2/N.$%
$\end{romanlist}
\end{lemma}

\begin{proof}
First suppose that $t = 1$.  Then $N = d_1 k$.  Hence
$k^2/N = k/d_1 = r_0d_1^{-1}$.  This proves (i).

Now suppose that $t = 2$.  Then $N = r_0d_1 + r_1$ and $r_0 = r_1d_2$.
We must show that:
$${r_0 - r_1 \over d_1} + {r_1 \over d_1 + d_2} \leq {r_0^2\over r_0d_1+r_1}.$%
$Substitute $r_0 = r_1d_2$, and cancel out $r_1$.  We must show:
$${d_2 - 1 \over d_1} + {1 \over d_1+d_2} \leq {d_2^2 \over d_1d_2 + 1}.$%
$Eliminating denominators, we find that we must show:
$$(d_1 - 1)(d_2 - 1) \geq 0,$%
$which is certainly true.

Finally, suppose that $t \geq 3$.  Then $N = r_0d_1 + r_1$ and
$r_0 = r_1d_2 + r_2$.  We must show that:
$${r_0-r_1 \over d_1} + {r_1-r_2 \over d_1+d_2} + {r_2 \over d_1+d_2+1}
     \leq {r_0^2 \over r_0d_1 + r_1}.$%
$Substitute $r_0 = r_1d_2 + r_2$.  We must show that:
$${r_1d_2 + r_2 - r_1 \over d_1} + {r_1 - r_2 \over d_1 + d_2} +
     {r_2 \over d_1 + d_2 + 1} \leq
     {(r_1d_2 + r_2)^2 \over r_1d_1d_2 + r_2d_1 + r_1}.$%
$Now cancel denominators.  (This is best done with the aid of a computer.)
We must show:
\splitdiagram{d_2r_1^2 - d_1^2d_2r_1^2 - d_1d_2^2r_1^2
  + d_1^2d_2^2r_1^2 - d_2^3r_1^2 + d_1d_2^3r_1^2 - d_1^2r_1r_2 - d_2r_1r_2%
}{ - d_1d_2r_1r_2
  + 2d_1^2d_2r_1r_2 - d_2^2r_1r_2 + d_1d_2^2r_1r_2 + d_1^2r_2^2 \geq 0.}
Equivalently, we must show that:
\splitdiagram{r_1^2d_2[ d_2^2(d_1-1) + d_1(d_1d_2 - d_1 - d_2)+1]%
}{ + r_1r_2[ d_1^2(d_2-1) + d_2^2(d_1-1)+ d_2(d_1^2-d_1-1)] + d_1^2r_2^2 \geq
0.}
If $d_1 \geq 2$ and $d_2 \geq 2$, this is clear.  Since $k \leq N/2$,
we have $d_1 \geq 2$.  Suppose that $d_2 = 1$.  Then the needed inequality
simplifies to:
$$r_1r_2(d_1^2 - 2) + d_1^2r_2^2 \geq 0,$%
$which is true.  \qed
\end{proof}

{}From \pref{yechh}, we obtain the following weaker statement, which we use in
\pref{rdp-one}:

\begin{corollary}\label{weirdness}
Fix positive integers $k$ and $N$, with $k \leq N/2$.  Let $t$ be the
smallest positive integer such that $rem_t(N,k) = 0$.  Let $r_i = \rem_i(N,k)$,
$d_i = \div_i(N,k)$, for various $i$.  Then:
$$(r_0-r_1)d_1^{-1} \many+ (r_{t-1}-r_t)(\svec d1t)^{-1} \leq k^2/N.$$
\end{corollary}

\begin{warning}
When we use \pref{weirdness}, the symbol $d_1$ will appear to have two
different values, differing by $1$: in the application \pref{rdp-one},
$d_1$ will be smaller by $1$.
\end{warning}

\begin{prop}\label{rdp-one}
Let $(S,C)$ correspond to a rational double point singularity.
Let $p_k = p_k(S,C)$, for each $k \in \N$.
Then:
$$\sum_{k=1}^\infty{1 \over k(k+1)}p_k \geq \Delta(S,C).$$
\end{prop}

\begin{proof}
We use \pref{delta-formulas} and \pref{fantastico}.
If $(S,C)$ is not of species $A$, then the proposition is proved by the
following table:

{\renewcommand{\arraystretch}{1.2}
\vspace*{0.1in}
\centerline{
\begin{tabular}{||c|c|c||}  \hline
singularity & $\sum_{k=1}^\infty {1\over k(k+1)} p_k$
          & $\Delta(S,C)$\\ \hline
$D_{n,1}$ & $1$ & $1$\\ \hline
$D_{n,n}$ ($n$ even) & $n/4$ & $n/4$\\ \hline
$D_{n,n}$ ($n$ odd)
          & ${n-1\over4}$ + $\sum_{k=2}^n {1\over k(k+1)}$
          & $n/4$\\ \hline
$E_{6,1}$ & $4/3$ & $4/3$\\ \hline
$E_{7,1}$ & $3/2$ & $3/2$\\ \hline
\end{tabular}    }
\vspace*{0.1in}}

Suppose that $(S,C) = A_{n,k}$ for some $n, k$.  We may assume that
$k \leq (n+1)/2$.  We must show that:
$$\sum_{j=1}^\infty {1 \over j(j+1)}\phi(n,k)_j \geq k(n-k+1)/(n+1).\eqno(*)$%
$Let $t$ be the largest integer such that $\rem_t(n-k+1,k) \not= 0$.  For
$i = 0, \ldots, t$, let $r_i = \rem_i(n-k+1,k)$,
$d_{i+1} = \div_{i+1}(n-k+1,k)$.  By \pref{key-rem}, we see that $(*)$ is
equivalent to:
\splitdisplay{r_0 \sum_{j=1}^{d_1} {1 \over j(j+1)} +
r_1 \sum_{j=d_1+1}^{d_1+d_2} {1 \over j(j+1)} \many+
r_t \sum_{j=\svec d1t + 1}^{\svec d1{t+1}} {1 \over j(j+1)}%
}{\geq {k(n-k+1) \over n+1}.%
}Note that for any $a,b \in \N$ with $a \leq b$,
$$\sum_{j=a+1}^b {1 \over j(j+1)} = {b \over b+1} - {a \over a+1}.$%
$Hence $(*)$ is equivalent to:
\splitdisplay{\left(r_0-r_1\right)\left({d_1 \over d_1+1}\right) \many+
      \left(r_{t-1}-r_t\right)\left({\svec d1t \over
      \svec d1t + 1}\right)%
}{+ r_t\left({\svec d1{t+1} \over \svec d1{t+1} + 1}\right)
\geq k(n-k+1)/(n+1).%
}This is equivalent to:
\splitdiagram{r_0 - (r_0 - r_1)(d_1 + 1)^{-1} - \cdots - (r_{t-1}-r_t)
     (\svec d1t + 1)^{-1}%
}{- r_t(\svec d1{t+1} + 1)^{-1} \geq k(n-k+1)/(n+1).%
}Since $r_0 = k$, this is equivalent to:
\splitdiagram{(r_0 - r_1)(d_1 + 1)^{-1} \many+ (r_{t-1}-r_t)
     (\svec d1t + 1)^{-1}}
{ + r_t(\svec d1{t+1} + 1)^{-1} \leq k^2/(n+1).%
}Let $N = n+1$.  Then this follows from \pref{weirdness}, and thence completes
the proof.  \qed
\end{proof}

\begin{definition}
Let $(S,C)$ be a local-geometric pair corresponding to a \RDP.  Then the
{\it deficiency\/} of $(S,C)$ is:
$$\DEF(S,C) = \Sigma(S,C) - \sum_{i=1}^\infty p_i(S,C).$$
\end{definition}

One always has $\DEF(S,C) \geq 0$, except for $D_{n,n}$, with $n$ odd,
$n \geq 5$.
\part{Iterated curve blowups}

\block{Intersection ring of a blow up}\label{oneblow}

\par\indent\indent
In this section we describe (without proof) the intersection ring of the
blow-up of a nonsingular variety along a nonsingular subvariety, following the
statements given in (\Lcitemark 7\Rcitemark \ 6.7, 8.3.9).  There are
two differences between the assertions we make
and the assertions made in\Lspace \Lcitemark 7\Rcitemark \Rspace{}.
Firstly, we work with cycles modulo algebraic equivalence,
rather than modulo rational equivalence.  Secondly, we have adjusted the signs
to reflect our convention regarding projective space bundles.

Let $X$ be a nonsingular closed subvariety of a nonsingular variety $Y$.
Let $d = \codim(X,Y)$, and assume that $d \geq 2$.  Let $N$ be the
normal bundle of $X$ in $Y$.

Let $\tY$ be the blow-up of $Y$ along $X$.  The exceptional divisor
is isomorphic to $\PP N^*$.  We use the following diagram to fix notation:

\diagramx{\PP N^*&\mapE{j}&\tY\cr
          \mapS{g}&&\mapS{f}\cr
          X&\mapE{i}&Y.\cr%
}Let $F = \Ker[ g^*(N^*)\ \mapE{\rm{can}}\ \O_{\PP N^*}(1)]$.  For each $k$,
there is a canonically split exact sequence:
\sesmaps{A^{k-d}(X)}{\delta}{A^{k-1}(\PP N^*) \o+ A^k(Y)}{\beta}{A^k(\tY)%
}of cycle groups modulo algebraic equivalence.  The maps are given by:

$$\delta(x) = (c_{d-1}(F) \cdot g^*(x), i_*(x))$%
$and
$$\beta(\ltX, y) = j_*(\ltX) + f^*(y).$%
$
This describes $A^*(\tY)$ as an abelian group.  The ring structure is
described by the following rules:

\begin{eqnarray*}
(f^*y) \cdot (f^*y') & = & f^*(y \cdot y') \\
(j_*\ltX) \cdot (j_*\ltX')&=&-j_*((c_1 \O_{\PP N^*}(1))\cdot\ltX \cdot \ltX')
\\
(f^*y) \cdot (j_*\ltX) & = & j_*((g^*i^*y) \cdot \ltX).
\end{eqnarray*}
\block{The intersection ring of an iterated curve blow-up}\label{iterate}

\par\indent\indent The result of this section is:

\begin{theorem}\label{iteration}
Let $Y_0 = \P3$.  Let $C_0 \IN Y_0$ be a nonsingular curve of degree $d$
and genus $g$.  Let $Y_1$ be the blow-up of $Y_0$ along $C_0$.  Choose a
smooth curve $C_1$ which lies on the exceptional divisor $E_1 \IN Y_1$
and which meets each ruling on $E_1$ exactly once.  Let $Y_2$ be the
blow-up of $Y_1$ along $C_1$.  Iterate this process: $Y_{k+1}$ is obtained
by blowing up a smooth curve $C_k \IN E_k \IN Y_k$.  We assume that $C_k$
meets each ruling on $E_k$ exactly once and that for all $k \geq 2$,
$C_k \not= E_k \cap E_{k-1,k}$, where $E_{k-1,k} \IN Y_k$ denotes the
strict transform of $E_{k-1}$.  Let $H \IN Y_0$ be a plane.
Let $\lbH = [H] \in A^1(Y_0)$.  Let $\lbE_k = [E_k] \in A^1(Y_k)$.
Let $\lbR_k \in A^1(E_k)$ denote the class of a ruling,
which we identify with its image in $A^2(Y_k)$.
Identify $\lbH$, $\lbE_k$ and $\lbR_k$ with their images
in the intersection ring $A^*(Y_n)$ of the $\th{n}$ iterated blow-up $Y_n$.
Then $A^k(Y_n)$ has as a basis:
$$[Y_n]\ \ (k=0);\ \lbH, \vec\lbE1n\ \ (k=1);\ \lbH^2, \vec\lbR1n\ \ (k=2);
       \ 1\ \ (k=3).$%
$This information, together with the following multiplication rules, completely
describe $A^*(Y_n)$ as a graded ring: $\lbH^3 = 1$, $\lbH \cdot \lbR_k = 0$,
$\lbH^2 \cdot \lbE_k = 0$, $\lbE_i \cdot \lbR_j = -\delta_{i,j}$,
$\lbH \cdot \lbE_k = d\lbR_k$,
$\lbE_i \cdot \lbE_j = -\beta_i\lbR_j\ \ (\hbox{if } i < j)$,
$$\lbE_k^2 = -d\lbH^2 - \alpha_{k-1}\lbR_k - \sum_{i=1}^{k-1} \beta_i\lbR_i,$%
$where $\alpha_k$ is determined by $[C_k] = c_1\O_{E_k}(1) - \alpha_k\lbR_k$
in $A^1(E_k)$, for $k \geq 1$, $\alpha_0 = 2-2g-4d$, and
$\beta_k = \alpha_{k-1} - \alpha_k$, for each $k \geq 1$.
\end{theorem}

We note the following generalization and conceptual reformulation of
\pref{iteration}, whose proof is omitted.  It will not be used again.

\begin{theorem}
Let $Y_0$ be a nonsingular complete three-fold.  Let $C_0 \IN Y_0$ be a
nonsingular curve.  Let $Y_1$ be the blow-up of $Y_0$ along $C_0$.  Choose a
smooth curve $C_1$ which lies on the exceptional divisor $E_1 \IN Y_1$
and which meets each ruling on $E_1$ exactly once.  Let $Y_2$ be the
blow-up of $Y_1$ along $C_1$.  Iterate this process: $Y_{k+1}$ is obtained
by blowing up a smooth curve $C_k \IN E_k \IN Y_k$.  We assume that $C_k$
meets each ruling on $E_k$ exactly once and that for all $k \geq 2$,
$C_k \not= E_k \cap E_{k-1,k}$, where $E_{k-1,k} \IN Y_k$ denotes the
strict transform of $E_{k-1}$.  Let $\lbE_k = [E_k] \in A^1(Y_k)$.
Let $\lbR_k \in A^1(E_k)$ denote the class of a ruling,
which we identify with its image in $A^2(Y_k)$.
Identify $\lbE_k$ and $\lbR_k$ with their images
in the intersection ring $A^*(Y_n)$ of the $\th{n}$ iterated blow-up $Y_n$.
Then $A^*(Y_n)$ is the graded $A^*(Y_0)$-algebra generated by
$\vec\lbE1n$ (degree $1$) and $\vec\lbR1n$ (degree $2$), modulo the relations:
$A^1(Y_0) \cdot \lbR_k = 0$, $A^2(Y_0) \cdot \lbE_k = 0$,
$\lbE_i \cdot \lbR_j = -\delta_{i,j}$,
$\lbH \cdot \lbE_k = (\lbH \cdot C_0)\lbR_k$\ (for all $\lbH \in A^1(Y_0)$),
$\lbE_i \cdot \lbE_j = -\beta_i\lbR_j\ \ (\hbox{if } i < j)$,
$$\lbE_k^2 = - [C_0] - \alpha_{k-1}\lbR_k - \sum_{i=1}^{k-1} \beta_i\lbR_i,$%
$where $\alpha_k$ is determined by $[C_k] = c_1\O_{E_k}(1) - \alpha_k\lbR_k$
in $A^1(E_k)$, for $k \geq 1$, $\alpha_0 = \deg(N_{C_0}^*)$, and
$\beta_k = \alpha_{k-1} - \alpha_k$, for each $k \geq 1$.
\end{theorem}

The remainder of this section breaks up into two parts.  First we
introduce various notations and conventions which we will use in the proof and
in subsequent sections.  Then we prove \pref{iteration}.

There are group homomorphisms \mapx[[ A^i(E_k) || A^{i+1}(Y_k) ]] and
injective ring homomorphisms:
\diagramx{A^*(Y_0)&\mapE{}&A^*(Y_1)&\mapE{}&\cdots&\mapE{}&A^*(Y_n).%
}We systematically identify various elements with their images, via these maps.
Since the latter maps are ring homomorphisms, it is not necessary to
distinguish
between multiplication in $A^*(Y_i)$ and $A^*(Y_j)$, for any $i,j$.  On the
other hand, since
the maps \mapx[[ A^i(E_k) || A^{i+1}(Y_k) ]] are not ring homomorphisms, it
is necessary to distinguish between multiplication in $A^*(Y_k)$ and
$A^*(E_k)$.
We do this by using a dot ($\cdot$) to denote multiplication in $A^*(Y_k)$
and brackets ($\inn,$) to denote multiplication in $A^*(E_k)$.  No problems
are introduced by the fact that $k$ does not occur explicitly in the bracket
notation.

Let $\lbC_k = [C_k]$.  This is an element of $A^2(Y_k)$, and it is an
element of $A^1(E_k)$ if $k \geq 1$.
Let $\lbD_k = c_1\O_{E_k}(1)$.  It is an element of $A^1(E_k)$.

For $k\leq n$,
let $E_{k,n} \IN Y_n$ denote the strict transform of $E_k$.  In $A^1(Y_n)$
we have $\lbE_k = [E_{k,n}] \many+ [E_{n,n}]$.  (This depends on our
assumption that $C_k \not= E_k \cap E_{k-1,k}$.)  {\it In particular, the
reader should observe the following insidious source of error:
$\lbE_k \not= [E_{k,n}]$.}  This same sort of error applies to other cycles
which we shall discuss.

In this section we do not fix a particular ruling $R_k \IN E_k$.  We do
so in the next section.  Having made such a choice, one can then discuss the
strict transform $R_{k,n} \IN Y_n$ of $R_k$.

Let $N_k$ be the normal bundle of $C_k$ in $Y_k$.  Then
$E_k \iso \PP(N_{k-1}^*)$.
For each $k = 0, \ldots, n$, we let $\alpha'_k = \deg(N_k^*)$.
For each $k = 1, \ldots, n$, we let $\beta'_k = \alpha'_{k-1} - \alpha_k$.
(We will show that $\alpha_k = \alpha'_k$ and hence that $\beta_k = \beta'_k$.)

\begin{proofnodot}
(of \ref{iteration}).
We make repeated use of the results of \S\ref{oneblow}, without explicitly
referring to them.  The abelian group structure of $A^*(Y_n)$ and
the assertions that $\lbH^3 = 1$, $\lbH \cdot \lbR_k = 0$,
$\lbH^2 \cdot \lbE_k = 0$, and $\lbE_i \cdot \lbR_j = -\delta_{i,j}$ are
left to the reader.

We compute $\lbH \cdot \lbE_i$.  Let $\mu_i = \lbH \cdot \lbC_i$.  Note that:
$$\lbH \cdot \lbE_i\ =\ \mu_{i-1} \lbR_i.$%
$We show that $\mu_k$ is independent of $k$, and in fact equals $d$.  First one
checks that $\lbH \cdot \lbC_0 = d$.  Now we have:
\begin{eqnarray*}
\mu_i & = & \lbH \cdot \lbC_i \\
      & = & \inn{ (\lbH \cdot \lbE_i), \lbC_i} \\
      & = & \inn{\mu_{i-1} \lbR_i, \lbC_i} \\
      & = & \inn{ \mu_{i-1} \lbR_i, \lbD_i - \alpha_i \lbR_i } \\
      & = & \mu_{i-1}.
\end{eqnarray*}
Hence $\mu_k = d$ for all $k$.  Hence $\lbH \cdot \lbC_i\ =\ d$ and
$\lbH \cdot \lbE_i\ =\ d\lbR_i$.

We now work on showing that $\alpha_k = \alpha'_k$.  In the process we
calculate
$\lbE_i \cdot \lbC_j$ for all $i \leq j$, a result we shall need later.
We have:

$$\inn{\lbD_i,\lbD_i} = c_1(N_{i-1}^*) = \alpha'_{i-1}.$%
$Further:
\begin{eqnarray*}
\lbE_i \cdot \lbC_i & = & -\inn{\lbD_i, \lbC_i} \\
                    & = & -\inn{\lbD_i, \lbD_i - \alpha_i\lbR_i} \\
                    & = & -(\alpha'_{i-1}-\alpha_i)\ =\ -\beta'_i.
\end{eqnarray*}
Using this we find:
\begin{eqnarray*}
\lbE_i \cdot \lbC_{i+1} & = & \inn{(\lbE_i \cdot
\lbC_i)\lbR_{i+1},\lbC_{i+1}}\\
                        & = & \lbE_i \cdot \lbC_i.
\end{eqnarray*}
Continuing in this manner, the reader may verify that for $i \leq j$,
$\lbE_i \cdot \lbC_j = -\beta'_i$.

The class of the canonical divisor on $Y_i$ is given by:
$$[K_{Y_i}] = -4\lbH + \lbE_1 \many+ \lbE_i.$%
$(This may be computed from the formula for the canonical divisor of a
blowup  --  see\Lspace \Lcitemark 9\Rcitemark \Rspace{}\ p.\ 608.)

For all $i \geq 0$, we have:
\begin{eqnarray*}\label{quirkalpha}
\alpha'_i & = & c_1 N_i^* \\
          & = & -c_1 \det(N_i)\ =\ -[[K_{C_i}] - [K_{Y_i}]|_{C_i}] \\
          & = & -[[K_{C_i}] - [K_{Y_i}] \ \cdot \lbC_i] \\
          & = & -[2g-2 - (-4\lbH + \lbE_1 \many+ \lbE_i) \cdot \lbC_i] \\
          & = & -[2g-2 + 4d + \beta'_1 \many+ \beta'_i].
\end{eqnarray*}

{}From this, and from the definition of the $\beta$'s and the $\alpha$'s,
we conclude:
$$\alpha'_k = \alpha_k\ \ \hbox{for all } k \geq 0.$%
$
\par\indent For $i < j$,
$$\lbE_i \cdot \lbE_j = (\lbE_i \cdot \lbC_{j-1})\lbR_j,$%
$so we obtain the formula $\lbE_i \cdot \lbE_j = -\beta_i\lbR_j$.

We proceed to calculate $\lbE_k^2$.  By the definition of the map $\delta$
given in \S\ref{oneblow}, we have:
\begin{eqnarray*}
\lbD_i & = & \alpha_{i-1}\lbR_i + \lbC_{i-1}\ \ \ (i \geq 1).
\end{eqnarray*}
Continuing to calculate, we find:
\begin{eqnarray*}\label{quirkC}
\lbC_i & = & \lbD_i - \alpha_i \lbR_i\ \ \ (i \geq 1) \\
\lbD_i & = & \alpha_{i-1}\lbR_i + \lbD_{i-1} - \alpha_{i-1}\lbR_{i-1}
                                                \ \ \ (i \geq 2) \\
\lbD_i - \lbD_{i-1} & = & \alpha_{i-1}\lbR_i - \alpha_{i-1}\lbR_{i-1}
                                                \ \ \ (i \geq 2) \\
\lbD_1 & = & \alpha_0 \lbR_1 + d\lbH^2 \\
\lbD_k & = & d\lbH^2 + \left( \sum_{i=1}^{k-1} (\alpha_{i-1} - \alpha_i)
             \lbR_i \right) + \alpha_{k-1}\lbR_k\ \ \ (k \geq 1) \\
\lbE_k^2 & = & -\lbD_k\ \ \ (k \geq 1) \\
         & = & - \left[ d\lbH^2 + \left( \sum_{i=1}^{k-1} \beta_i\lbR_i \right)
               + \alpha_{k-1}\lbR_k\right].
\end{eqnarray*}
\vspace*{-0.25in}

\qed
\end{proofnodot}
\block{The strict transform of a ruling}

\par\indent\indent The results of this section will be used in the proof of
theorem II \pref{thmII}.
The notations introduced in \S\ref{iterate} remain in effect in this section.

Fix a particular ruling $R_k \IN E_k$, where $1 \leq k \leq n$.  We compute
the class of $R_{k,n}$ in $A^2(Y_n)$.  A priori, this is a $\Z$-linear
combination of $\lbH^2, \vec\lbR1n$, which depends on the particular
choice of $R_k$.

We use the term {\it graph\/} to mean an undirected graph, which we shall
formally view as a reflexive, symmetric relation.
By an {\it augmented graph}, we shall mean a graph, together with a mapping
from the set of vertices of that graph to $\Z$.  If the augmentation map is
injective, we shall refer to the graph as a {\it labeled graph}, with the
obvious connotations.

Let $\Gamma$ be a labeled graph, which we suppose has a maximum vertex $m$.
We define various labeled graphs, coming from $\Gamma$, with maximum vertex
$m+1$.

First we define a labeled graph $\Gamma^+$ by
$\vertices(\Gamma^+) = \vertices(\Gamma) \cup \setof{m+1}$ and
$\edges(\Gamma^+) = \edges(\Gamma) \cup \setof{\edge(m,m+1)}$.

Now suppose that $\edge(l,m) \in \Gamma$.  We define a graph
$\Gamma^l$ by $\vertices(\Gamma^l) = \vertices(\Gamma) \cup \setof{m+1}$ and
$$\edges(\Gamma^l) = \edges(\Gamma) \cup \setof{\edge(l, m+1), \edge(m, m+1)}
                  - \setof{\edge(l,m)}.$%
$Intuitively, this construction may be thought of as adding a vertex $(m+1)$
``in the middle'' of the edge from $l$ to $m$.

\begin{definition}
A {\it standard operation\/} is an operation on a labeled graph of the
form $\Gamma \mapsto \Gamma^+$ or $\Gamma \mapsto \Gamma^l$ for some $l$.
A {\it standard labeled graph\/} is a labeled graph obtained from
a one-vertex labeled graph by a finite sequence of standard operations.
\end{definition}

It is not hard to see that given a standard labeled graph, one may
compute the last standard operation which was performed, and thence undo
that operation.  It follows that:

\begin{prop}\label{unique-operations}
Let $G$ be a standard labeled graph.  Then there is a unique sequence of
standard operations which gives rise to $G$.
\end{prop}

Fix integers $k$ and $m$ with $1 \leq k \leq m \leq n$.
Let $R_k \IN E_k$ be a ruling.
We will show how to associate a certain standard labeled graph
$\Gamma_m(R_k)$ to $R_k$, in such a way that
$[R_{k,m}] \in A^2(Y_m)$ depends only on $\Gamma_m(R_k)$.

To do this, consider the set of all curves $H \IN Y_m$ which are the
strict transforms of some ruling $R_l$ on $E_l$, for some $l$ with
$k \leq l \leq m$.  To each such $H$, we may associate an integer, namely $l$.
It may be that $H \IN E_{l',m}$, for some $l'$ with $l' \not= l$ and
$k \leq l' \leq m$, but this does not matter to us.  The set of all such
curves $H$ may be viewed as the vertices of a graph $\Gamma_m(k)$: two distinct
vertices are connected by an edge \IFF\ the corresponding two curves on $Y_m$
meet.  There is an augmentation on $\Gamma_m(k)$ given by $H \mapsto l$ as
above.

Define $\Gamma_m(R_k)$ to be the maximal connected subgraph of $\Gamma_m(k)$
which contains $R_{k,m}$.  The augmentation on $\Gamma_m(k)$ induces an
augmentation on $\Gamma_m(R_k)$.  We shall prove shortly \pref{snooker} that
$\Gamma_m(R_k)$ is a labeled graph, and that in fact it is a standard
labeled graph.

\begin{lemma}\label{lemma1}
If two distinct curves $H_1, H_2 \in \Gamma_m(k)$ meet, then they meet at a
unique point, and they meet transversally.
\end{lemma}

\begin{proof}
What we need to show is that if $p \leq q$ are integers
($k \leq p,q \leq m$), and if
$R_p \IN E_p$ and $R_q \IN E_q$ are rulings, and if $R_{p,m}$ meets
$R_{q,m}$ (but $R_{p,m} \not= R_{q,m}$), then in fact $R_{p,m}$ meets
$R_{q,m}$ at a unique point and they
do so transversally.  It suffices to show that $R_{p,q}$ meets
$R_q$ in this way.  We may assume that $p < q$.

Indeed if $R_{p,q}$ met $R_q$ at more than one point, or if they did not
meet transversally, then the image of $R_{p,q}$ under the map
\mapx[[ Y_q || Y_p ]] would be singular, because this map contracts $R_q$.
\qed
\end{proof}

\begin{lemma}\label{morsel}
No three distinct curves $H_1, H_2, H_3 \in \Gamma_m(k)$ meet at a common
point.
\end{lemma}

\begin{proof}
We may reduce to showing the following: if $p < q < r$
($k \leq p,q,r \leq m$) and
$R_p \IN E_p$, $R_q \IN E_q$, and $R_r \IN E_r$ are rulings, then
$R_{p,r} \cap R_{q,r} \cap R_r = \emptyset$.  We proceed by contradiction:
let $x \in R_{p,r} \cap R_{q,r} \cap R_r$.
We may assume that $r$ is minimal with respect to this assertion.

Let $y$ be the image of $x$ under the map \mapx[[ Y_r || Y_{r-1} ]].  Then
$y \in R_{p,r-1} \cap R_{q,r-1} \cap C_{r-1}$.  If $q < r-1$, then
for some ruling $R_{r-1} \IN E_{r-1}$, we have
$y \in R_{p,r-1} \cap R_{q,r-1} \cap R_{r-1}$, thereby contradicting the
minimality of $r$.  Hence we may assume that $q = r-1$.

To prove the lemma, it suffices to show that
$\T_y(R_{p,r-1}) + \T_y(R_{q,r-1}) + \T_y(C_{r-1}) = \T_y(Y_{r-1})$.
Since by \pref{lemma1} $R_{p,r-1}$ meets $R_{q,r-1}$ at a unique point,
this will imply that $R_{p,r} \cap R_{q,r} = \emptyset$, thereby yielding
a contradiction.  Substituting $q = r-1$, we must show:
$$\T_y(R_{p,r-1}) + \T_y(R_{r-1}) + \T_y(C_{r-1}) = \T_y(Y_{r-1}).\eqno(*)$%
$The curves $R_{r-1}$ and $C_{r-1}$ meet transversally at $y$, tangentially
spanning $\T_y(E_{r-1})$.  Therefore, to prove $(*)$, and hence the lemma,
it suffices to show that $R_{p,r-1}$ meets $E_{r-1}$ transversally.  This
may be deduced by repeated application of the following two facts, applied
to integers $t$ with $p \leq t \leq r-2$:
\begin{itemize}
\item if $R_{p,t}$ meets $C_t$ transversally (on $Y_t$), then
      $R_{p,t+1}$ meets $E_{t+1}$ transversally (on $Y_{t+1}$);
\item if $R_{p,t}$ meets $E_t$ transversally (on $Y_t$), then
      $R_{p,t}$ meets any smooth curve on $E_t$ transversally (if at all).
\qed
\end{itemize}
\end{proof}

\begin{prop}\label{snooker}
Let $k, m \in \Z$, with $1 \leq k \leq m \leq n$.
Let $R_k \IN E_k$ be a ruling.  Let
$\Gamma = \Gamma_m(R_k)$.  Then $\Gamma$ is a standard labeled graph with
vertices $[k,m] \cap \Z$, and provided that $m < n$,
$\Gamma_{m+1}(R_k)$ is obtained from $\Gamma$ by a single standard operation.
\end{prop}

\begin{proof}
By induction, we may assume that $\Gamma$ is a standard labeled graph with
vertices $[k,m] \cap \Z$.  For each $q$ between $k$ and $m$, let
let $R_q \IN E_q$ be the ruling corresponding to the vertex $q \in \Gamma$.

First we show $(*)$ that if $l$ is such that $k \leq l < m$ and $R_{l,m}$ meets
$C_m$, then in fact $R_{l,m}$, $R_m$, and $C_m$ meet at a common point.
Suppose otherwise: $R_{l,m} \cap R_m \cap C_m = \emptyset$.  We will obtain
a contradiction.  We may choose $m$ to be as small as possible.  There are
two cases.

Case (a).  We have $l = m-1$.  Since $R_{m-1}$ meets $C_{m-1}$ transversally
at a single point, $R_{m-1,m}$ meets $E_m$ at a single point.  Since $\Gamma$
is a standard labeled graph,
it is clear that $R_{m-1,m}$ meets $R_m$.  Since $R_{m-1,m}$ meets $C_m$,
we see that $R_{m-1,m}$ meets $E_m$ at two distinct points: contradiction.
This proves case (a).

Case (b).  We have $l < m-1$.  Since $R_{l,m}$ meets $C_m$
(and a fortiori $R_{l,m}$ meets $E_m$), it follows that $R_{l,m-1}$ meets
$C_{m-1}$.  By the minimality of $m$,
$R_{l,m-1} \cap R_{m-1} \cap C_{m-1} \not= \emptyset$.
It follows that $R_{l,m}$ and $R_{m-1,m}$ meet a common ruling on $E_m$.
Since $R_{m-1,m}$ meets $R_m$, it is clear that this ruling must be $R_m$.
Hence $R_{l,m}$ meets $R_m$.
Thus $R_{l,m}$ meets both $R_m$ and $C_m$, but the three curves do not
meet at a common point.  Hence $R_{l,m}$ meets two distinct rulings on
$E_m$.  Hence
$R_{l,m-1}$ meets $C_{m-1}$ at $\geq 2$ distinct points, so $R_l$ meets
$C_l$ at $\geq 2$ distinct points: contradiction.  This proves case (b),
and hence $(*)$.

We now proceed with the proof of the proposition.  There are two cases.

Case I.  For no $l$ (with $k \leq l < m$) is it true that
$R_{l,m}$, $R_m$ and $C_m$ have a point in common.  We claim
that $\Gamma_{m+1}(R_k) = \Gamma^+$.  It suffices to show that $R_m$ is the
unique curve in $\Gamma$ which meets $C_m$.  This follows from $(*)$.

Case II.  For some $l$ (with $k \leq l < m$), $R_{l,m}$, $R_m$ and $C_m$ have
a point (say $x$) in
common.  We claim that $\Gamma_{m+1}(R_k) = \Gamma^l$.  To prove this, we need
to prove two things:
\begin{romanlist}
\item for any $q$ such that $k \leq q < m$ and $q \not= l$,
      $R_{q,m}$ does not meet $C_m$;
\item $\T_x(R_{l,m}) + \T_x(R_m) + \T_x(C_m) = \T_x(Y_m)$.
\end{romanlist}

The first assertion follows immediately from $(*)$ and from \pref{morsel}.
The second assertion follows from the proof of \pref{morsel}.  \qed
\end{proof}

\begin{lemma}\label{pine-cone}
Let $G$ be a  standard labeled graph, constructed from the single vertex
graph \setof{k} by a sequence of standard operations $+, k^\br{p},\vec o1r$,
for some $p, r \geq 0$, such that $o_1 \not= k$.
Then $G - \setof{k}$ is a standard labeled graph,
which can be constructed from the single vertex graph \setof{k+1} by the
sequence of standard operations:
$$\cases{+^\br{p},\vec o1r,&if $p \geq 1;$\cr
         +, \vec o2r,&if $p=0$ and $r \geq 1;$\cr
         \emptyset,&if $p = r = 0$.}$$
\end{lemma}

The proof of this lemma is left to the reader.

Let $G$ be a standard labeled graph, with smallest vertex $k$, having at
least two vertices.  It is clear that there is a unique $r > k$ such that
$\edge(k,r)$ is in $G$.  We define the {\it order\/} of $G$ to be $r-k$.
Moreover, if $G$ has order $p$, then $G$ is constructed from the single
vertex graph $\setof{k}$ by a sequence of operations which begins with
$+, k^\br{p-1}$, and whose next operation (if any) is not $k$.

Let $G$ be any standard labeled graph, with vertices $k, \ldots, m$.
We associate a function \mp[[ \mu_G || G || \N ]], defined by inducting
on $G$: if $G$ is a single vertex graph, then $\mu_G(k) = 1$.  If $G$
is any standard labeled graph, then $\mu_{G^+}(j) = \mu_G(j)$ and
$\mu_{G^l}(j) = \mu_G(j)$ for all $j$ with $k \leq j \leq m$, and
$\mu_{G^+}(m+1) = \mu_G(m)$, $\mu_{G^l}(m+1) = \mu_G(m) + \mu_G(l)$.
The fact that $\mu_G$ is well-defined follows from \pref{unique-operations}.

\begin{prop}\label{spitup}
Let $G$ be a standard labeled graph, with smallest vertex $k$, having
at least two vertices.  Then:
$$\mu_G = \mu_{\setof{k}} + \sum_{i=1}^{\ord(G)}
          \mu_{G - \setof{k, \ldots, k+i-1}},$%
$where the functions on the \RHS\ are viewed as functions on $G$, via
extension by zero.
\end{prop}

\begin{sketch}
Use \pref{pine-cone}.  If $p = \ord(G)$, then
$$G\ \longleftrightarrow\ +, k^\br{p-1}, *$%
$where $*$ is a sequence of standard operations (possibly empty), not beginning
with $k$.  The case $p = 1$ is left to the reader.  For $p \geq 2$:
$$G - \setof{k}\ \longleftrightarrow\ +^\br{p-1}, *$$
$$G - \setof{k,k+1}\ \longleftrightarrow\ +^\br{p-2}, *$%
$and so forth:
$$G - \setof{k,\ldots,k+p-2}\ \longleftrightarrow\ +, *$%
$$$G - \setof{k,\ldots,k+p-1}\ \longleftrightarrow\ *'$%
$where $*'$ can be determined from \pref{pine-cone}.
We compute $\mu$ in a special case, namely when $*$ is empty.  Then:
$$G\ \longleftrightarrow\ (1,1,2,3,\ldots,p)$$
$$G - \setof{k}\ \longleftrightarrow\ (0^\br{1}, 1^\br{p})$$
$$\cdots$$
$$G - \setof{k,\ldots,k+p-1}\ \longleftrightarrow\ (0^\br{p}, 1^\br{1}),$%
$where the sequences on the right are $(\mu(k), \ldots, \mu(k+p))$.  In this
case the proposition is clear.  The general case is left to the reader.  \qed
\end{sketch}

\begin{corollary}
Fix an integer $k$ with $1 \leq k \leq n$.  Let $R_k \IN E_k$ be a ruling.
If $k = n$ then $[R_{k,n}] = \lbR_n$, and if
$k < n$, then there exists an integer $l$, with $k < l \leq n$, such that
$$[R_{k,n}] = \lbR_k - \sum_{i=k+1}^l \lbR_i.$$
\end{corollary}

\begin{sketch}
For each integer $m$ with $k \leq m \leq n$, let $R_m \IN E_m$ be the
ruling which enters into $\Gamma_n(R_k)$.  For each integer $l$ with
$k \leq l \leq n$, write:
$$\mu_{\Gamma_n(R_l)} = (\vec bln).$%
$By considering the scheme-theoretic inverse image of $R_l$ under the
map \mapx[[ Y_n || Y_l ]], one can show that:
$$\lbR_l = b_l [R_{l,n}] \many+ b_n [R_{n,n}].$%
$The result then follows from \pref{spitup}.  \qed
\end{sketch}

\begin{corollary}
Let $H \IN Y_n$ be a cycle which is a sum of strict transforms of rulings.
Then $[H]$ is a positive $\Z$-linear combination of the classes:
$$\lbR_n, (\lbR_{n-1}-\lbR_n), (\lbR_{n-2}-\lbR_{n-1}-\lbR_n), \ldots,
          (\lbR_1 - \lbR_2 - \cdots - \lbR_n).$$
\end{corollary}

\begin{corollary}\label{snort-snort-snort}
Let $H \IN Y_n$ be a cycle which is a sum of strict transforms of rulings.
Then there exists integers $\vec a1n$ such that
$[H] = a_1 \lbR_1 \many+ a_n \lbR_n$ and for each integer $k$ with
$1 \leq k \leq n$, we have:
$$\left(\sum_{i=1}^{k-1} 2^{k-i-1} a_i\right) + a_k \geq 0.$$
\end{corollary}
\part{Application to \STCIs}

\block{Theorems I and II}\label{section9}

\par\indent\indent
Let $S, T \IN \PP^3$ be surfaces of degrees $s$ and $t$, respectively.
Write $S_0 = S$, $T_0 = T$.
Assume that $S \cap T$ is set-theoretically a smooth curve $C = C_0$.
Let $d = \deg(C)$.  Then $d | st$.  Let $n = st/d$.
Assume that $C \notIN \Sing(S)$ and that $C \notIN \Sing(T)$.
Let $Y_0 = \PP^3$.  Let $Y_1$ be the blow-up of $Y_0$ along $C_0$.  Let
$S_1 \IN Y_1$ be the strict transform of $S$.  There is a unique curve
$C_1 \IN S_1$ which maps isomorphically onto $C_0$.  Let $Y_2$ be the blow-up
of $Y_1$ along $C_1$.  Iterate this process.  This puts us in the situation
of \pref{iteration}.  Let $p_i = p_i(S,C)$, for $i = 1, \ldots, n-1$.

For each $k = 1, \ldots, n$, let $S_k$ and $T_k$ denote the strict transforms
of $S$ and $T$ on $Y_k$.  Since $C \notIN \Sing(S)$ and $C \notIN \Sing(T)$,
it follows that for each $k$ with $0 \leq k \leq n$, $S_k$ meets $T_k$ along
$C_k$ with multiplicity $n - k$.  As consequences of this, we see that
$S_n \cap T_n$ is a union of strict transforms of rulings, and that
$[S_k] = s\lbH - \sum_{i=1}^k \lbE_i$, $[T_k] = t\lbH - \sum_{i=1}^k \lbE_i$.

\label{turnipgreen}
First we derive the formula:
$$\beta_k = ds + (2-4d-2g) - p_k.\eqno(*)$%
$We have $[S_k] \cdot \lbE_k = \lbC_k + p_k\lbR_k$.  (See \ref{goober-peas}.)
Combining this with
$[S_k] = s\lbH - \sum_{i=1}^k \lbE_i$, $\lbH \cdot \lbE_k = d\lbR_k$
(from \ref{iteration}), and $\lbC_k = \lbD_k - \alpha_k\lbR_k$
(from p.\ \pageref{quirkC}), we obtain:
$$ds\lbR_k - \sum_{i=1}^k (\lbE_i \cdot \lbE_k) = \lbD_k
     - \alpha_k\lbR_k + p_k\lbR_k.$%
$Combine this with $\lbE_i \cdot \lbE_k = -\beta_i\lbR_k$ (if $i < k$)
(from \ref{iteration}) and $\lbE_k^2 = -\lbD_k$ (from p.\ \pageref{quirkC})
to obtain:
$$ds\lbR_k + (\svec\beta1{k-1})\lbR_k = -\alpha_k\lbR_k + p_k\lbR_k.$%
$Combine this with the formula:
$$\alpha_k = (2-4d-2g) - (\svec\beta1k)\eqno(**)$%
$from p.\ \pageref{quirkalpha}, to obtain $(*)$.

Since $S_n \cap T_n$ is a union of strict transforms of rulings, it follows
from \pref{snort-snort-snort} that:
$$(s\lbH - \sum_{i=1}^n \lbE_i)(t\lbH - \sum_{i=1}^n \lbE_i)
     = \sum_{l=1}^n a_l \lbR_l, \eqno(*{*}*)$%
$for some integers $a_l$ such that for each $k$ with $1 \leq k \leq n$, we
have:
$$\left(\sum_{m=1}^{k-1} 2^{k-m-1} a_m\right) + a_k \geq 0.$%
$
We proceed to analyze the consequences of this.  The left hand side of
$(*{*}*)$ equals:

$$ st\lbH^2 - d(s+t)(\sum_{i=1}^n \lbR_i)
           - \sum_{1 \leq i < j \leq n} \beta_i \lbR_j
           - \sum_{1 \leq j < i \leq n} \beta_j \lbR_i$$
$$         - \sum_{k=1}^n \left[ d\lbH^2
                + \left( \sum_{i=1}^{k-1} \beta_i \lbR_i \right)
                + \alpha_{k-1} \lbR_k \right].$%
$
Then for each $m$ with $1 \leq m \leq n$:
$$-a_m = d(s+t) + 2\sum_{i=1}^{m-1} \beta_i + (n-m) \beta_m + \alpha_{m-1}.$%
$Substituting
$\alpha_k = (2-4d-2g) - (\beta_1 \many+ \beta_k)$, we obtain:
$$-a_m = d(s+t) + \sum_{i=1}^{m-1} \beta_i + (n-m)\beta_m + (2-4d-2g).$%
$
\par\indent In the special case where
$\Sing(S) \cap \Sing(T) = \emptyset$, we have
$a_m = 0$ for all $m$, with $1 \leq m \leq n$.
Now substitute $\beta_i = ds + (2-4d-2g) - p_i$.  We obtain:
$$\left(\sum_{i=1}^{m-1} p_i\right) + (n-m)p_m = d[n(s-4) + t] + (2-2g)n.$%
$This implies:

\begin{theorem}\label{thmI}
{\bf (``I'')}
Let $C \IN \P3$ be a smooth curve of degree $d$ and genus $g$.  Suppose that
$C = S \cap T$, where $S$ and $T$ are surfaces of degree $s$ and $t$
respectively.  Assume that $\Sing(S) \cap \Sing(T) = \emptyset$.
Let $n = st/d$.
Let $p_i = p_i(S,C)$, for each $i = 1, \ldots, n-1$, as defined in
\S\ref{measure}.  Then:
$$p_1 = \cdots = p_{n-1} = {1 \over n-1}
        \left\{d[n(s-4)+t] + (2-2g)n\right\}.$$
\end{theorem}

\begin{example}\label{examplex-char-two}
If $s = t = 4$, $d = 4$, $g = 0$, we obtain $p_1 = p_2 = p_3 = 8$.  This
can occur in characteristic two, at least.  Indeed, let $(S,C)$ be as in
\pref{example-char-two}, and let $T$ be given by $z^4-xw^3=0$.
\end{example}

We now return to the general case.  For each $k$ with $1 \leq k \leq n$,
we have:
\splitdisplay{\left(\sum_{m=1}^{k-1} 2^{k-m-1}
[d(s+t) + \sum_{i=1}^{m-1} \beta_i + (n-m) \beta_m + (2-4d-2g)] \right)%
}{+ [d(s+t) + \sum_{i=1}^{k-1} \beta_i + (n-k) \beta_k +(2-4d-2g)]\leq 0.%
}A simplification yields:
$$2^{k-1}[d(s+t-4)+2-2g]
   + \left( \sum_{i=1}^{k-1} 2^{k-i-1}(n-i+1) \beta_i \right)
   + (n-k) \beta_k \leq 0.$%
$Note that:
$$\sum_{i=1}^{k-1} 2^{k-i-1}(n-i+1) \ =\ (n-1)2^{k-1} + k - n.$%
$Substitute $\beta_i = ds + (2-4d-2g) - p_i$.  We obtain:

\begin{theorem}\label{thmII}
{\bf (``II'')}
Let $C \IN \P3$ be a smooth curve of degree $d$ and genus $g$.  Suppose that
$C = S \cap T$, where $S$ and $T$ are surfaces of degree $s$ and $t$
respectively.  Assume that $C \notIN \Sing(S)$ and $C \notIN \Sing(T)$.
Let $n = st/d$.
Let $p_i = p_i(S,C)$, for each $i = 1, \ldots, n-1$, as defined in
\S\ref{measure}.  Then for each $k = 1, \ldots, n-1$, we have:
$$\sum_{i=1}^{k-1} 2^{k-i-1}(n-i+1) p_i + (n-k)p_k
\geq 2^{k-1}\left\{dt + n[d(s-4)+2-2g]\right\}.$$
\end{theorem}

\begin{examples}
\
\begin{itemize}
\item $s = 2$, $t = 3$, $d = 3$, $g = 0$: then the theorem yields the
      single inequality $p_1 \geq 1$;
\item $s = 2$, $t = 2$, $d = 1$, $g = 0$: as above the theorem yields
      $p_1 \geq 1$;
\item $s = 4$, $t = 4$, $d = 4$, $g = 0$: the theorem yields three
inequalities:
      \begin{romanlist}
      \item $p_1 \geq 8$;
      \item $2p_1 + p_2 \geq 24$;
      \item $8p_1 + 3p_2 + p_3 \geq 96$.
      \end{romanlist}
\end{itemize}
\end{examples}
\block{Theorems III, Q, and B}

\begin{theorem}\label{thmIII}
{\bf (``III'')}
Let $C \IN \C\P3$ be a smooth curve of degree $d$ and genus $g$ which is
the \STCI\ of two surfaces $S$, $T$ of degrees $s$, $t$, respectively.
Let $n = st/d$.  Let $p_k = p_k(S,C)$, for each $k \in \N$.  Assume that
$S$ has only rational singularities.  Assume that $C \notIN \Sing(T)$.  Then:
$$\sum_{k=1}^\infty{1 \over k(k+1)}p_k \geq {d^2 \over s} + d(s-4) + 2-2g.$$
\end{theorem}

\begin{proof}
By (\Lcitemark 18\Rcitemark \ 1.1), we know that
$\Delta(S,C) = d^2/s + d(s-4) + 2-2g$.  Apply \pref{rdp-one}.  \qed
\end{proof}

\begin{remark}
In the statement of \pref{thmIII}, we do not know if $\sum_{k=1}^\infty$
can be replaced by $\sum_{k=1}^{n-1}$.  From (\ref{interesting}),
we see that this can at least be done if $(S,C)$ does not contain any
singularities of type $D_{n,n}$ (with $n$ odd).
\end{remark}

\begin{lemma}\label{bound-formula}
Let $S \IN \C\P3$ be a surface of degree $s$ having only rational
singularities.
Let \mp[[ \pi || \tS || S ]] be a minimal resolution.  Let $N$ be the number
of exceptional curves on $\tS$.  Then:
$$N \leq {s \over 3}(2s^2 - 6s + 7) - 1.$$
\end{lemma}

\begin{proof}
Clearly $N \leq \rank \NS(\tS) - 1$.  Also $\rank \NS(\tS) \leq h^{1,1}(\tS)$,
so it suffices to show that $h^{1,1}(\tS) = {s \over 3}(2s^2-6s+7)$.  By
simultaneous resolution of rational double points
\Lcitemark 2\Rcitemark \Rspace{}, and deformation invariance of
Hodge numbers, we may reduce to showing that
$h^{1,1}(S) = {s \over 3}(2s^2-6s+7)$ if $S$ is itself nonsingular.  We have:
\begin{eqnarray*}
h^{1,1}(S) & = & h^2(S,\Q) - 2h^2(S,\O_S) \\
           & = & [\chi_{\op top}(S) + 4h^1(S,\O_S) - 2] - 2h^2(S,\O_S).
\end{eqnarray*}
Using the fact that the top Chern class of the tangent bundle equals the Euler
characteristic (see e.g.{\ }\Lcitemark 1\Rcitemark \ 11.24, 20.10.6), and using
Riemann-Roch, we find:
$$\chi_{\op top}(S) = c_2(S) = 12\chi(S) - c_1^2(S)
  = 12(1-h^1(S,\O_S)+h^2(S,\O_S)) - (4-s)^2s.$%
$The formula for $h^{1,1}$ follows from $h^1(S,\O_S) = 0$ and
\formulaqed{h^2(S,\O_S) = {s-1 \choose 3}.}
\end{proof}

\begin{remark}
We do not know if the lemma remains valid if $\C$ is replaced by an
algebraically closed field of positive characteristic.
\end{remark}

\begin{lemma}\label{bungobungo}
Let $(p_k)_{k \in \N}$ be a sequence of nonnegative integers.  Let $n$
be a nonnegative integer.  Assume that:
\begin{romanlist}
\item $p_1 \leq 9 - {2\over5}n$
\item $p_1 \geq p_2 \geq p_3 \geq \cdots$
\item $\sum_{k=1}^\infty p_k \leq 19 - n$
\item $n/4 + \sum_{k=1}^\infty {1 \over k(k+1)} p_k \geq 6$.
\end{romanlist}
\par\noindent Then $n=0$ and
$(p_k) \in \setof{ (9,8,2), (9,9), (9,9,1) }.$
\end{lemma}

\begin{proof}
Constraints (i), (iii), and (iv) imply that
$$n/4 + \hbox{$1 \over 2$}\floor{9-\hbox{$2 \over 5$}n}
      + \hbox{$1 \over 6$}(19-n-\floor{9-\hbox{$2 \over 5$}n} )
        \geq 6.$%
$It follows that $n \in \setof{0,2}$.

Suppose that $n=2$.  Then the same constraints imply that
$$1/2 + \hbox{$1 \over 2$}p_1 + \hbox{$1 \over 6$}(17-p_1)
  \geq 6,$%
$so $p_1 = 8$.  Now we see that the \LHS\ of (iv) is maximized when
$(p_k) = (8,8,1)$.  In that case, the \LHS\ of (iv)
is $5{11\over12}$: contradiction.

Hence $n=0$.  Then
$$\hbox{$1 \over 2$}p_1 + \hbox{$1 \over 6$}(19-p_1) \geq 6,$%
$so $p_1 = 9$.  If $p_2 \leq 7$,
then the sum in (iv) is bounded by the sum obtained when
$(p_k) = (9,7,3)$.  This sum is $< 6$, so $p_2 \not\leq 7$.  Hence
$p_2 \in \setof{8,9}$.  Etc.  \qed
\end{proof}

\begin{prop}\label{kformula}
Let $C \IN \C\P3$ be a smooth curve of degree $d$ and genus $g$, which
lies on a surface $S \IN \C\P3$ of degree $s$.  Assume that $C \notIN
\Sing(S)$.
Let $p_1 = p_1(S,C)$.  Let $N$ be the normal bundle of $C$ in $\CP^3$, and
let $l$ be the maximum degree of a sub-line-bundle
of $N$.  Let $k = 3d + (2g-2) - l$.  Then $p_1 \leq d(s-1) - k$.
\end{prop}

\begin{proof}
We use the notation of \pref{iteration}.  We also use various facts from
\S\ref{section9}, which although apparently dependent on another surface $T$,
actually make sense in this context.
We have $\inn{\lbC_1, \lbC_1} \geq \deg(N) - 2l$.
Since $\deg(N) = 4d+2g-2$ and $k = 3d + (2g-2) - l$, we have:
$$\inn{\lbC_1,\lbC_1} \geq -2d+2-2g+2k.\eqno(\dag)$%
$Now $\lbC_1 = \lbD_1 - \alpha_1\lbR_1$ in $A^1(E_1)$, and
$\inn{\lbD_1,\lbD_1} = 2-2g-4d$, so:
$$\inn{\lbC_1,\lbC_1}\ =\ \inn{\lbD_1,\lbD_1} - 2\alpha_1\ =
     \ 2 - 2g - 4d - 2\alpha_1.$%
$Combining this with $(\dag)$, we obtain $d + \alpha_1 \leq -k$.  The
formulas $(*)$ and $(**)$ from \S\ref{section9} imply that
$\alpha_1 = p_1 - ds$.  Hence $p_1 \leq d(s-1)-k$.  \qed
\end{proof}

\begin{remark}
This result \pref{kformula} is a strengthening of the very elementary fact
that:
$$\abs{\Sing(S) \cap C} \leq d(s-1).$$
\end{remark}

\begin{theorem}\label{thmQ}
{\bf (``Q'')}
Let $C \IN \CP^3$ be a curve.  Assume that $C = S \cap T$
set-theoretically for some surfaces $S$ and $T$.  Assume that $S$ is normal.
Assume that $\deg(C) > \deg(S)$.  Then $C$ is linearly normal.
\end{theorem}

\begin{proof}
To any Weil divisor $E$ on a normal surface $S$, one can associate a reflexive
$\O_S$-module $\O_S(E)$.  We recall the following result of Sakai from
\Lcitemark 29\Rcitemark \Rspace{}, which is a slightly less general version of
theorem
5.1 of that paper:
\begin{quote}
{\it Let $S$ be a normal projective surface.  Let $D$ be a nef Weil divisor on
$S$ with $D^2 > 0$.  Then $H^1(S, \O_S(-D)) = 0$.}
\end{quote}
Since the canonical map \mapx[[ H^0(\P3, \O_{\P3}(1)) || H^0(S, \O_S(1)) ]] is
surjective, it suffices to show that the canonical map
\mapx[[ H^0(S, \O_S(1)) || H^0(S, \O_C(1)) ]] is surjective.  Let $H$ be
a hyperplane section of $S$.  From the long exact sequence coming from
\sescomma{\O_S(H-C)}{\O_S(H)}{\O_C(H)%
}we see that it is sufficient to show that $H^1(S, \O_S(H-C)) = 0$.

Let $d = \deg(C)$.  Let $s = \deg(S)$, $t = \deg(T)$, and let $n$ be the
multiplicity of intersection of $S$ with $T$ along $C$, $n = st/d$.  Since
$d > s$, we have $t > n$.  Hence $(t-n)H$ is a very ample Cartier
divisor.  Since $n(C-H) \sim (t-n)H$, the theorem follows from Sakai's
result.  \qed
\end{proof}

\begin{corollary}
Let $C \IN \CP^3$ be a smooth curve.  Assume that $C$ is the set theoretic
complete intersection of two normal surfaces $S$ and $T$, with multiplicity
$\leq 3$.  Then $C$ is linearly normal.
\end{corollary}

\begin{proof}
Using the notation of the proof of \pref{thmQ}, we are done
if either $s$ or $t$ is bigger than $n$.  Otherwise, $d \leq 3$, and so
$C$ is linearly normal anyway.  \qed
\end{proof}

\begin{remark}
For the case of multiplicity $4$, we must have $C$ linearly normal, except
possibly for the case where $C$ is a rational quartic, which is the \STCI\ of
two normal quartic surfaces.  It is not known if this is possible.
\end{remark}

\begin{theorem}\label{thmB}
{\bf (``B'')}
Let $S, T \IN \C\P3$ be surfaces.  Assume that $S \cap T$ is set-theoretically
a smooth curve.  Assume that $\deg(S) = 4$ and that $S$ has only rational
singularities.  Then $C$ is linearly normal.
\end{theorem}

\begin{proof}
Let $C$ have degree $d$ and genus $g$.  By \pref{thmQ}, \WMAT\ $d=4$ and $g=0$.
By\Lspace \Lcitemark 18\Rcitemark \Rspace{}, \WMAT\ $\deg(T) \geq 4$.
Since $\deg(S) \leq \deg(T)$, we may assume that $C \notIN \Sing(T)$, as
follows.  Suppose that $C \IN \Sing(T)$.
Write $S = V(f)$, $T = V(g)$.  Choose $h$ so that
$\deg(fh) = \deg(g)$, and so that $C \notIN V(h)$.  Then
$C \notIN \Sing(V(fh + g))$.  Hence we may replace $T$ by $V(fh + g)$.

Write $(S,C) = (S',C') + (S'',C'')$, where
$$(S'',C'') = D_{n_1,n_1} \many+ D_{n_r,n_r},$%
$$\vec n1r$ are odd integers $\geq 5$, and $(S',C')$ is a configuration
which does not involve any such singularities.  Let $p_i = p_i(S',C')$.
Let $n = \svec n1r$.  We show that the hypotheses
of \pref{bungobungo} are satisfied.

We apply \pref{kformula}, using that fact\Lspace \Lcitemark 5\Rcitemark
\Rspace{} that
$l=7$, concluding that $p_1(S,C) \leq 9$.  We have
$p_1(\sum D_{n_i,n_i}) = \sum (n_i-1) / 2)$ by \pref{fantastico}, and
$n_i - 1 \geq {4\over5} n_i$, so $p_1(\sum D_{n_i,n_i}) \geq {2\over5}n$.
Thus hypothesis (i) is satisfied.  Hypothesis (ii) holds.  Hypothesis
(iii) follows from \pref{bound-formula} and from the fact that $(S',C')$
contains no $D_{m,m}$ pairs with $m$ odd, $m \geq 5$, so that
$\DEF(S',C') \geq 0$.  To prove hypothesis (iv), we would like to use
(\ref{thmIII} = ``III''), but that is not quite good enough.  By
\pref{rdp-one},
$$\sum_{k=1}^\infty {1\over k(k+1)}p_k \geq \Delta(S',C').$%
$Let $s = \deg(S) = 4$.  By (\Lcitemark 18\Rcitemark \ 1.1), we know
that:
$$\Delta(S,C) = d^2/s + d(s-4) + 2 - 2g = 6.$%
$Then $\Delta(S',C') = \Delta(S,C) - \Delta(S'',C'')$.  By
\pref{delta-formulas}, $\Delta(S'',C'') = n/4$.  Hypothesis (iv) follows.

By \pref{bungobungo}, we conclude that $n = 0$ and that:
$$\type(S,C) \in \setof{ (9,8,2), (9,9), (9,9,1) }.$%
$We will use \pref{fantastico}, \pref{potato-1}, and \pref{potato-2}.

Suppose that $\type(S,C) = (9,9,1)$.  Then for some $p \in C$,
$\type(S,C)_p  = (r,1,1)$ for some $r \geq 1$.  The case $r > 1 $ is
impossible because $\type(S,C) - \type(S,C)_p = (9-r,9-1)$ and
$9-r \geq 9-1$.  Hence $\type(S,C)_p = (1,1,1)$.  Hence $(S,C)_p = A_{3,1}$.
Since $\Sigma(S,C) \leq 19$ by \pref{bound-formula}, and since
$9+9+1 = 19$, we have $\DEF(S,C) = 0$.  Therefore,
since the other singularities of $S$ along $C$ have type
$(k,k)$ for some $k \leq 8$, we see that the other singularities must be
$A_{3k-1,k}$ for some $k \in \setof{1, \ldots, 8}$, depending on the
singular point.  Amongst these, only $A_{2,1}$ has deficiency zero.
Hence $(S,C) = 8A_{2,1} + A_{3,1}$.  Hence
$\Delta(S,C) = 8({2\over3}) + {3\over4} \not= 6$: contradiction.

Now suppose that $\type(S,C) = (9,8,2)$.  Then for some $p \in C$,
$$\type(S,C)_p \in \setof{ (1,1,1), (2,1,1), (2,2,2) }.$%
$These types are realized by the singularities $A_{3,1}$, $A_{4,2}$,
and $A_{7,2}$, respectively, and by no others.  Since $A_{7,2}$ has nonzero
deficiency, it can be excluded.  Hence $(S,C)_p \in \setof{A_{3,1},A_{4,2}}$.
If $(S,C)_p = A_{4,2}$, we find (by analogy with the $(9,9,1)$ case) that
$(S,C) = 6A_{2,1} + A_{3,1} + A_{4,2}$.  Hence $\Delta(S,C) =
6({2\over3}) + {3\over4} + {6\over5} \not= 6$: contradiction.  If
$(S,C)_p = A_{3,1}$, then \WMAT\ $(S,C) = 2A_{3,1} + \hbox{other}$, where
the ``other'' part must have
type $(7,6)$.  The only zero-deficiency \RDP\ configuration which realizes
this type is $A_{1,1} + 6A_{2,1}$.  Hence
$(S,C) = A_{1,1} + 6A_{2,1} + 2A_{3,1}$.  By\Lspace \Lcitemark 24\Rcitemark
\Rspace{},
the sum of the contributions of the singularities must not exceed
$(2/3)\deg(S)(\deg(S)-1)^2 = 24$, where each singularity $p$ contributes
$e(E) - 1/\abs{G}$, $e(E)$ is the topological Euler characteristic of the
exceptional fiber in the minimal resolution of $p$, and $G$ is the order of the
group which defines $p$ as a quotient singularity.  In particular, an
$A_n$ singularity contributes $(n+1) - (n+1)^{-1}$.  Then the sum of the
contributions is $25$: contradiction.

Suppose that $\type(S,C) = (9,9)$.  Then each singularity of $S$ along $C$
must have type $(k,k)$ for some $k$, depending on the singular point.  Hence
$(S,C)$ must be built up from $E_{6,1}$ and $A_{3k-1,k}$ for various $k$.
Since $\DEF(E_{6,1}) = 2$, we may rule out that case.  In fact, there are
only two configurations with deficiency $\leq 1$: either
$(S,C) = 9A_{2,1}$ or else $(S,C) = 7A_{2,1} + A_{5,2}$.  In both cases,
$\order(S,C) = 3$.  Hence we may assume that $\deg(T) = 3$.  By
\Lcitemark 18\Rcitemark \Rspace{}, we know that this is impossible.  \qed
\end{proof}
\block{Theorem A}

\begin{lemma}\label{murky-algebra}
Let $s,t,d,g \in \Z$.  Assume that $t \geq s \geq 4$, $d \geq 1$,
and $g \geq 0$.  Assume that $d \vert st$.  Let $n = st/d$.  Assume
that $n \geq 2$.  Let $r = d[n(s-4)+t] + (2-2g)n$.  Assume that
$r \leq {s \over 3}(2s^2-6s+7)-1$.  Then $d \leq g + 3$.
\end{lemma}

\begin{proof}
We assume that $d \geq g + 4$, working toward a contradiction.  We have
$2 - 2g \geq 10 - 2d$, so:
$$d[n(s-4) + t] + (10-2d)n \leq {s\over3}(2s^2-6s+7) - 1.\eqno(*)$%
$
\par First suppose that $s = 4$.  Then $t \geq 4$, $t \geq d/2$, and
$dt + (10-2d)n \leq 19$.  Substituting $n = st/d = 4t/d$ and simplifying,
we obtain:
$$(d^2 - 8d + 40)t \leq 19d.\eqno(\dag)$%
$Since $t \geq d/2$, we have $(d^2 - 8d + 40)(1/2) \leq 19$.  It follows
that $d \leq 7$.  Hence $d \in \setof{4,5,6,7}$.  In each case, $(\dag)$
gives us an upper bound $t_{\rm max}$ for $t$:

{\renewcommand{\arraystretch}{1.0}
\vspace*{0.1in}
\centerline{
\begin{tabular}{||c|c||}  \hline
$d$ & $t_{\rm max}$\\ \hline
$4$ & $3$\\ \hline
$5$ & $3$\\ \hline
$6$ & $4$\\ \hline
$7$ & $4$\\ \hline
\end{tabular}    }
\vspace*{0.1in}}

The cases $d \in \setof{4,5}$ contradict $t \geq 4$.  In case
$d \in \setof{6,7}$, we have $t = 4$, which contradicts our assumption that
$d | st$.  Hence $s > 4$.

Now suppose that $s = 5$.  Then $t \geq 5$, $t \geq 2d/5$, and
$d(t+n) + (10-2d)n \leq 44$.  Substituting $n = st/d = 5t/d$ and simplifying,
we obtain:
$$t \leq 44d / (d^2 - 5d + 50).$%
$This implies that $t < 5$: contradiction.  Hence $s \geq 6$.

Since $n = st/d \geq s^2/d$, it follows from $(*)$ that:
$$d[{s^2 \over d}(s-6) + s] + 10{s^2 \over d} \leq {s \over 3}(2s^2-6s+7).$%
$This implies that:
$$s(s-6) + d + 10{s \over d} \leq {1 \over 3}(2s^2 - 6s + 7)\eqno(**)$%
$and in particular that:
$$s(s-6) \leq {1 \over 3}(2s^2 - 6s + 7).$%
$It follows that $s \leq 12$.  Hence $6 \leq s \leq 12$.  If $s = 12$,
$(**)$ implies that $d + 120/d \leq 2{1\over3}$.  This is absurd.  In a similar
manner, one may eliminate the cases where $6 \leq s \leq 11$.  \qed
\end{proof}

\begin{theorem}\label{thmA}
{\bf (``A'')}
Let $S, T \IN \C\PP^3$ be surfaces.
Assume that $S$ has only rational singularities.
Assume that $\deg(S) \leq \deg(T)$.
Assume that $S \cap T$ is set-theoretically a smooth curve $C$ of
degree $d$ and genus $g$.
Assume that $\Sing(S) \cap \Sing(T) = \emptyset$.  Then $d \leq g + 3$.
\end{theorem}

\begin{proof}
Let $s = \deg(S)$, $t = \deg(T)$, $n = st/d$.
Let $p_i = p_i(S,C)$.  By (\ref{thmI} = ``I''), we have:
$$p_1 = \cdots = p_{n-1} = {1 \over n-1}
        \left\{d[n(s-4)+t] + (2-2g)n\right\}.$$

We show that $S$ has no singularities of type $D_{t,t}$
(with $t$ odd, $t \geq 5$), lying on $C$.  There are two cases.  If $n = 2$,
then $\order(S,C) | 2$.  But by \pref{order-calc}, the order of $D_{t,t}$
(as above) is $4$.  Hence $n > 2$.  Hence $p_1 = p_2$.  But
$p_1(D_{t,t}) > p_2(D_{t,t})$ (by \ref{fantastico}), so
``$D_{t,t} \notin (S,C)$'' for $t$ odd, $t \geq 5$.

{}From this, it follows that $\svec p1{n-1} \leq \Sigma(S,C)$.  Let
$r = d[n(s-4)+t] + (2-2g)n$.  By \pref{bound-formula}, we conclude that
$r \leq {s\over3}(2s^2-6s+7)-1$.

The case $n = 1$ corresponds to a complete intersection, and the theorem is
easily verified in this case.  Therefore \WMAT\ $n \geq 2$.
By\Lspace \Lcitemark 18\Rcitemark \Rspace{}, it follows that if $s \leq 3$,
then
$d \leq g+3$.  Hence we may assume that $s \geq 4$.  Therefore
\pref{murky-algebra} applies, and we conclude that $d \leq g + 3$.  \qed
\end{proof}

\begin{corollary}
Let $S, T \IN \C\PP^3$ be surfaces.
Assume that $S$ and $T$ have only rational singularities.
Assume that $S \cap T$ is set-theoretically a smooth curve $C$ of
degree $d$ and genus $g$.
Assume that $\Sing(S) \cap \Sing(T) = \emptyset$.  Then $d \leq g + 3$.
\end{corollary}
\block{Theorem X}

\par\indent\indent
As a corollary of theorem (I), we show:

\begin{theorem}\label{thmX}
{\bf (``X'')}
Let $C \IN \P3$ be a smooth curve.  Assume that $C$ is not a complete
intersection.  Suppose that $C = S \cap T$ as sets, where $S, T \IN \P3$ are
surfaces.  Assume that $\Sing(S) \cap \Sing(T) = \emptyset$.  Then:
$$\deg(S), \deg(T) < 2 \cdot \deg(C)^4.$$
\end{theorem}

\par\noindent First we make a few remarks.

\begin{arabiclist}

\item The proof of theorem (X) depends primarily on the fact that the numbers
$p_k$ in theorem (I) must be integers.

\item The importance of theorem (X) is that an upper bound
is given for the degrees of $S$ and $T$, that this bound is computable,
and that this bound depends only on
the degree of $C$.  In the proof, we give the better bounds
$\deg(S) < 2 \cdot \deg(C)^2$, $\deg(T) < 2 \cdot \deg(C)^4$, provided
that $\deg(S) \leq \deg(T)$.

\item Via the bounds in theorem (X), it becomes a computer triviality to
find all possible degrees for $S$ and $T$ which are consistent with the
integrality of the numbers $p_k$ in theorem (I).

\item Doing this when $\deg(C) = 4$, $\genus(C) = 0$, and assuming for
efficiency that $\deg(S) \leq \deg(T)$, we find:
$$(\deg(S), \deg(T)) \in \{ (3,4), (3,8), (4,4), (4,7), (6,26), (9,48),
     (10,28),$$
$$(12,18), (13,16), (17,220), (18,118), (19,84), (20,67), (22,50), (28,33)
\}.$$
[We have excluded the cases where $\deg(S)$ is $1$ or $2$, which cannot occur.]

\item We do not know which of these pairs of integers can be realized
by pairs of surfaces, as in theorem (X).  All that we know is that
$(3,4)$ and $(3,8)$ cannot be realized in characteristic zero, and that
$(4,4)$ can be realized in characteristic two.

\item Theorem (X) is false without the hypothesis that $C$ is
a complete intersection.  Counterexample: for any $s \in \N$,
one can find a smooth curve $D \IN \P2$ of degree $s$ and a line $L \IN \P2$
such that $D \cap L$ is a single point, set-theoretically.  Let $S$ and $T$ be
cones over $D$ and $L$, with the same vertex.  Then $S \cap T$ is a line,
set-theoretically.  Theorem (X) is also false without the hypothesis that
$\Sing(S) \cap \Sing(T) = \emptyset$.

\end{arabiclist}

Before proceeding with the proof of theorem (X), we need the following
lemma, which was known in characteristic zero, and for
the smooth case, was known in all characteristics.  (See proof for
references.)

\begin{lemma}\label{torsion-free}
Let $S \IN \P3$ be a normal surface.  Then $\Pic(S)/\Pic(\P3)$ is
torsion-free.
\end{lemma}

Before proceeding with the proof, we recall some standard material on
differentials for which we do not have a good reference.  First of all,
for any scheme $X$, there is a map of sheaves of abelian groups:
\dmap[[ \dlog || \O_X^* || \Omega_X ]]%
given by $f \mapsto df/f$.   (All sheaves we shall discuss are sheaves
on the Zariski site.)

Now suppose that $X$ is a normal proper variety, defined
over an algebraically closed field $k$ of positive characteristic $p$.  Then we
have an exact sequence:
\diagramx{0&\mapE{}&\O_X^*&\mapE{F}&\O_X^*&\mapE{\dlog}&\Omega_X\cr%
}of sheaves of abelian groups on $X$, where $F$ denotes the Frobenius map.
The exactness in the middle depends on normality, and may be deduced e.g.\ from
(\Lcitemark 16\Rcitemark \ I\ 4.2).  Let $\shD$ be the image of
$\dlog$.  Since $X$ is
proper, $H^0(X, \O_X^*) = k^*$, so $H^0(X, F)$ is an isomorphism, and we obtain
an isomorphism $H^0(X, \shD) \iso \ker H^1(X, F)$.  We have
$\ker H^1(X, F) \iso {}_p\Pic(X)$.  Composing with the canonical injection
\mapx[[ H^0(X, \shD) || H^0(X, \Omega_X) ]], we obtain an injective group
homomorphism:
\dmap[[ \psi_X || {}_p\Pic(X) || H^0(X, \Omega_X). ]]%

\begin{proofnodot}
(of \ref{torsion-free}).
First we show $(*)$ that $\Pic(S)/\Pic(\P3)$ has no torsion, except possibly
for $p$-torsion, when the ground field has positive characteristic $p$.
These arguments are very similar to those given by
Lang\Lspace \Lcitemark 21\Rcitemark \Rspace{}.  The methods were invented by
Grothendieck
(\Lcitemark 10\Rcitemark \ Expos\'e  XI), and further studied
by Hartshorne (\Lcitemark 13\Rcitemark \ \S4.3).  We refer the
reader to\Lspace \Lcitemark 21\Rcitemark \Rspace{} or\Lspace \Lcitemark
13\Rcitemark \Rspace{}
for details.

Let $S_n$ be the \th{n} infinitesimal neighborhood of $S$ in $\P3$.  Then:
$$\Pic(\P3) \iso \inverselim \Pic(S_n).$%
$Moreover, for each $n$ there is an exact sequence of abelian groups:
\les{\Pic(S_{n+1})}{\Pic(S_n)}{H^2(S, \shJ^n/\shJ^{n+1}),%
}where $\shJ$ is the ideal sheaf of $S$ in $\P3$.  Since the $H^2$ term
is a vector space, $(*)$ follows.

{}From now on we may assume that the ground field has positive characteristic
$p$.  A standard calculation shows that $H^0(S, \Omega_S) = 0$.
Up to now, we have not used the hypothesis that $S$ is normal.  We
now use this hypothesis.  Via the map $\psi_S$, defined immediately above
this proof, we see that ${}_p\Pic(S) = 0$.  (This argument is essentially
that used in\Lspace \Lcitemark 21\Rcitemark \Rspace{}.)

Finally, to complete the proof, we must show that $[O_S(1)]$ does not
have a \th{p} root in $\Pic(S)$.  The argument given here is essentially
the argument given in (\Lcitemark 4\Rcitemark \ 1.8).  For any variety $X$,
there is a natural group homomorphism
\mp[[ H^1(\dlog) || \Pic(X) || H^1(X, \Omega_X) ]].
Consider this map when $X = S$ and when $X = \P3$.  A standard calculation
shows that the map \mapx[[ H^1(\P3, \Omega_{\P3}) || H^1(S, \Omega_S) ]] is
injective.  Moreover, one knows that the image of $[\O_{\P3}(1)]$ in
$H^1(\P3, \Omega_{\P3})$ is not zero.
(See e.g.{\ }\Lcitemark 14\Rcitemark \ Chapter 3, exercise 7.4.)
Hence the image of $[\O_S(1)]$ in
$H^1(S, \Omega_S)$ is not zero.  Hence $[O_S(1)]$ does not have a \th{p}
root in $\Pic(S)$.  \qed
\end{proofnodot}

\begin{remark}
Over an algebraically closed field of positive characteristic, let
$S \IN \P3$ be a surface, not necessarily normal.  We do not know if
$\Pic(S)/\Pic(\P3)$ is torsion-free, or even if $\Pic(S)$ is torsion-free.
Answers to these questions might be obtained from a general structure
theorem for \snap\hbox{$\Ker[\Pic(S)\ \mapE{}\ \Pic(\nor{S})]$}, where $S$ is
an arbitrary projective variety.
\end{remark}

\begin{corollary}\label{smooth-stci}
In $\P3$, suppose that $C = S \cap T$ as sets, where $C$ is a curve,
and $S, T$ are surfaces.  Assume that $C$ does not meet $\Sing(S)$.  Then
there exists a surface $T' \IN \P3$ such that $C = S \cap T'$,
scheme-theoretically.
\end{corollary}

\begin{proof}
Since $T \cap \Sing(S) = \emptyset$, $S$ is normal.  By \pref{torsion-free},
$\Pic(S)/\Pic(\P3)$ is torsion-free.  Hence $[\O_S(C)] = 0$ in
$\Pic(S)/\Pic(\P3)$.  Hence $\O_S(C) \iso \O_S(t)$ for some $t \in \N$.
Since the canonical map \mapx[[ H^0(\P3, \O_{\P3}(t)) || H^0(S, \O_S(t)) ]]
is surjective, it follows that there exists a surface $T'$ of degree $t$ as
claimed.  \qed
\end{proof}

\begin{remark}
Robbiano\Lspace \Lcitemark 28\Rcitemark \Rspace{} proved this in the case
where $S$ is smooth and the ground field has characteristic zero.
\end{remark}

\begin{proofnodot}
(of theorem X).
Let $d = \deg(C)$, $g = \genus(C)$, $s = \deg(S)$, $t = \deg(T)$.  We may
assume that $s \leq t$.  We show that $s < 2d^2$ and $t < 2d^4$.  Let
$n = st/d$.  By theorem (I), we know that:
$$(n-1) \kern3pt | \kern3pt \setof{ d [ n (s-4) + t ] + (2-2g)n }.$%
$A proof of this fact, independent of (I), is given at the end of this paper.
By \pref{smooth-stci}, we know that $C$ meets $\Sing(S)$.  Hence
$p_1(S,C) > 0$.  Hence the \RHS\ is positive.
Write $d = d_s d_t$, where $d_s, d_t \in \N$, $d_s|s$, and $d_t|t$.  Let
$s_1 = s/d_s$, $t_1 = t/d_t$.  Then $n = s_1t_1$, so:
$$(s_1t_1 - 1) | \setof{ d [ s_1t_1 (d_s s_1 - 4) + d_t t_1 ]
     + (2-2g) s_1 t_1 }.$%
$The \RHS\ is divisible by $t_1$, and $\gcd(s_1t_1 - 1, t_1) = 1$, so:
$$(s_1t_1 - 1) | \setof{ d [ s_1 (d_s s_1 - 4) + d_t ] + (2-2g) s_1}.\eqno(*)$%
$Thus for some $k \in \N$, we have:
$$(s_1t_1 - 1)k = d [ s_1 (d_s s_1 - 4) + d_t ] + (2-2g) s_1.$%
$Reorganizing, we find:
$$(s_1t_1 - 1)k = (d d_s)s_1^2 + (2-2g-4d)s_1 + d d_t.  \eqno(**)$%
$Now we have $t \geq s$, so $t_1 \geq (d_s/d_t)s_1$.  Hence:
$$\left[s_1^2\left({d_s \over d_t}\right) - 1\right]k
     \leq (d d_s)s_1^2 + (2-2g-4d)s_1 + d d_t.$%
$It is conceivable that the \LHS\ of this inequality is negative.  This will
not effect the following argument.
Suppose that $k \geq d d_t$.  After a short calculation, one finds that
$s_1 \leq d d_t / (2d+g-1)$, and hence that $s \leq d^2/(2d+g-1)$.  This
implies that $s < 2d^2$.  Hence, in order to prove our assertion that
$s < 2d^2$, we may assume that $k < d d_t$.

\par\noindent From $(**)$ we obtain:
$$(d d_s)s_1^2 + (2 - 2g - 4d - t_1 k)s_1 + (k + d d_t) = 0.$%
$Hence $s_1 | (k + d d_t)$.  Hence $s_1 \leq k + d d_t$.  Hence
$s_1 < 2d d_t$.  Hence $s < 2d^2$.

To complete the proof, we must show that $t < 2d^4$.
The \RHS\ of $(*)$ is nonzero, so:
$$s_1 t_1 - 1 \leq d[s_1(d_s s_1 - 4) + d_t] + (2 - 2g)s_1.$%
$Dividing by $s_1$ and isolating $t_1$, we find:
$$t_1 \leq d [ d_s s_1 - 4 + d_t s_1^{-1} ] + 2 - 2g + s_1^{-1}.$%
$Taking account of $t_1 = t d_t^{-1}$ and $s_1 = s d_s^{-1}$, we obtain:
$$t \leq d_t \setof{ d [ s - 4 + ds^{-1} ] + 2 - 2g } + ds^{-1}.$%
$Since $s < 2d^2$, it follows (with a little work) that $t < 2d^4$.  \qed
\end{proofnodot}

\begin{remark}
We give here an alternate proof of the main ingredient of the proof of (X),
namely that
$$(n-1) \kern3pt | \kern3pt \setof{ d [ n (s-4) + t ] + (2-2g)n }.
  \eqno(\dag)$%
$Let $\tC$ be the scheme-theoretic complete intersection of $S$ and $T$.
Let $\shJ$ be the ideal sheaf of $C$ in $\tC$.  Let $p$ be a closed point
of $C$.  If $\O_{S,p}$ is regular, then near $p$, $\tC$ and $C$ are Cartier
divisors on $S$, with $\tC = nC$.  Choose an isomorphism
$\O_{S,p} \iso k[[x,y]]$, such that $C$ corresponds to $V(x)$.  Then $\tC$
corresponds to $V(x^n)$.  Therefore the algebra of conormal invariants
$$\shA\ =\ \O_{\tC}/\shJ \o+ \shJ/\shJ^2 \o+ \shJ^2/\shJ^3 \o+ \cdots$%
$is a locally free $\O_C$-module near $p$.  But we similarly get the same
conclusion if $\O_{T,p}$ is regular, so in fact $\shA$ is locally free since
$\Sing(S) \cap \Sing(T) = \emptyset$.

Moreover, there is a line bundle $\shL$ on $C$ such that
$\shA \iso \O_C \o+ \shL \o+ \shL^2 \manyo+ \shL^{n-1}$.  Hence
$\chi(\shA) = n(1-g) + {n \choose 2}\deg(\shL)$.  On the other hand,
$\chi(\shA) = \chi(\O_{\tC})$, which (via $\tC = S \cap T$) is easily
computed to be $st(4-s-t)/2$.  Hence
$$n(1-g) + {n \choose 2}\deg(\shL) = {st(4-s-t) \over 2}.$%
$Hence
$${n \choose 2} \kern3pt \left| \kern3pt {st(4-s-t) \over 2} - n(1-g).\right.$%
$It is not difficult to verify that this is equivalent to $(\dag)$.
\end{remark}

\vspace{0.25in}

\section*{References}
\addtocontents{toc}{\protect\vspace*{2.25em}}
\addcontentsline{toc}{special}{References}
\ \par\noindent\vspace*{-0.25in}
\hfuzz 5pt

\bgroup\Resetstrings%
\def\Loccittest{}\def\Abbtest{ }\def\Capssmallcapstest{}\def\Edabbtest{
}\def\Edcapsmallcapstest{}\def\Underlinetest{ }%
\def\NoArev{1}\def\NoErev{0}\def\Acnt{2}\def\Ecnt{0}\def\acnt{0}\def\ecnt{0}%
\def\Ftest{ }\def\Fstr{1}%
\def\Atest{ }\def\Astr{Bott\Revcomma R\Initper %
  \Aand L\Initper \Initgap W\Initper  Tu}%
\def\Ttest{ }\def\Tstr{Differential Forms in Algebraic Topology}%
\def\Itest{ }\def\Istr{Spring\-er-Ver\-lag}%
\def\Ctest{ }\def\Cstr{New York}%
\def\Dtest{ }\def\Dstr{1982}%
\def\Qtest{ }\def\Qstr{access via "bott tu"}%
\Refformat\egroup%

\bgroup\Resetstrings%
\def\Loccittest{}\def\Abbtest{ }\def\Capssmallcapstest{}\def\Edabbtest{
}\def\Edcapsmallcapstest{}\def\Underlinetest{ }%
\def\NoArev{1}\def\NoErev{0}\def\Acnt{1}\def\Ecnt{0}\def\acnt{0}\def\ecnt{0}%
\def\Ftest{ }\def\Fstr{2}%
\def\Atest{ }\def\Astr{Brieskorn\Revcomma E\Initper }%
\def\Ttest{ }\def\Tstr{Die Aufl\"osung der rationalen Singularit\"aten
holomorpher Abbildungen}%
\def\Jtest{ }\def\Jstr{Math. Ann.}%
\def\Vtest{ }\def\Vstr{178}%
\def\Dtest{ }\def\Dstr{1968}%
\def\Ptest{ }\def\Pcnt{ }\def\Pstr{255--270}%
\def\Qtest{ }\def\Qstr{access via "brieskorn simultaneous resolution annalen"}%
\def\Xtest{ }\def\Xstr{Not on file.  The author studies simultaneous resolution
of rational singularities.}%
\Refformat\egroup%

\bgroup\Resetstrings%
\def\Loccittest{}\def\Abbtest{ }\def\Capssmallcapstest{}\def\Edabbtest{
}\def\Edcapsmallcapstest{}\def\Underlinetest{ }%
\def\NoArev{1}\def\NoErev{0}\def\Acnt{1}\def\Ecnt{0}\def\acnt{0}\def\ecnt{0}%
\def\Ftest{ }\def\Fstr{3}%
\def\Atest{ }\def\Astr{Catanese\Revcomma F\Initper }%
\def\Ttest{ }\def\Tstr{Babbage's conjecture, contact of surfaces, symmetric
determinantal varieties and applications}%
\def\Jtest{ }\def\Jstr{Invent. Math.}%
\def\Vtest{ }\def\Vstr{63}%
\def\Dtest{ }\def\Dstr{1981}%
\def\Ptest{ }\def\Pcnt{ }\def\Pstr{433--465}%
\def\Qtest{ }\def\Qstr{access via "catanese contact"}%
\Refformat\egroup%

\bgroup\Resetstrings%
\def\Loccittest{}\def\Abbtest{ }\def\Capssmallcapstest{}\def\Edabbtest{
}\def\Edcapsmallcapstest{}\def\Underlinetest{ }%
\def\NoArev{1}\def\NoErev{0}\def\Acnt{1}\def\Ecnt{0}\def\acnt{0}\def\ecnt{0}%
\def\Ftest{ }\def\Fstr{4}%
\def\Atest{ }\def\Astr{Deligne\Revcomma P\Initper }%
\def\Ttest{ }\def\Tstr{{\rm\tolerance=1000 Cohomologie des intersections
compl\`etes,  \EXPOSE\ XI in {\it\SGA} (SGA 7)}}%
\def\Stest{ }\def\Sstr{Lecture Notes in Mathematics}%
\def\Itest{ }\def\Istr{Spring\-er-Ver\-lag}%
\def\Ctest{ }\def\Cstr{New York}%
\def\Vtest{ }\def\Vstr{340}%
\def\Dtest{ }\def\Dstr{1973}%
\def\Ptest{ }\def\Pcnt{ }\def\Pstr{39--61}%
\def\Qtest{ }\def\Qstr{access via "deligne intersections"}%
\Refformat\egroup%

\bgroup\Resetstrings%
\def\Loccittest{}\def\Abbtest{ }\def\Capssmallcapstest{}\def\Edabbtest{
}\def\Edcapsmallcapstest{}\def\Underlinetest{ }%
\def\NoArev{1}\def\NoErev{0}\def\Acnt{2}\def\Ecnt{0}\def\acnt{0}\def\ecnt{0}%
\def\Ftest{ }\def\Fstr{5}%
\def\Ven{Van de Ven}{}%
\def\Atest{ }\def\Astr{Eisenbud\Revcomma D\Initper %
  \Aand A\Initper  \Ven}%
\def\Ttest{ }\def\Tstr{On the normal bundles of smooth rational space curves}%
\def\Jtest{ }\def\Jstr{Math. Ann.}%
\def\Vtest{ }\def\Vstr{256}%
\def\Dtest{ }\def\Dstr{1981}%
\def\Ptest{ }\def\Pcnt{ }\def\Pstr{453--463}%
\def\Qtest{ }\def\Qstr{access via "eisenbud normal bundles"}%
\def\Xtest{ }\def\Xstr{Not on file.}%
\Refformat\egroup%

\bgroup\Resetstrings%
\def\Loccittest{}\def\Abbtest{ }\def\Capssmallcapstest{}\def\Edabbtest{
}\def\Edcapsmallcapstest{}\def\Underlinetest{ }%
\def\NoArev{1}\def\NoErev{0}\def\Acnt{1}\def\Ecnt{0}\def\acnt{0}\def\ecnt{0}%
\def\Ftest{ }\def\Fstr{6}%
\def\Atest{ }\def\Astr{Fossum\Revcomma R\Initper \Initgap M\Initper }%
\def\Ttest{ }\def\Tstr{The Divisor Class Group of a Krull Domain}%
\def\Itest{ }\def\Istr{Spring\-er-Ver\-lag}%
\def\Ctest{ }\def\Cstr{New York}%
\def\Dtest{ }\def\Dstr{1973}%
\def\Qtest{ }\def\Qstr{access via "fossum"}%
\def\Xtest{ }\def\Xstr{I don't have this.}%
\Refformat\egroup%

\bgroup\Resetstrings%
\def\Loccittest{}\def\Abbtest{ }\def\Capssmallcapstest{}\def\Edabbtest{
}\def\Edcapsmallcapstest{}\def\Underlinetest{ }%
\def\NoArev{1}\def\NoErev{0}\def\Acnt{1}\def\Ecnt{0}\def\acnt{0}\def\ecnt{0}%
\def\Ftest{ }\def\Fstr{7}%
\def\Atest{ }\def\Astr{Fulton\Revcomma W\Initper }%
\def\Ttest{ }\def\Tstr{Intersection Theory}%
\def\Itest{ }\def\Istr{Spring\-er-Ver\-lag}%
\def\Ctest{ }\def\Cstr{New York}%
\def\Dtest{ }\def\Dstr{1984}%
\def\Qtest{ }\def\Qstr{access via "fulton intersection theory"}%
\Refformat\egroup%

\bgroup\Resetstrings%
\def\Loccittest{}\def\Abbtest{ }\def\Capssmallcapstest{}\def\Edabbtest{
}\def\Edcapsmallcapstest{}\def\Underlinetest{ }%
\def\NoArev{1}\def\NoErev{0}\def\Acnt{1}\def\Ecnt{0}\def\acnt{0}\def\ecnt{0}%
\def\Ftest{ }\def\Fstr{8}%
\def\Atest{ }\def\Astr{Gallarati\Revcomma D\Initper }%
\def\Ttest{ }\def\Tstr{Ricerche sul contatto di superfiche algebriche lungo
curve}%
\def\Jtest{ }\def\Jstr{Acad\'emie royale de Belgique, Classe des Sciences,
M\'emoires, Collection in-$8^0$}%
\def\Vtest{ }\def\Vstr{32}%
\def\Dtest{ }\def\Dstr{1960}%
\def\Ptest{ }\def\Pcnt{ }\def\Pstr{1--78}%
\def\Qtest{ }\def\Qstr{access via "gallarati belgique"}%
\Refformat\egroup%

\bgroup\Resetstrings%
\def\Loccittest{}\def\Abbtest{ }\def\Capssmallcapstest{}\def\Edabbtest{
}\def\Edcapsmallcapstest{}\def\Underlinetest{ }%
\def\NoArev{1}\def\NoErev{0}\def\Acnt{2}\def\Ecnt{0}\def\acnt{0}\def\ecnt{0}%
\def\Ftest{ }\def\Fstr{9}%
\def\Atest{ }\def\Astr{Griffiths\Revcomma P\Initper %
  \Aand J\Initper  Harris}%
\def\Ttest{ }\def\Tstr{Principles of Algebraic Geometry}%
\def\Itest{ }\def\Istr{John Wiley \& Sons}%
\def\Ctest{ }\def\Cstr{New York}%
\def\Dtest{ }\def\Dstr{1978}%
\def\Qtest{ }\def\Qstr{access via "griffiths harris principles" (was "griffiths
harris")}%
\Refformat\egroup%

\bgroup\Resetstrings%
\def\Loccittest{}\def\Abbtest{ }\def\Capssmallcapstest{}\def\Edabbtest{
}\def\Edcapsmallcapstest{}\def\Underlinetest{ }%
\def\NoArev{1}\def\NoErev{0}\def\Acnt{1}\def\Ecnt{0}\def\acnt{0}\def\ecnt{0}%
\def\Ftest{ }\def\Fstr{10}%
\def\Atest{ }\def\Astr{Grothendieck\Revcomma A\Initper }%
\def\Ttest{ }\def\Tstr{{\rm Cohomologie Locale des Faisceaux Coherents et
Theorems de Lefschetz Locaux et Globaux, in {\it \SGA} (SGA 2)}}%
\def\Itest{ }\def\Istr{North-Holland}%
\def\Dtest{ }\def\Dstr{1968}%
\def\Qtest{ }\def\Qstr{access via "grothendieck lefschetz"}%
\Refformat\egroup%

\bgroup\Resetstrings%
\def\Loccittest{}\def\Abbtest{ }\def\Capssmallcapstest{}\def\Edabbtest{
}\def\Edcapsmallcapstest{}\def\Underlinetest{ }%
\def\NoArev{1}\def\NoErev{0}\def\Acnt{2}\def\Ecnt{0}\def\acnt{0}\def\ecnt{0}%
\def\Ftest{ }\def\Fstr{11}%
\def\Atest{ }\def\Astr{Grothendieck\Revcomma A\Initper %
  \Aand J\Initper \Initgap A\Initper  Dieudonn\'e}%
\def\Ttest{ }\def\Tstr{El\'ements de \geom\ alg\'e\-brique II}%
\def\Jtest{ }\def\Jstr{Inst. Hautes \'Etudes Sci. Publ. Math.}%
\def\Vtest{ }\def\Vstr{8}%
\def\Dtest{ }\def\Dstr{1961}%
\def\Qtest{ }\def\Qstr{access via "EGA2"}%
\Refformat\egroup%

\bgroup\Resetstrings%
\def\Loccittest{}\def\Abbtest{ }\def\Capssmallcapstest{}\def\Edabbtest{
}\def\Edcapsmallcapstest{}\def\Underlinetest{ }%
\def\NoArev{1}\def\NoErev{0}\def\Acnt{2}\def\Ecnt{0}\def\acnt{0}\def\ecnt{0}%
\def\Ftest{ }\def\Fstr{12}%
\def\Atest{ }\def\Astr{Grothendieck\Revcomma A\Initper %
  \Aand J\Initper \Initgap A\Initper  Dieudonn\'e}%
\def\Ttest{ }\def\Tstr{El\'ements de \geom\ alg\'e\-brique IV (part four)}%
\def\Jtest{ }\def\Jstr{Inst. Hautes \'Etudes Sci. Publ. Math.}%
\def\Vtest{ }\def\Vstr{32}%
\def\Dtest{ }\def\Dstr{1967}%
\def\Qtest{ }\def\Qstr{access via "EGA4-4"}%
\def\Astr{\Underlinemark}%
\Refformat\egroup%

\bgroup\Resetstrings%
\def\Loccittest{}\def\Abbtest{ }\def\Capssmallcapstest{}\def\Edabbtest{
}\def\Edcapsmallcapstest{}\def\Underlinetest{ }%
\def\NoArev{1}\def\NoErev{0}\def\Acnt{1}\def\Ecnt{0}\def\acnt{0}\def\ecnt{0}%
\def\Ftest{ }\def\Fstr{13}%
\def\Atest{ }\def\Astr{Hartshorne\Revcomma R\Initper }%
\def\Ttest{ }\def\Tstr{Ample Subvarieties of Algebraic Varieties}%
\def\Stest{ }\def\Sstr{Lecture \snap Notes in Mathematics}%
\def\Itest{ }\def\Istr{Spring\-er-Ver\-lag}%
\def\Ctest{ }\def\Cstr{New York}%
\def\Vtest{ }\def\Vstr{156}%
\def\Dtest{ }\def\Dstr{1970}%
\def\Qtest{ }\def\Qstr{access via "hartshorne ample subvarieties"}%
\Refformat\egroup%

\bgroup\Resetstrings%
\def\Loccittest{}\def\Abbtest{ }\def\Capssmallcapstest{}\def\Edabbtest{
}\def\Edcapsmallcapstest{}\def\Underlinetest{ }%
\def\NoArev{1}\def\NoErev{0}\def\Acnt{1}\def\Ecnt{0}\def\acnt{0}\def\ecnt{0}%
\def\Ftest{ }\def\Fstr{14}%
\def\Atest{ }\def\Astr{Hartshorne\Revcomma R\Initper }%
\def\Ttest{ }\def\Tstr{Algebraic Geometry}%
\def\Itest{ }\def\Istr{Spring\-er-Ver\-lag}%
\def\Ctest{ }\def\Cstr{New York}%
\def\Dtest{ }\def\Dstr{1977}%
\def\Qtest{ }\def\Qstr{access via "hartshorne algebraic geometry"}%
\def\Astr{\Underlinemark}%
\Refformat\egroup%

\bgroup\Resetstrings%
\def\Loccittest{}\def\Abbtest{ }\def\Capssmallcapstest{}\def\Edabbtest{
}\def\Edcapsmallcapstest{}\def\Underlinetest{ }%
\def\NoArev{1}\def\NoErev{0}\def\Acnt{1}\def\Ecnt{0}\def\acnt{0}\def\ecnt{0}%
\def\Ftest{ }\def\Fstr{15}%
\def\Atest{ }\def\Astr{Hartshorne\Revcomma R\Initper }%
\def\Ttest{ }\def\Tstr{Complete intersections in characteristic $p > 0$}%
\def\Jtest{ }\def\Jstr{Amer. J. Math.}%
\def\Vtest{ }\def\Vstr{101}%
\def\Dtest{ }\def\Dstr{1979}%
\def\Ptest{ }\def\Pcnt{ }\def\Pstr{380--383}%
\def\Qtest{ }\def\Qstr{access via "hartshorne complete intersections
characteristic"}%
\def\Xtest{ }\def\Xstr{had: (Based on a talk given at the 1964 Woods Hole
conference)}%
\def\Astr{\Underlinemark}%
\Refformat\egroup%

\bgroup\Resetstrings%
\def\Loccittest{}\def\Abbtest{ }\def\Capssmallcapstest{}\def\Edabbtest{
}\def\Edcapsmallcapstest{}\def\Underlinetest{ }%
\def\NoArev{1}\def\NoErev{0}\def\Acnt{1}\def\Ecnt{0}\def\acnt{0}\def\ecnt{0}%
\def\Ftest{ }\def\Fstr{16}%
\def\Atest{ }\def\Astr{Iversen\Revcomma B\Initper }%
\def\Ttest{ }\def\Tstr{Generic Local Structure in Commutative Algebra}%
\def\Stest{ }\def\Sstr{Lecture Notes in Mathematics}%
\def\Itest{ }\def\Istr{Spring\-er-Ver\-lag}%
\def\Ctest{ }\def\Cstr{New York}%
\def\Vtest{ }\def\Vstr{310}%
\def\Dtest{ }\def\Dstr{1973}%
\def\Qtest{ }\def\Qstr{access via "iversen generic local structure"}%
\Refformat\egroup%

\bgroup\Resetstrings%
\def\Loccittest{}\def\Abbtest{ }\def\Capssmallcapstest{}\def\Edabbtest{
}\def\Edcapsmallcapstest{}\def\Underlinetest{ }%
\def\NoArev{1}\def\NoErev{0}\def\Acnt{1}\def\Ecnt{0}\def\acnt{0}\def\ecnt{0}%
\def\Ftest{ }\def\Fstr{17}%
\def\Atest{ }\def\Astr{Jaffe\Revcomma D\Initper \Initgap B\Initper }%
\def\Ttest{ }\def\Tstr{Space curves which are the intersection of a cone with
another surface}%
\def\Jtest{ }\def\Jstr{Duke Math. J.}%
\def\Vtest{ }\def\Vstr{57}%
\def\Dtest{ }\def\Dstr{1988}%
\def\Ptest{ }\def\Pcnt{ }\def\Pstr{859--876}%
\def\Qtest{ }\def\Qstr{access via "jaffe another"}%
\def\Htest{ }\def\Hstr{1}%
\Refformat\egroup%

\bgroup\Resetstrings%
\def\Loccittest{}\def\Abbtest{ }\def\Capssmallcapstest{}\def\Edabbtest{
}\def\Edcapsmallcapstest{}\def\Underlinetest{ }%
\def\NoArev{1}\def\NoErev{0}\def\Acnt{1}\def\Ecnt{0}\def\acnt{0}\def\ecnt{0}%
\def\Ftest{ }\def\Fstr{18}%
\def\Atest{ }\def\Astr{Jaffe\Revcomma D\Initper \Initgap B\Initper }%
\def\Ttest{ }\def\Tstr{On set theoretic complete intersections in $\P3$}%
\def\Jtest{ }\def\Jstr{Math. Ann.}%
\def\Vtest{ }\def\Vstr{285}%
\def\Dtest{ }\def\Dstr{1989}%
\def\Ptest{ }\def\Pcnt{ }\def\Pstr{165--173, 175, 174, 176}%
\def\Qtest{ }\def\Qstr{access via "jaffe on set theoretic annalen"}%
\def\Htest{ }\def\Hstr{2}%
\def\Astr{\Underlinemark}%
\Refformat\egroup%

\bgroup\Resetstrings%
\def\Loccittest{}\def\Abbtest{ }\def\Capssmallcapstest{}\def\Edabbtest{
}\def\Edcapsmallcapstest{}\def\Underlinetest{ }%
\def\NoArev{1}\def\NoErev{0}\def\Acnt{1}\def\Ecnt{0}\def\acnt{0}\def\ecnt{0}%
\def\Ftest{ }\def\Fstr{19}%
\def\Atest{ }\def\Astr{Jaffe\Revcomma D\Initper \Initgap B\Initper }%
\def\Ttest{ }\def\Tstr{Smooth curves on a cone which pass through its vertex}%
\def\Jtest{ }\def\Jstr{Manu\-scripta Math.}%
\def\Vtest{ }\def\Vstr{73}%
\def\Dtest{ }\def\Dstr{1991}%
\def\Ptest{ }\def\Pcnt{ }\def\Pstr{187--205}%
\def\Qtest{ }\def\Qstr{access via "jaffe vertex"}%
\def\Htest{ }\def\Hstr{5}%
\def\Astr{\Underlinemark}%
\Refformat\egroup%

\bgroup\Resetstrings%
\def\Loccittest{}\def\Abbtest{ }\def\Capssmallcapstest{}\def\Edabbtest{
}\def\Edcapsmallcapstest{}\def\Underlinetest{ }%
\def\NoArev{1}\def\NoErev{0}\def\Acnt{1}\def\Ecnt{0}\def\acnt{0}\def\ecnt{0}%
\def\Ftest{ }\def\Fstr{20}%
\def\Atest{ }\def\Astr{Jaffe\Revcomma D\Initper \Initgap B\Initper }%
\def\Ttest{ }\def\Tstr{Local geometry of smooth curves passing through rational
double points}%
\def\Jtest{ }\def\Jstr{Math. Ann.}%
\def\Vtest{ }\def\Vstr{294}%
\def\Dtest{ }\def\Dstr{1992}%
\def\Ptest{ }\def\Pcnt{ }\def\Pstr{645--660}%
\def\Qtest{ }\def\Qstr{access via "jaffe local geometry"}%
\def\Htest{ }\def\Hstr{6}%
\def\Astr{\Underlinemark}%
\Refformat\egroup%

\bgroup\Resetstrings%
\def\Loccittest{}\def\Abbtest{ }\def\Capssmallcapstest{}\def\Edabbtest{
}\def\Edcapsmallcapstest{}\def\Underlinetest{ }%
\def\NoArev{1}\def\NoErev{0}\def\Acnt{1}\def\Ecnt{0}\def\acnt{0}\def\ecnt{0}%
\def\Ftest{ }\def\Fstr{21}%
\def\Atest{ }\def\Astr{Lang\Revcomma W\Initper \Initgap E\Initper }%
\def\Ttest{ }\def\Tstr{Remarks on $p$-torsion of algebraic surfaces}%
\def\Jtest{ }\def\Jstr{Compositio Math.}%
\def\Vtest{ }\def\Vstr{52}%
\def\Dtest{ }\def\Dstr{1984}%
\def\Ptest{ }\def\Pcnt{ }\def\Pstr{197--202}%
\def\Qtest{ }\def\Qstr{access via "william lang torsion"}%
\Refformat\egroup%

\bgroup\Resetstrings%
\def\Loccittest{}\def\Abbtest{ }\def\Capssmallcapstest{}\def\Edabbtest{
}\def\Edcapsmallcapstest{}\def\Underlinetest{ }%
\def\NoArev{1}\def\NoErev{0}\def\Acnt{1}\def\Ecnt{0}\def\acnt{0}\def\ecnt{0}%
\def\Ftest{ }\def\Fstr{22}%
\def\Atest{ }\def\Astr{Lipman\Revcomma J\Initper }%
\def\Ttest{ }\def\Tstr{Rational singularities, with applications to algebraic
surfaces and unique factorization}%
\def\Jtest{ }\def\Jstr{Inst. Hautes \'Etudes Sci. Publ. Math.}%
\def\Vtest{ }\def\Vstr{36}%
\def\Dtest{ }\def\Dstr{1969}%
\def\Ptest{ }\def\Pcnt{ }\def\Pstr{195--279}%
\def\Qtest{ }\def\Qstr{access via "lipman rational"}%
\Refformat\egroup%

\bgroup\Resetstrings%
\def\Loccittest{}\def\Abbtest{ }\def\Capssmallcapstest{}\def\Edabbtest{
}\def\Edcapsmallcapstest{}\def\Underlinetest{ }%
\def\NoArev{1}\def\NoErev{0}\def\Acnt{1}\def\Ecnt{0}\def\acnt{0}\def\ecnt{0}%
\def\Ftest{ }\def\Fstr{23}%
\def\Atest{ }\def\Astr{Matsumura\Revcomma H\Initper }%
\def\Ttest{ }\def\Tstr{Commutative Algebra, second edition}%
\def\Itest{ }\def\Istr{Benjamin/Cum\-mings}%
\def\Ctest{ }\def\Cstr{Reading, Massachusetts}%
\def\Dtest{ }\def\Dstr{1980}%
\def\Qtest{ }\def\Qstr{access via "matsumura commutative algebra"}%
\Refformat\egroup%

\bgroup\Resetstrings%
\def\Loccittest{}\def\Abbtest{ }\def\Capssmallcapstest{}\def\Edabbtest{
}\def\Edcapsmallcapstest{}\def\Underlinetest{ }%
\def\NoArev{1}\def\NoErev{0}\def\Acnt{1}\def\Ecnt{0}\def\acnt{0}\def\ecnt{0}%
\def\Ftest{ }\def\Fstr{24}%
\def\Atest{ }\def\Astr{Miyaoka\Revcomma Y\Initper }%
\def\Ttest{ }\def\Tstr{The maximal number of quotient singularities on surfaces
with given numerical invariants}%
\def\Jtest{ }\def\Jstr{Math. Ann.}%
\def\Vtest{ }\def\Vstr{268}%
\def\Dtest{ }\def\Dstr{1984}%
\def\Ptest{ }\def\Pcnt{ }\def\Pstr{159--171}%
\def\Qtest{ }\def\Qstr{access via "miyaoka quotient"}%
\Refformat\egroup%

\bgroup\Resetstrings%
\def\Loccittest{}\def\Abbtest{ }\def\Capssmallcapstest{}\def\Edabbtest{
}\def\Edcapsmallcapstest{}\def\Underlinetest{ }%
\def\NoArev{1}\def\NoErev{0}\def\Acnt{1}\def\Ecnt{0}\def\acnt{0}\def\ecnt{0}%
\def\Ftest{ }\def\Fstr{25}%
\def\Atest{ }\def\Astr{Morrison\Revcomma D\Initper \Initgap R\Initper }%
\def\Ttest{ }\def\Tstr{The birational geometry of surfaces with rational double
points}%
\def\Jtest{ }\def\Jstr{Math. Ann.}%
\def\Vtest{ }\def\Vstr{271}%
\def\Dtest{ }\def\Dstr{1985}%
\def\Ptest{ }\def\Pcnt{ }\def\Pstr{415--438}%
\def\Qtest{ }\def\Qstr{access via "morrison rational double points"}%
\def\Xtest{ }\def\Xstr{MR 87e:14011. \par The author considers the problem of
when a curve on a surface can be contracted to a surface having only rational
double points.  Has appendix on rational double points.}%
\Refformat\egroup%

\bgroup\Resetstrings%
\def\Loccittest{}\def\Abbtest{ }\def\Capssmallcapstest{}\def\Edabbtest{
}\def\Edcapsmallcapstest{}\def\Underlinetest{ }%
\def\NoArev{1}\def\NoErev{0}\def\Acnt{1}\def\Ecnt{0}\def\acnt{0}\def\ecnt{0}%
\def\Ftest{ }\def\Fstr{26}%
\def\Atest{ }\def\Astr{Mumford\Revcomma D\Initper }%
\def\Ttest{ }\def\Tstr{The topology of normal singularities of an algebraic
surface and a criterion for simplicity}%
\def\Jtest{ }\def\Jstr{Inst. Hautes \'Etudes Sci. Publ. Math.}%
\def\Vtest{ }\def\Vstr{9}%
\def\Dtest{ }\def\Dstr{1961}%
\def\Ptest{ }\def\Pcnt{ }\def\Pstr{229--246}%
\def\Qtest{ }\def\Qstr{access via "mumford topology of normal"}%
\Refformat\egroup%

\bgroup\Resetstrings%
\def\Loccittest{}\def\Abbtest{ }\def\Capssmallcapstest{}\def\Edabbtest{
}\def\Edcapsmallcapstest{}\def\Underlinetest{ }%
\def\NoArev{1}\def\NoErev{0}\def\Acnt{1}\def\Ecnt{0}\def\acnt{0}\def\ecnt{0}%
\def\Ftest{ }\def\Fstr{27}%
\def\Atest{ }\def\Astr{Rao\Revcomma A\Initper \Initgap P\Initper }%
\def\Ttest{ }\def\Tstr{On self-linked curves}%
\def\Jtest{ }\def\Jstr{Duke Math. J.}%
\def\Vtest{ }\def\Vstr{49}%
\def\Dtest{ }\def\Dstr{1982}%
\def\Ptest{ }\def\Pcnt{ }\def\Pstr{251--273}%
\def\Qtest{ }\def\Qstr{access via "rao self linked"}%
\Refformat\egroup%

\bgroup\Resetstrings%
\def\Loccittest{}\def\Abbtest{ }\def\Capssmallcapstest{}\def\Edabbtest{
}\def\Edcapsmallcapstest{}\def\Underlinetest{ }%
\def\NoArev{1}\def\NoErev{0}\def\Acnt{1}\def\Ecnt{0}\def\acnt{0}\def\ecnt{0}%
\def\Ftest{ }\def\Fstr{28}%
\def\Atest{ }\def\Astr{Robbiano\Revcomma L\Initper }%
\def\Ttest{ }\def\Tstr{A problem of complete intersections}%
\def\Jtest{ }\def\Jstr{Nagoya Math. J.}%
\def\Vtest{ }\def\Vstr{52}%
\def\Dtest{ }\def\Dstr{1973}%
\def\Ptest{ }\def\Pcnt{ }\def\Pstr{129--132}%
\def\Qtest{ }\def\Qstr{access via "robbiano complete intersections nagoya"}%
\def\Xtest{ }\def\Xstr{Not on file.}%
\Refformat\egroup%

\bgroup\Resetstrings%
\def\Loccittest{}\def\Abbtest{ }\def\Capssmallcapstest{}\def\Edabbtest{
}\def\Edcapsmallcapstest{}\def\Underlinetest{ }%
\def\NoArev{1}\def\NoErev{0}\def\Acnt{1}\def\Ecnt{0}\def\acnt{0}\def\ecnt{0}%
\def\Ftest{ }\def\Fstr{29}%
\def\Atest{ }\def\Astr{Sakai\Revcomma F\Initper }%
\def\Ttest{ }\def\Tstr{Weil divisors on normal surfaces}%
\def\Jtest{ }\def\Jstr{Duke Math. J.}%
\def\Vtest{ }\def\Vstr{51}%
\def\Dtest{ }\def\Dstr{1984}%
\def\Ptest{ }\def\Pcnt{ }\def\Pstr{877--887}%
\def\Qtest{ }\def\Qstr{access via "sakai weil divisors"}%
\Refformat\egroup%

\bgroup\Resetstrings%
\def\Loccittest{}\def\Abbtest{ }\def\Capssmallcapstest{}\def\Edabbtest{
}\def\Edcapsmallcapstest{}\def\Underlinetest{ }%
\def\NoArev{1}\def\NoErev{0}\def\Acnt{1}\def\Ecnt{0}\def\acnt{0}\def\ecnt{0}%
\def\Ftest{ }\def\Fstr{30}%
\def\Atest{ }\def\Astr{Samuel\Revcomma M\Initper \Initgap P\Initper }%
\def\Ttest{ }\def\Tstr{Handwritten notes (one page)}%
\def\Dtest{ }\def\Dstr{July 1964}%
\def\Qtest{ }\def\Qstr{access via "samuel unpublished"}%
\Refformat\egroup%

\end{document}